\documentclass[12pt,letterpaper]{elsarticle}
	\pdfoutput=1
	
	\usepackage{graphics,epsfig,graphicx,color}
	\graphicspath{{/EPSF/}{../Figures/}{Figures/}}
	
	\usepackage{amssymb,latexsym}
	
	\usepackage{setspace}
	\usepackage{amsmath}
	
	\usepackage{mathrsfs,amsfonts}
	
	\usepackage{amsthm,amsxtra}
	
	\usepackage[font={small,it}]{caption}

	\usepackage[utf8]{inputenc}

	\usepackage{graphics,epsfig,graphicx,color}
	\usepackage{subcaption}
	
	\newtheorem{remark}{Remark}[section]

	\newcommand{\beqa}{\begin{eqnarray}}
	\newcommand{\eeqa}[1]{\label{#1}\end{eqnarray}}
	\newcommand{\beq}{\begin{equation}}
	\newcommand{\eeq}[1]{\label{#1}\end{equation}}

	\newcommand{\GD}{\Delta}

	\def\Be{{\bf e}}

	\def\Bn{{\bf n}}

	\def\Bx{{\bf x}}
	\def\By{{\bf y}}

	\def\B0{{\bf 0}}

	\def \RR {{\mathbb R}}
	
	\def \ba {\begin{array}}
	\def \ea {\end{array}}
	\newtheorem {Thm} {Theorem} [section]

	\newtheorem {Adef} [Thm] {Definition}
	
	\newtheorem {Arem} [Thm] {Remark}
	
	\newtheorem {Aexa} [Thm] {Example}
	
	\newtheorem {Anot} [Thm] {Notation}

	\def \refe #1.{(\ref{#1})}
	\def \reff #1.{figure~\ref{#1}}
	\def \refs #1.{section~\ref{#1}}
	\def \refss #1.{subsection~\ref{#1}}
	\def \refD #1.{Definition~\ref{#1}}
	\def \refT #1.{Theorem~\ref{#1}}
	\def \refL #1.{Lemma~\ref{#1}}
	\def \refC #1.{Corollary~\ref{#1}}
	\def \refP #1.{Proposition~\ref{#1}}
	\def \refR #1.{Remark~\ref{#1}}
	\def \refE #1.{Example~\ref{#1}}
	\def \refN #1.{Notation~\ref{#1}}
	\newif\ifPDF
	\ifx\pdfoutput\undefined
	    \PDFfalse
	\else
	    \ifnum\pdfoutput > 0
	        \PDFtrue
	    \else
	        \PDFfalse
	    \fi
	\fi
	
	\ifPDF
	    \usepackage{pdftricks}
	    \begin{psinputs}
	        \usepackage{pstricks}
	        \usepackage{pstcol}
	        \usepackage{pst-plot}
	        \usepackage{pst-tree}
	        \usepackage{pst-eps}
	        \usepackage{multido}
	        \usepackage{pst-node}
	        \usepackage{pst-eps}
	    \end{psinputs}
	\else
	    \usepackage{pstricks}
	\fi
	
	\ifPDF
	    \usepackage[debug,pdftex,colorlinks=true,
	    linkcolor=blue, bookmarksopen=false,
	    plainpages=false,pdfpagelabels]{hyperref}
	\else
	    \usepackage[dvips]{hyperref}
	\fi
	
	\usepackage{fancyhdr}
	
\usepackage{calc}

\title{Sensitivity analysis for the active manipulation of Helmholtz fields in 3D}
\journal{refereed journal}
\author[negarguin]{Neil Jerome A. Egarguin}
\author[donofrei]{Daniel Onofrei}
\author[eplatt]{Eric Platt}

\address[negarguin]{University of Houston, Department of Mathematics and the University of the Philippines Los Ba\~nos, Institute of Mathematical Sciences and Physics: naegarguin1@up.edu.ph}
\address[donofrei]{University of Houston, Department of Mathematics: onofrei@math.uh.edu}
\address[eplatt]{University of Houston, Department of Mathematics: eplatt@math.uh.edu}

\newcommand{\figsizeB}{0.4750\textwidth}
\newcommand{\figsizeC}{0.425\textwidth} 

\begin{document}

\begin{frontmatter}


\begin{abstract}
	
	In this paper we continue to study the feasibility of active manipulation of Helmholtz fields and by using an improved and more robust numerical strategy we present a detailed sensitivity analysis for the method proposed in our previous works \cite{DOactivemanipulation},\cite{DO_EP}. In this regard, we study the behavior of physically relevant parameters (i.e. source power, control accuracy, stability) with respect to variations in: the type of control regions (bounded or unbounded), relative position of the control regions, distances between the control regions and the source, frequency range and fields to be approximated. We produce strong numerical evidence indicating the accuracy of our scheme and in the same time develop a better understanding of several important challenges for its physical implementation.

\end{abstract}


\end{frontmatter}

\section{Introduction and related results}

The active control of Helmholtz fields can be studied using various techniques such as those discussed in \cite{fuller}, \cite{peake}, \cite{newactr5}.  Common strategies used in applications include active noise control (see early works \cite{nelson}, \cite{olson} and \cite{ANC1}), sound field reproduction (\cite{kim}, \cite{rabenstein}, \cite{ahrens1} and monograph \cite{WFS1} with its references) and active acoustic cloaking (see \cite{miller}, \cite{AC1}, \cite{AC2}, \cite{AC3}, \cite{AC4}, \cite{broadband_vasquez},\cite{active_vasquez} and references therein).  In \cite{cheer}, a comprehensive comparison and analysis of the physical limitations of these approaches was done in the context of acoustic scattering.  

Methods of determining active controls for acoustic cloaking based on the Green representation theorem for the Helmholtz equation were proposed in \cite{vasquez} and \cite{miller} while a construction using generalized Calderon's potentials and boundary projection operators was proposed in \cite{loncaric2} (see also \cite{loncaric} where optimization of the sources of the active controls with respect to some quadratic functions of merit was performed).

In recent years many researches in the field focused their efforts on addressing the sound zone problem and its applications to personal audio systems. The goal of such efforts is to design inputs to a loud speaker array so that in a free space divided into several zones, a desired acoustic signature is synthesized in each zone, without affecting the quality of the sound produced in other zones. Among the approaches developed to solve the problem we can recall the use of the weighted pressure matching method presented in \cite{olivieri}, energy focusing, field cancellation, and synthesis approaches compared in \cite{coleman}, least squares minimization of acoustic control criteria discussed in \cite{cai} and modal-domain analysis proposed in \cite{zhang} (see also \cite{Cheerps}, \cite{Poletti}, for a review on strategies for the realization of personal audio zones).

Recently, building up and extending the earlier results presented in \cite{vasquez}, in \cite{DOactivemanipulation}, \cite{DO_EP} (see also \cite{Onofrei-S} for the quasistatics case), a novel strategy was developed and discussed for the problem of controlling radiating solutions of the Helmholtz equation in two and three dimensions. More explicitly, our effort is to characterize continuous secondary surface sources for the approximation of desired Helmholtz potentials in prescribed exterior regions of space (possibly including the far field region). The sources we construct , which are modeled in our analysis as continuous boundary data, could then, for their physical implementation, be discretized in monopole or dipole arrays ( see for example \cite{Doicu}).

In the same functional framework considered in  \cite{DOactivemanipulation} and \cite{DO_EP}, in \cite{hubenthal_DO}  the authors present a sensitivity study for the problem of characterizing surface sources with vanishingly small far fields and controllable near fields with numerical simulations in \cite{hubenthal_DO} performed only for the two dimensional case. 

In this paper, we present a detailed sensitivity study for the problem of three dimensional exterior control of scalar radiating Helmholtz fields explored in \cite{DOactivemanipulation}, \cite{DO_EP} by analyzing the behavior of relevant physical measures (e.g. source power, control accuracy, stability) with respect to variations in: the type of control regions (bounded or unbounded), relative position of the control regions, distances between the control regions and the exterior source, frequency and fields to be approximated. We note that, the problem considered can be also cast as an inverse source problem and this immediately reveals one of the major challenges we faced: the ill posed character such problems have which was highlighted also in our previous works.

The paper is organized as follows: Section \ref{theoretical} provides a discussion of relevant theoretical results obtained in \cite{DOactivemanipulation} and \cite{DO_EP} while Section \ref{numframework} introduces the numerical framework and optimization schemes used in the experiments. Section \ref{numresults} presents the results of the sensitivity analysis for two main configurations: 

\vspace{0.2cm}

$I.$ almost non-radiating sources with controllable near fields, i.e., radiators which approximate desired Helmholtz potentials in specified near field exterior regions while keeping a very low profile beyond a given fixed radius;

\vspace{0.2cm}
$II.$ sources approximating desired Helmholtz potentials in a collection of specified exterior compact disjoint regions of space. 

\vspace{0.2cm}
Sensitivity analysis results for configuration $I.$ above are presented in Subsection \ref{SSANR}. Sensitivity analysis results for configuration $II.$ above  are presented in Subsection \ref{SSnull} for a fixed frequency and in Subsection \ref{fourier} for a range of superimposed frequencies.  Section \ref{future} presents our conclusions and announces some of our future research reports on this subject.

\section{Introduction of the problem}
\label{theoretical}

In this section we present a general description of the scalar control problem we discuss describing the geometric and functional framework together with several essential theoretical results from \cite{DOactivemanipulation} and \cite{DO_EP}. 

Although our results will apply for any homogeneous isotropic medium and arbitrary number of exterior control regions, for illustrative purposes we consider a free space environment containing a source $D_a \subset \RR^3$ (single compact region) and only two control regions $D_1$ and $D_2$ that are mutually disjoint smooth domains. We require that the control domains and the source are mutually "well separated", i.e. $(D_1\cup D_2)\cap D_a=\emptyset$. Under these assumptions, in what follows we will focus on two geometric configurations, namely:
\begin{eqnarray}
&\;i)\;& D_1 \mbox{ bounded}, D_2 \mbox{ bounded},\label{geom}\\
&\;ii)& D_1 \mbox{ bounded}, D_2 = \RR^3\!\setminus\!B_R(\B0),\nonumber
\end{eqnarray}
where here and throughout the rest of the paper $B_R(\B0)$ denotes the ball centered in the origin and radius $R$. Also, by convention, throughout the paper we assume an $e^{-i\omega t}$ time dependence of the fields. 

Let $u_1$ and $u_2$ be the solutions of the Helmholtz equation in neighborhoods of $D_1$ and $D_2$, respectively. We will sometimes refer to such functions as Helmholtz potentials. Without loss of generality, assuming the source $D_a$ is mathematically modeled as a surface input (i.e., boundary data on $\partial D_a$) we focus on the case when one of the fields $u_1$ or $u_2$ is zero. The main goal is the characterization of boundary inputs (Neumann or Dirichlet data) on the surface of $D_a$ so that the radiating solution of the Helmholtz equation exterior to $D_a$ approximates $u_1$ in $D_1$ and $u_2$ in $D_2$. That is, assuming for exemplification that $u_2=0$, the main question is to characterize the Neumann data $v_n$ (or Dirichlet data $p_b$) on the boundary $\partial D_a$ of $D_a$ such that
\begin{equation}
\label{P1a}
\vspace{0.15cm}\left\{\vspace{0.15cm}\begin{array}{llll}
\GD u+k^2u=0 \mbox{ in }\RR^3\!\setminus \!\overline{D}_a\vspace{0.15cm},\\
\nabla u\cdot \Bn=v_n, (\mbox{ or } u=p_b)\mbox{ on }\partial D_a\vspace{0.15cm},\\
\displaystyle\left({\hat{\Bx}},\nabla u(\Bx)\right)\! -\!iku(\Bx)\!=\!o\left(\frac{1}{|\Bx|}\!\right)\!,\mbox{ as }|\Bx|\rightarrow\infty
\mbox{ uniformly for all ${\hat{\Bx}}$, }\end{array}\right.
\end{equation}
and \begin{equation} \label{P1b} u \approx u_1 \text{ in } D_1 \text{ and } u \approx u_2=0 \text{ in } D_2, \end{equation} in the sense of smooth norms, (e.g., $C^2$ norms), where here $\Bn$ denotes the outward normal to $\partial D_a$ and here and throughout the rest of the paper ${\hat{\Bx}}=\frac{\Bx}{|\Bx|}$ denotes the unit vector along the direction $\Bx$.

In \cite{DOactivemanipulation}, it was shown that under the assumptions and geometric configurations stated above, the problem \eqref{P1a} and \eqref{P1b} admits a solution except for a discrete family of  $k$ values.  For the rest of the paper we will assume that $k$ will be outside this discrete family. In this context, the existence of an infinite class of smooth functions $w$ characterizing the desired boundary inputs was established in \cite{DOactivemanipulation}. More explicitly, for any given smooth domain $D_{a'}$, with $D_{a'} \Subset D_a$, defining $\Bn_\By$ to be the outward normal to $\partial D_{a'}$ at $\By \in D_{a'} $, $\rho$ to be the density of the surrounding medium, $c$ to be the wave speed and $\Phi$ to be the free space fundamental solution of the Helmholtz equation, it was established in \cite{DOactivemanipulation} (see also \cite{DO_EP}) that there exists an infinity of smooth functions $w$ such that the Neumann data $v_n$ or the Dirichlet data $p_b$ required in problem (\ref{P1a}-\ref{P1b}) are given respectively by,
\begin{eqnarray}
& v_n(\Bx) =&\displaystyle \frac{-i}{\rho c k}\frac{\partial}{\partial\Bn_{\Bx}}\int_{\partial D_{a'}}w(y)\frac{\partial
	\Phi(\Bx,\By)}{\partial \Bn_{\By}}ds_{\By} \text{ and }\label{eqnvn}\\
\nonumber \\ 
& p_b(\Bx)=&\displaystyle \int_{\partial D_{a'}}w(y)\frac{\partial
	\Phi(\Bx,\By)}{\partial \Bn_{\By}}ds_{\By},\label{eqnpb}
\end{eqnarray}
for all $\Bx\in\partial D_a$. The utility of the introduction of the fictitious source region $D_{a'}$ is justified in the following two remarks.

\begin{remark}
	\label{rem1}
	We make the observation here that the real source on the boundary of which the desired inputs $v_n$ or $p_b$ are to be prescribed is $D_a$ while, as it can be observed, for their characterization in \eqref{eqnvn}, \eqref{eqnpb}  a ``fictitious" non-physical domain $D_{a'}$ was used. So while for simplicity and ease of computations it is assumed that $D_{a'}$ is smooth (actually spherical in the present paper) the physical source $D_a$ can have any shape as long as it is Lipschitz, compactly embeds $D_{a'}$, and $(D_1\cup D_2)\cap D_a=\emptyset$ as mentioned above.
\end{remark}

\begin{remark}
	\label{rem2}
	From (\ref{eqnvn}) and (\ref{eqnpb}), the smoothness of $D_{a'}$ together with the fact that $D_{a'} \Subset D_a$ (where here $\Subset$ denotes compact embedding) implies the smoothness of the boundary inputs $v_n$ and $p_b$ on $\partial D_a$. This makes the minimal assumption that $\partial D_a$ is Lipschitz sufficient for the exterior problem to be well-posed. 
\end{remark}

The next remark highlights the versatility of our strategy with regards to the type of field propagator employed. 
\begin{remark}
	\label{rem3}
	The boundary inputs given in (\ref{eqnvn}) and (\ref{eqnpb}) generate a solution $u$ of (\ref{P1a}-\ref{P1b}) as a double layer potential given by 
	
	\begin{equation}
	\label{remark3}
	u(\Bx)=\displaystyle \int_{\partial D_{a'}}w(y)\frac{\partial
		\Phi(\Bx,\By)}{\partial \Bn_{\By}}ds_{\By},\end{equation}
	for all $\Bx\in \RR^3\!\setminus \!{\overline D}_a$. On the other hand, as described in \cite{DO_EP}, our analysis permits the characterization of solutions $u$ of (\ref{P1a}-\ref{P1b}) that are complex linear  combinations of single and double layer potentials. This form is efficient when dealing with real wave numbers $k$.
\end{remark}

\section{Numerical framework and optimization scheme}
\label{numframework}

In this section, we briefly present the numerical framework and optimization scheme we employ in the study of the problems (\ref{P1a}-\ref{P1b}). We also describe the set-up and parameters used in the numerical simulations. The sensitivity analysis performed later in the paper is based on the optimization scheme described below.

In \cite{hubenthal_DO}, the $L^2$-optimization and sensitivity analysis for the solutions of the two-dimensional analog of (\ref{P1a}-\ref{P1b}) in the case of  geometry (\ref{geom}-{\it ii}) showed that a good approximation for a stable solution with minimal power budget is achieved when $D_1$ is very near to $D_a$. Meanwhile, the three-dimensional formulation was solved using the Tikhonov regularization with the Morozov discrepancy principle. Results from \cite{DOactivemanipulation} and \cite{hubenthal_DO} show that in order to find approximate smooth controls in $D_1$ and $D_2$, it suffices to have $L^2$ controls on the boundaries of slightly larger sets $W_1$ and $W_2$ with $D_1 \Subset W_1$, $D_2 \Subset W_2$, $W_1\cap W_2=\emptyset$ and $(W_1\cup W_2)\cap D_a=\emptyset$. In this regard, solution $u$ can be approximated by the following ansatz 
\begin{equation}
\mathcal{D}(w_\alpha)(\mathbf{x}) = \eta_{1}\int_{\partial D_{a'}} w_\alpha(\mathbf{y}) \frac{\partial \Phi(\mathbf{x}, \mathbf{y})}{\partial \Bn_{\mathbf{y}}}\,dS_{\mathbf{y}} + i\eta_{2} \int_{\partial D_{a'}} w_\alpha(\mathbf{y}) \Phi(\mathbf{x},\mathbf{y})\,dS_{\mathbf{y}}, \label{uapprox}
\end{equation}
where $\eta_{1}, \eta_{2} \in \mathbb{R}$ are fixed parameters indicating the weight assigned to the double and single layer potential terms. The function $w_\alpha$ is the Tikhonov regularization solution, i.e., the minimizer of the discrepancy functional,
\begin{equation}
F(w)=\frac{1}{\|w\|_{L^{2}(\partial D_{a})}^{2}}\| \mathcal{D}(w) - u_1\|_{L^{2}(\partial W_{1})}^{2} + \mu\|\mathcal{D}(w)\|_{L^{2}(\partial W_{2})}^{2} + \alpha\|w\|_{L^{2}(\partial D_{a'})}^{2}, \label{eq:tikhonovobjective}
\end{equation}
with the regularization parameter $\alpha$ (computed following the Morozov Discrepancy principle) representing the penalty weight for the power required by the solution and with the weight $\mu$ given by
\begin{equation}\mu=\left\{\begin{array}{lll}
\label{mu}
1, & \mbox{ if } D_2 \mbox{ is bounded }\\
\frac{1}{4\pi R^2}, & \mbox{ if } D_2 =\RR^3\setminus B_R(\B0)
\end{array}\right.,\end{equation}
where, here and further in the paper, $B_R(\B0)$ denotes the ball centered on the origin with a large enough radius $R$ such that $D_a\cup W_1\Subset B_R(\B0)$.  

Numerical simulations for the synthesis of the prescribed patterns on the regions $D_1$ and $D_2$ are performed using the spherical harmonics decomposition of $w_\alpha$ with $L$ harmonic orders, i.e.,
\begin{equation}
\label{density}
w_\alpha(\By)=\sum_{l=0}^{L}\sum_{p=-l}^{l}\alpha_{pl} Y_l^p({\hat{\By}}), \mbox{for }\By\in\partial D_{a'}
\end{equation}
for different prescribed values of $L$, where $Y_l^p$ above form the orthonormal family of spherical harmonics as considered in (\cite{Colton-Kress}, Chapter 2., \cite{stegun}).
For illustrative purposes we assumed in the numerics that the fictitious source , i.e., $D_{a'}$ above,  is a single sphere $B_{a'}(\B0)$ with $a'=0.01m$ discretized into 20, 000 points with 200 by 100 equidistant azimuthal and polar increments, respectively. Unless stated otherwise, for the initial geometries described in (\ref{geom}),  the near field or primary control region $D_1$ is defined as:
\begin{equation}
\label{D1}
D_1=\left \{(r,\theta, \phi) : r \in [0.011,0.015], \theta \in \left [-\frac{\pi}{4}, \frac{\pi}{4} \right ] , \phi \in \left [\frac{3\pi}{4}, \frac{5\pi}{4} \right] \right   \},
\end{equation} 
and, for configuration described at (\ref{geom}-{\it i}), $D_2$ is a secondary bounded control region defined by 
\begin{equation}
\label{D2_secondary}
D_2\!=\!\left \{\!(r,\theta, \phi) : r \!\in \![0.011,0.015], \theta \in \left [-\frac{\pi}{4}, \frac{\pi}{4} \right ] , \phi \in \left [\frac{7\pi}{4}, 2\pi\right]\!\cup\!\left[0,\frac{\pi}{4} \right]  \!\right \} \!+ \!(0.09,0,0).  
\end{equation} 
while for the configuration introduced at (\ref{geom}-{\it ii}) it is given by, \begin{equation} \label {D2_farfield} D_2= \mathbb R^3 \setminus B_{10}(\B0). \end{equation}

In \cite{DO_EP}, the method of moments together with the Tikhonov regularization procedure with Morozov discrepancy principle were used towards a solution for the above problems. Specialized Gauss - Legendre quadrature procedures were employed to numerically compute the moments represented by ${\cal{D}}(Y_l^p)$, where ${\cal{D}}$ is the integral operator propagator defined at \eqref{uapprox}. This strategy proved to be very sensitive near the source and cost inefficient when attempting to use more harmonic orders in the description of the density function $w_\alpha$ or perform a time domain simulation through related Fourier synthesis. 

For the current research effort we developed a more direct approach where the moments ${\cal{D}}(Y_l^p)$ (with ${\cal{D}}$ defined at \eqref{uapprox}) are evaluated explicitly by employing a truncated series obtained by making use of the addition theorem for the representation of the fundamental solution $\Phi(\mathbf{x}, \mathbf{y})$. More explicitly, from the addition theorem (\cite{Colton-Kress}, Chapter 2) we have,
\begin{equation}
\label{addition}
\Phi(\Bx,\By)=ik\sum_{n=0}^{\infty}\sum_{m=-n}^{n}h_n^{(1)}(k|\Bx|)Y_n^m({\hat{x}})j_n(k|\By|){\overline {Y_n^m({\hat{y}})}},
\end{equation}
where $h_n^{(1)}, j_n$ are the spherical Haenkel of first kind and spherical Bessel functions of order $n$ (\cite{Colton-Kress}, Chapter 2, see also \cite{watson} and \cite{stegun}) and where $Y_n^m$ form the orthonormal family as above in \eqref{density}. Then, considering the expansion for $w_\alpha$ defined in \eqref{density} together with \eqref{addition} and orthogonality of the spherical harmonics we obtain,
\begin{eqnarray}
& \displaystyle \int_{\partial D_{a'}}w_\alpha(\By)\frac{\partial
	\Phi(\Bx,\By)}{\partial \Bn_{\By}}ds_{\By}=& ik^2{a'}^2\sum_{l=0}^{L}\sum_{p=-l}^{l}\alpha_{pl}j'_lk(|\Bx|)h_l^{(1)}(k|\Bx|)Y_l^p({\hat{\Bx}}),\label{eqnddl}\\
\nonumber \\ 
& \displaystyle \int_{\partial D_{a'}}w_\alpha(\By)
	\Phi(\Bx,\By)ds_{\By}=& ik{a'}^2\sum_{l=0}^{L}\sum_{p=-l}^{l}\alpha_{pl}j_l(k|\Bx|)h_l^{(1)}(k|\Bx|)Y_l^p({\hat{\Bx}}).\label{eqndsl}
\end{eqnarray}

Expressions \eqref{eqnddl} and \eqref{eqndsl} are then used in the regularization routine for a much faster and more accurate computational tool to perform our current sensitivity analysis. 
\section{Numerical results and sensitivity analysis}
\label{numresults}

In this section, we present the sensitivity results obtained using the scheme discussed in the preceding sections. We consider two major geometries: 

\vspace{0.2cm}

{\textit i) Sources approximating a prescribed Helmholtz potential in $D_1$ defined at \eqref{D1} with very low field amplitudes in $D_2= \mathbb R^3 \setminus B_{10}(\B0)$}, (i.e., exterior of the ball of radius $10$m).

 \vspace{0.1cm}
 
 For the case when the source is geometrically modeled as a single compact region, we considered the problem \eqref{P1a} and \eqref{P1b} and addressed the sensitivity of the control results with respect to the size of $D_1$ and its distance from the fictitious source $D_{a'}$. 
 
 \vspace{0.2cm}
 
{\textit ii) Sources approximating a prescribed Helmholtz potential in $D_1$ defined at \eqref{D1} and another prescribed Helmholtz potential in $D_2$ defined at \eqref{D2_secondary}.}

\vspace{0.1cm}

   For this geometry, we assume the source modeled geometrically as a single compact region and for the problem \eqref{P1a} and \eqref{P1b} we perform sensitivity analysis of the control scheme 
   with respect to the relative position of the two control regions $D_1$ and $D_2$. We also study the sensitivity of our control scheme when a small amplitude field is generated in the near field region $D_1$ and a plane wave is approximated in region $D_2$ (i.e., problem \eqref{P1a} and \eqref{P1b} with $u_1=0$). We conclude the discussion for this geometry with the sensitivity analysis with respect to frequency and simulate, through a Fourier synthesis procedure, a time-domain source approximating an outgoing pulse in region $D_1$ and a null in region $D_2$.
	
\vspace{0.2cm}

Our sensitivity analysis will consider the behavior of the following physically relevant quantities with respect to various specified geometrical parameters:  $||w_\alpha||_{L^2(\partial D_{a'})}$ as an indicator of the overall power on the source, stability of the solution $w_\alpha$ (see definition \eqref{sensitivity_measure} below), $L^2$ and $L^\infty$ relative error (i.e., relative difference in the respective norms) in region where we approximate a plane wave and absolute $L^2$ and respectively $L^\infty$ error in the null region (error computed on $\partial D_2$ when $D_2= \mathbb R^3 \setminus B_{10}(\B0)$).
The stability measure is the $L^2$ relative norm of the difference of antenna density patterns:
\begin{equation} \label{sensitivity_measure} S(w, w_\epsilon)= \dfrac{\|w_\alpha-w_{\alpha_\epsilon}\|_{L^2(\partial D_{a'})}}{\|w_\alpha\|_{L^2(\partial D_{a'})}}, \end{equation}
where $w_{\alpha_\epsilon}$ is the solution of problem \eqref{P1a}, \eqref{P1b} with noise of magnitude $\epsilon$ in the data $u_1,u_2$. In this paper, our noise is modeled as a uniformly distributed additive perturbation of order $\epsilon$. 

Unless otherwise specified, each single-source sensitivity test was performed on two synthesis results: one with 15 and another with 30 harmonic orders used in \eqref{density}.

\subsection{Almost non-radiating source with controllable near field}
\label{SSANR}

\begin{figure}[!htbp] \centering
	\includegraphics[scale=0.4]{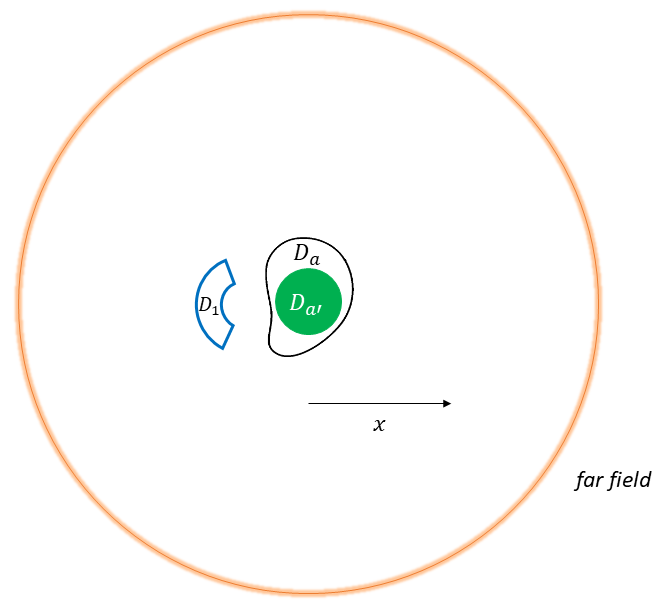}
	\caption{Initial geometry for the sensitivity experiments with an almost non-radiating single source.}
	\label{initial}
\end{figure}

In the following tests, we consider a single almost non-radiating source $D_a$ containing the following  fictitious region  $D_{a'}=B_{0.01}(\B0)$, i.e., a sphere centered at the origin and with radius $1cm$. We also consider the near control region $D_1$ as defined in \eqref{D1} where we match an {\textit{incoming}} plane wave $u_1(\Bx)= e^{i\Bx \cdot (10 \hat {\Be}_1)}$ propagating along the positive $x$-axis (i.e. towards the source $D_a$), with  $\hat {\Be}_1=\langle 1, 0, 0\rangle$ and the wave number $k=10$, and the far field region the exterior of a ball centered at the origin of radius $10m$, i.e., $D_2= \mathbb R^3 \setminus B_{10}(\B0)$ given in \eqref{D2_farfield}, where we maintain a near zero signature i.e., $u_2 (\Bx) \approx 0$. This initial geometry is shown in Figure \ref{initial}.

We recall that as stated above and proved in \cite{DOactivemanipulation} our method is such that only good boundary controls are needed for smooth interior controls. In this spirit, in Sections \ref{var1}, \ref{var2}, \ref{var3}, the boundary of region $D_1$ is uniformly discretized into 2,400 points while the far field boundary, i.e.,  $\partial D_2=\partial B_{10}(\B0)$ is uniformly  discretized with 20,000 points with 200 in the azimuthal and 100 in the polar direction respectively.

This geometry was considered in  \cite{DO_EP} where we obtained a density $w_\alpha$ on  $D_{a'}$ (and also necessary boundary inputs as a consequence of \eqref{eqnvn} and \eqref{eqnpb}) so that the relative control error in the region $D_1$ was of order $O(10^{-3})$ while maintaining field values of $O(10^{-3})$ in $D_2$. Recalling the ansatz for the propagating field described above at \eqref{remark3} in \cite{DO_EP} 
we produced numerical support showing the optimal density $w_\alpha$ in \eqref{remark3} exhibiting maximum values of order $O(10^5)$ and small and fast oscillations on the side facing $D_1$ with large and slower oscillation on the opposite side and at the poles. 

\begin{figure}[!htb] \centering
	\!\!\!\!\vspace{-1cm}\begin{subfigure}{\figsizeC}
		\includegraphics[width=1.4\textwidth]{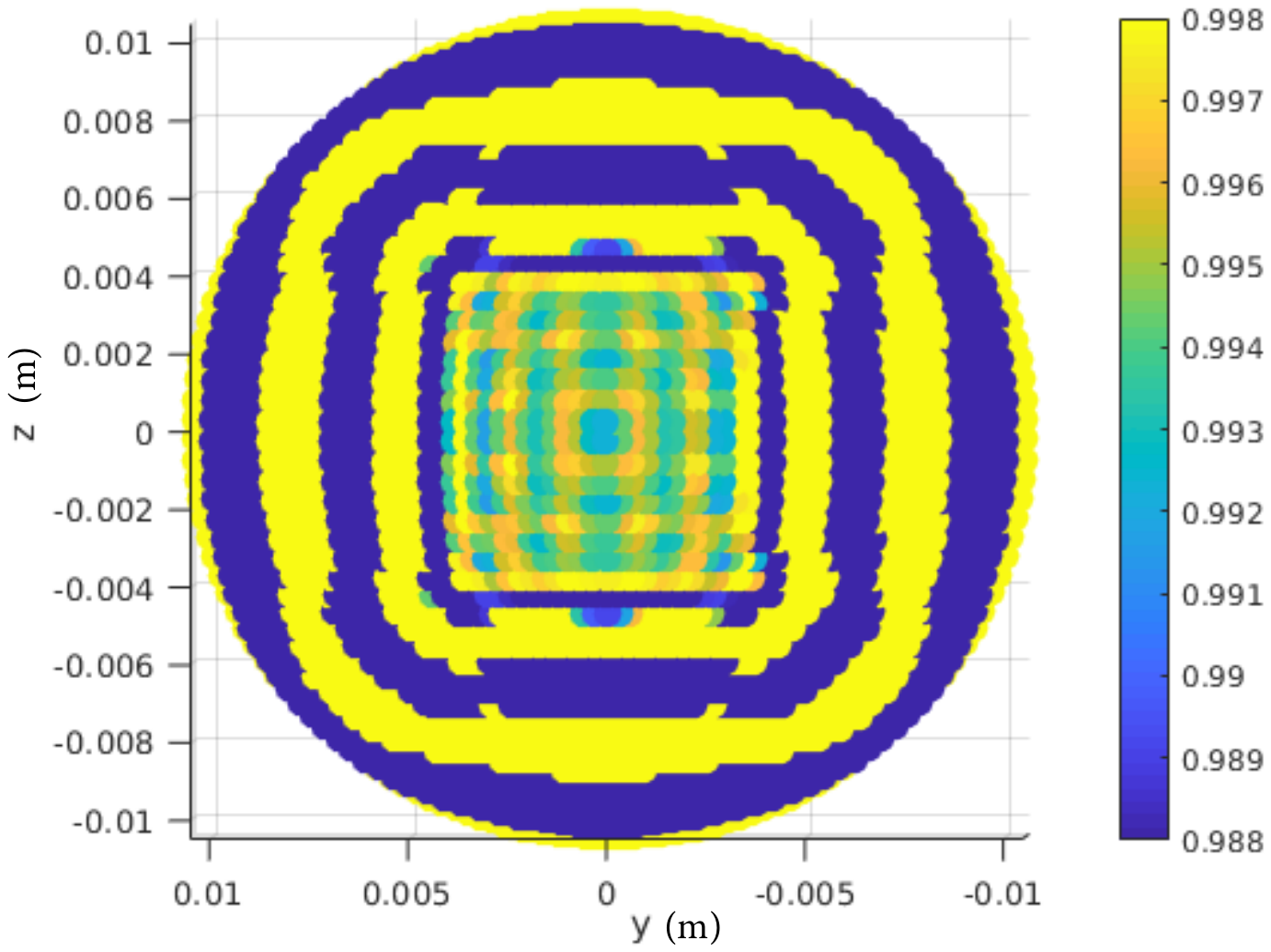}
		\label{kw_anr_front}
		\vspace{-4cm}
		\caption{front}
	\end{subfigure}
	\;\;\;\;\;\;\;\;\;\;\;\;\;\;\begin{subfigure}{\figsizeC}
		\includegraphics[width=\textwidth]{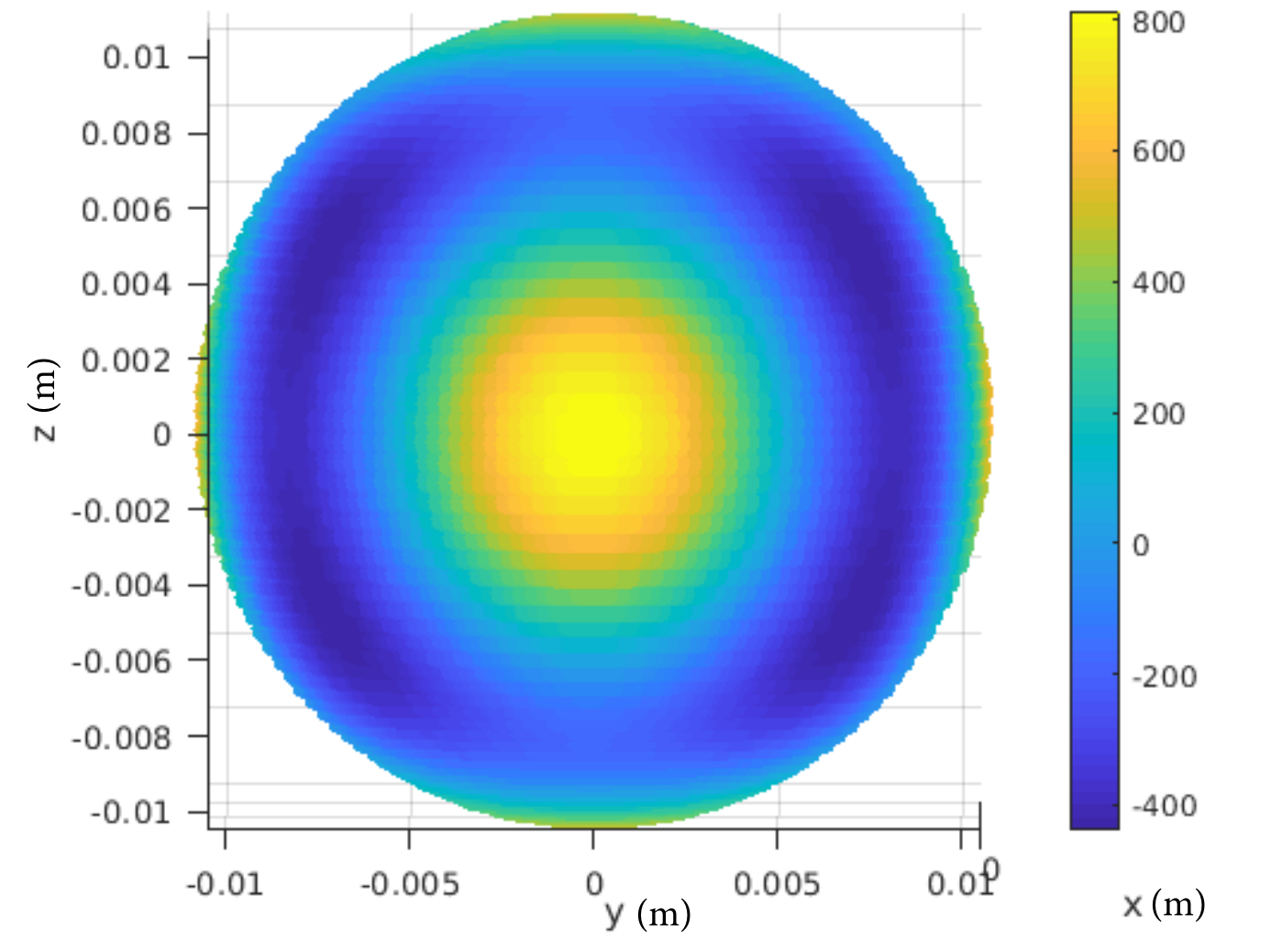}
		\label{kw_anr_front}
		\caption{back}
	\end{subfigure}
	
	\begin{subfigure}{\figsizeB}
		\includegraphics[width=\textwidth]{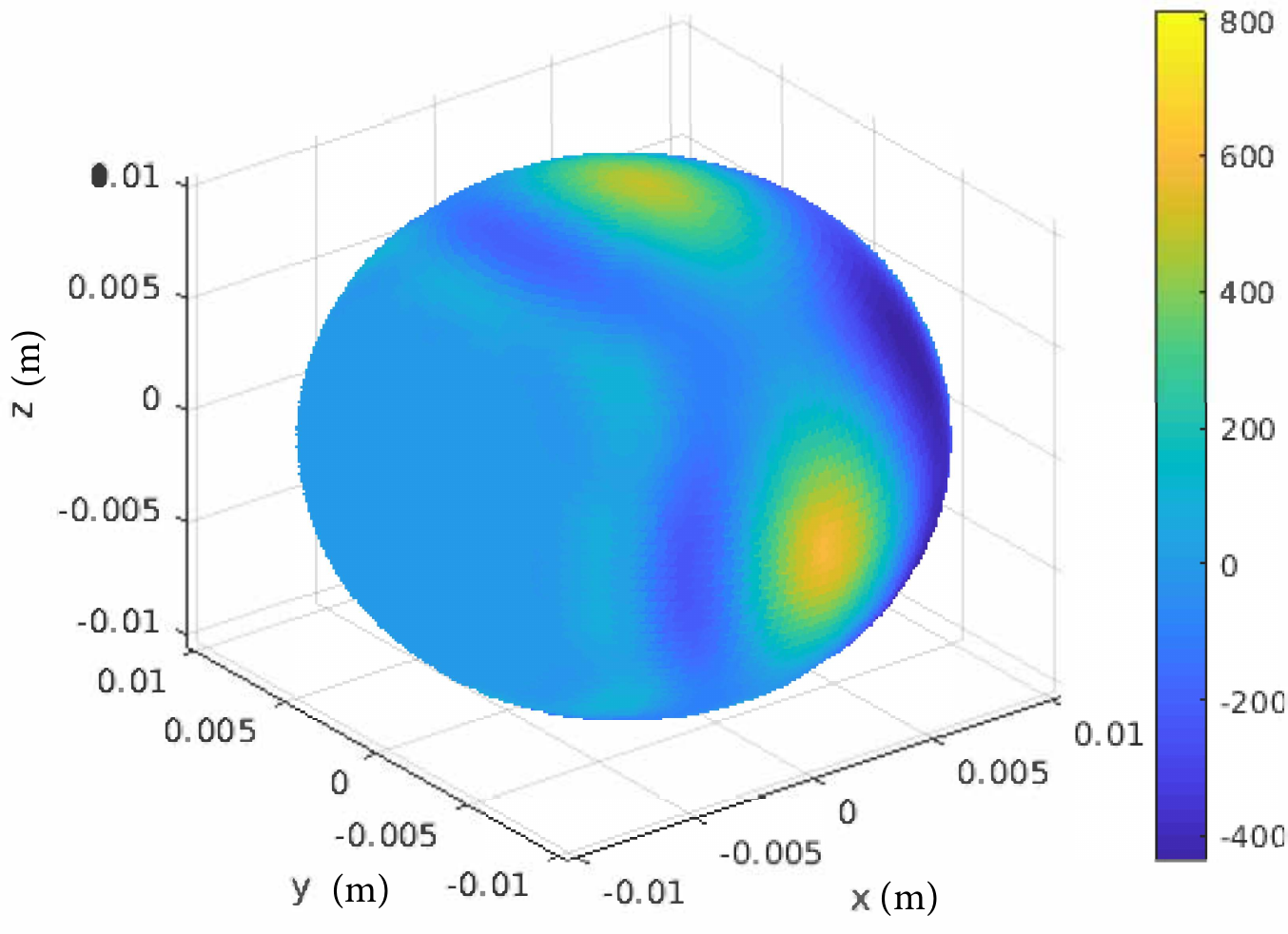}
		\label{kw_anr_angle}
		\caption{side}
	\end{subfigure}
	
	\caption{Different views of the surface field pattern on the actual source $\partial B_{0.0105}(\mathbf 0)$.}
	\label{kw_anr}
\end{figure}

Recently, with the current updated scheme we can accurately compute and evaluate the surface field pattern on any nearby surfaces surrounding $D_1$. Also, as suggested in Remark \ref{rem1} the actual physical source boundary $\partial D_a$ can be chosen as one needs as long as $D_a$ compactly embeds $D_{a'}$ and $(D_1\cup D_2)\cap D_a=\emptyset$. For example, in the case when the actual physical source is chosen to be the nearby sphere centered in the origin and radius $a=0.015cm$, i.e., $D_a=B_{0.0105}(\mathbf 0)$ the Dirichlet boundary input required for the desired control level is, as showed in \eqref{eqnpb}, the restriction of the field $u(\Bx)$ described at \eqref{remark3} to $\partial B_{0.0105}(\mathbf 0)$ and is shown in Figure \ref{kw_anr}. 
As one can observe in this picture, the surface pattern necessary on the physical source  $D_a=B_{0.0105}(\mathbf 0)$ for the desired control effect, i.e., approximation of an incoming plane wave pattern in $D_1$ with very low field values in $D_2$, exhibits a similar behavior as was recalled above about the respective density $w_\alpha$ (high and fast oscillations on the side opposite to region $D_1$ (back side in the figure) with small and faster oscillations on the side facing $D_1$ (front side in the figure)). The important difference is that the required power on $\partial B_{0.0105}(\mathbf 0)$ would be a few orders of magnitude smaller then $||w_\alpha||_{L^2}$ (which is $O(10^{5})$ as shown in Figure \ref{rotantnorm}) and has less complexity in the pattern with much smaller contrast in the level of oscillations throughout its surface. This is another fact which motivates us to believe that an optimal shape design for the actual source boundary could lead to a less challenging source pattern required for a good control.

\subsubsection{Varying the distance of the near control}
\label{var1}
\begin{figure}[h] \centering
	\includegraphics[scale=0.4]{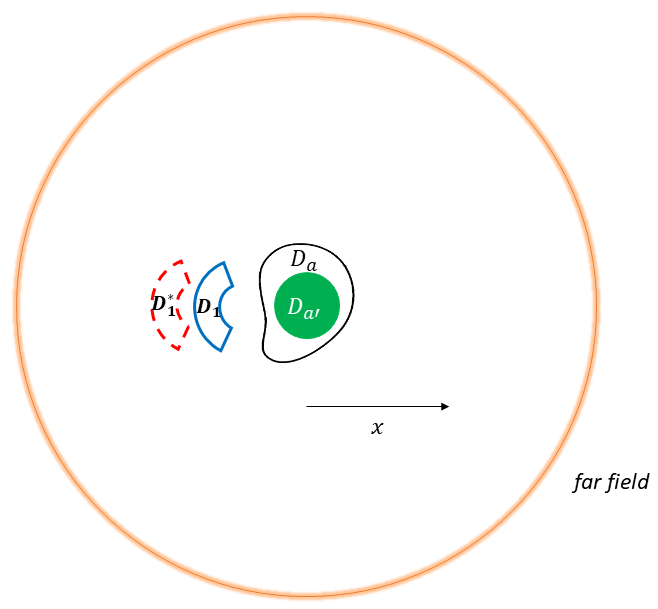}
	\caption{The near control region after an outward shift from the source.}
	\label{varyingdistance}
\end{figure}

The first test of our sensitivity analysis considers the dependence of our scheme (i.e., focusing on physical quantities defined above in the beginning of Section \ref{numresults}) on the distance between $D_1$ and $D_{a'}$. Two iterations of this experiment are shown in Figure \ref{varyingdistance}. $D_1^*$ is the near control region after an outward shift from the antenna.

In what follows in this subsection, in each figure, the sensitivity plots are presented as function of the distance between $D_1$ and $D_{a'}$, first up to a distance of $0.025m$ (left plot) and then going further to a distance of $0.28m$ (right plot).

\begin{figure}[h] \centering
	
	\begin{subfigure}{\figsizeB}
		\includegraphics[width=\textwidth]{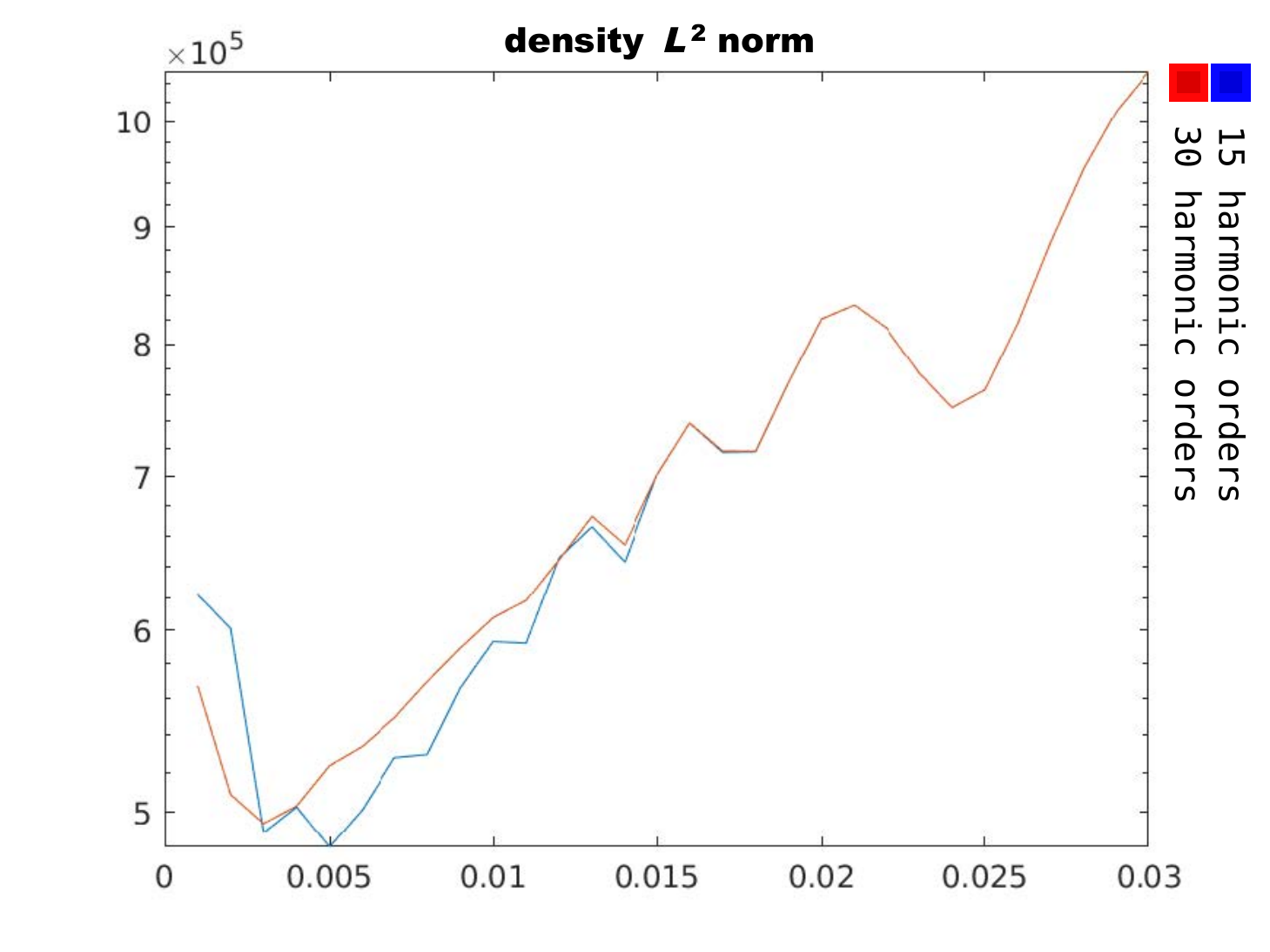}
		\label{chart:anr10hShift001_antnorm}
		\vspace{-0.5cm}
		\caption{Control near the source}
	\end{subfigure}
	\begin{subfigure}{\figsizeB}
		\includegraphics[width=\textwidth]{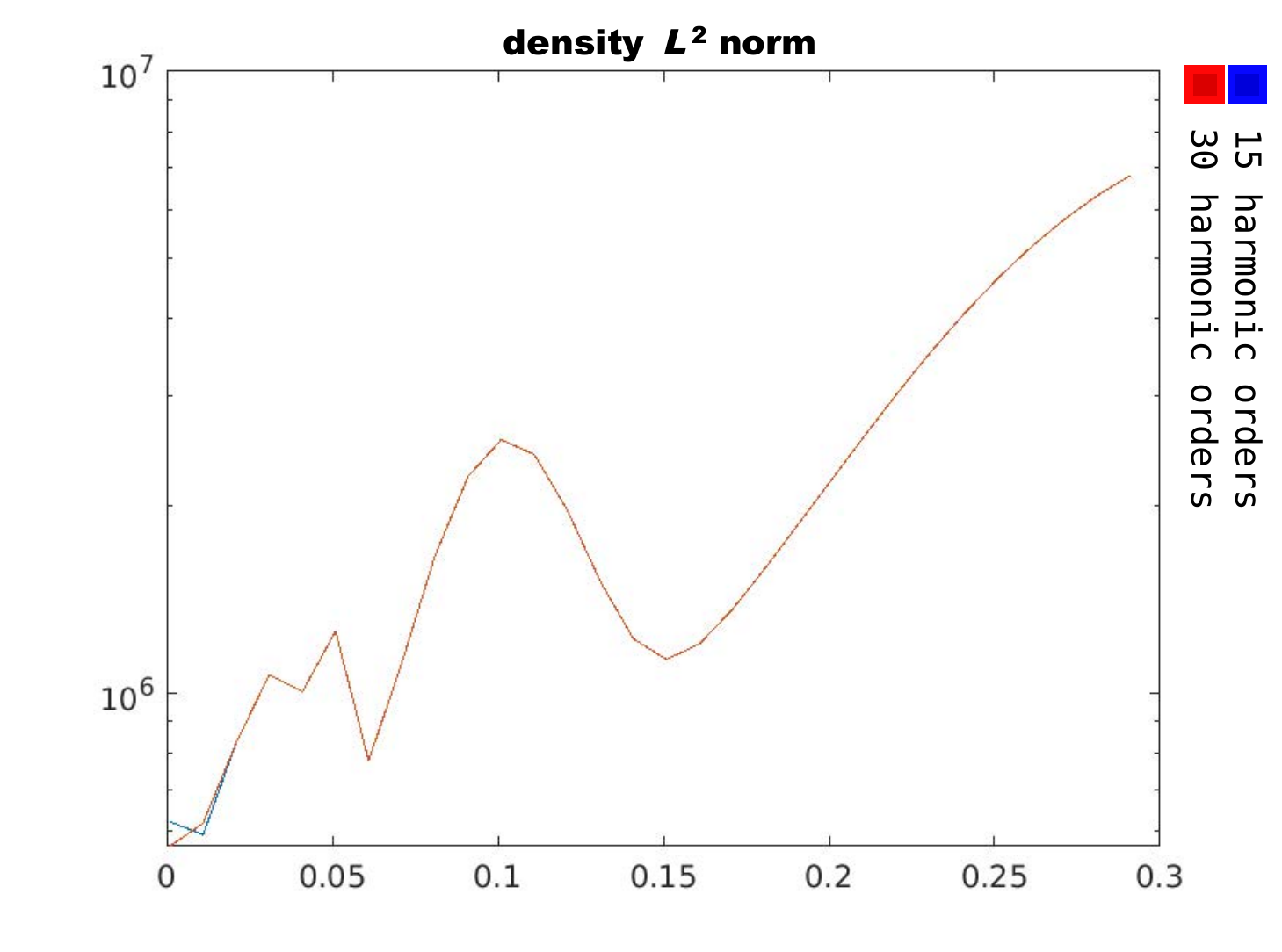}
		\label{chart:anr10hShift01_antnorm}
		\vspace{-0.5cm}
		\caption{Control further from the source}
	\end{subfigure}
	\caption{$L^2$ norm of the source density $w_\alpha$ as a function of the distance between $D_1$ and $D_{a'}$. }
	\label{2.42}
\end{figure}

In Figure $\ref{2.42}$, the $L^2$ norm of the density $w_\alpha$ is shown as a function of the distance between $D_1$ and $D_{a'}$. This quantity is an indication of the source power and it can be seen that it is bigger for the larger number of harmonics but in both scenarios (i.e., 15 or 30 harmonic orders)  grows exponentially. 

\begin{figure}[!htb] \centering
	\begin{subfigure}{\figsizeB}
		\includegraphics[width=\textwidth]{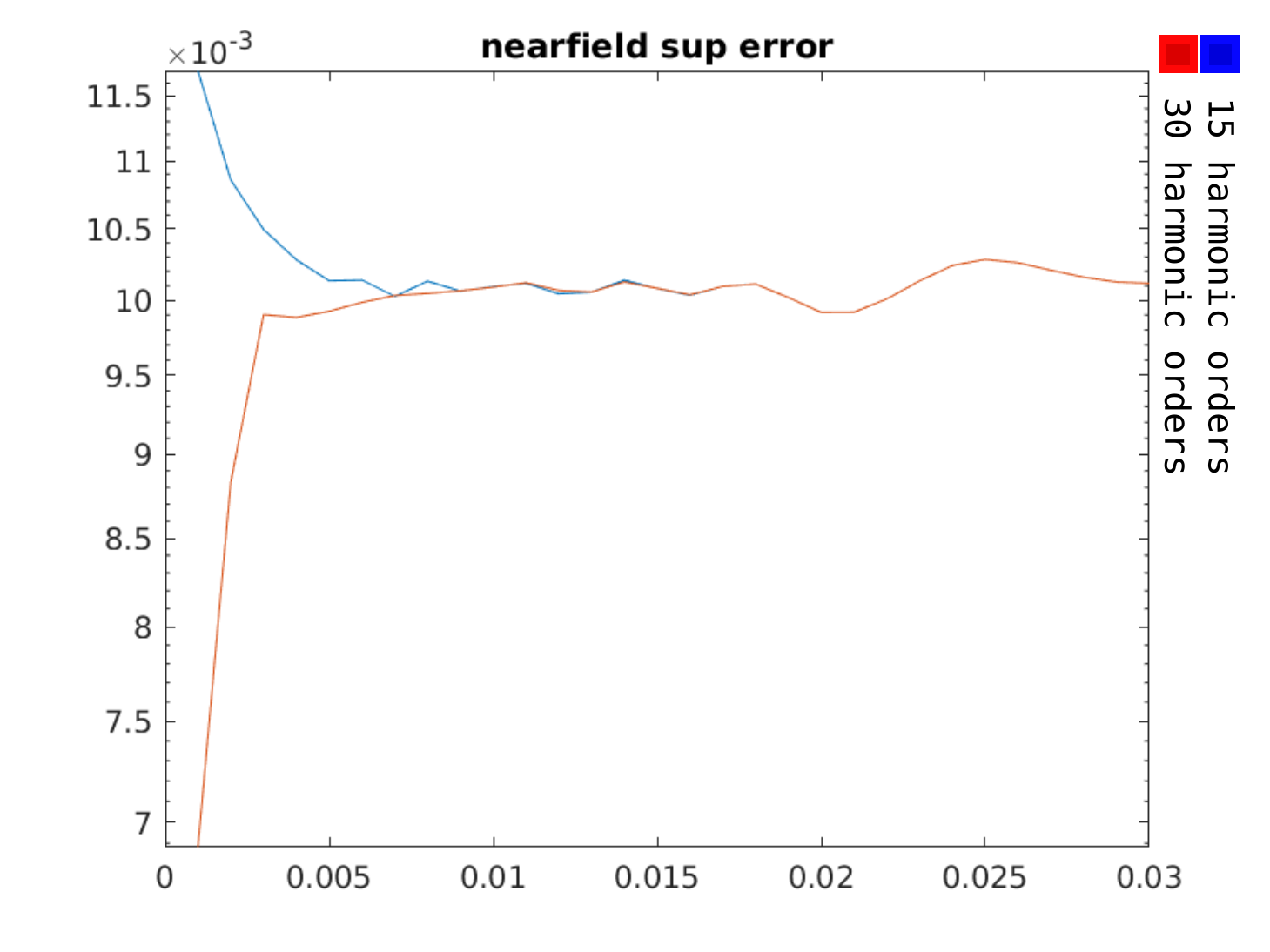}
		\label{chart:anr10hShift001_nearsup}
		\vspace{-0.5cm}
		\caption{Control near the source}
	\end{subfigure}
	\begin{subfigure}{\figsizeB}
		\includegraphics[width=\textwidth]{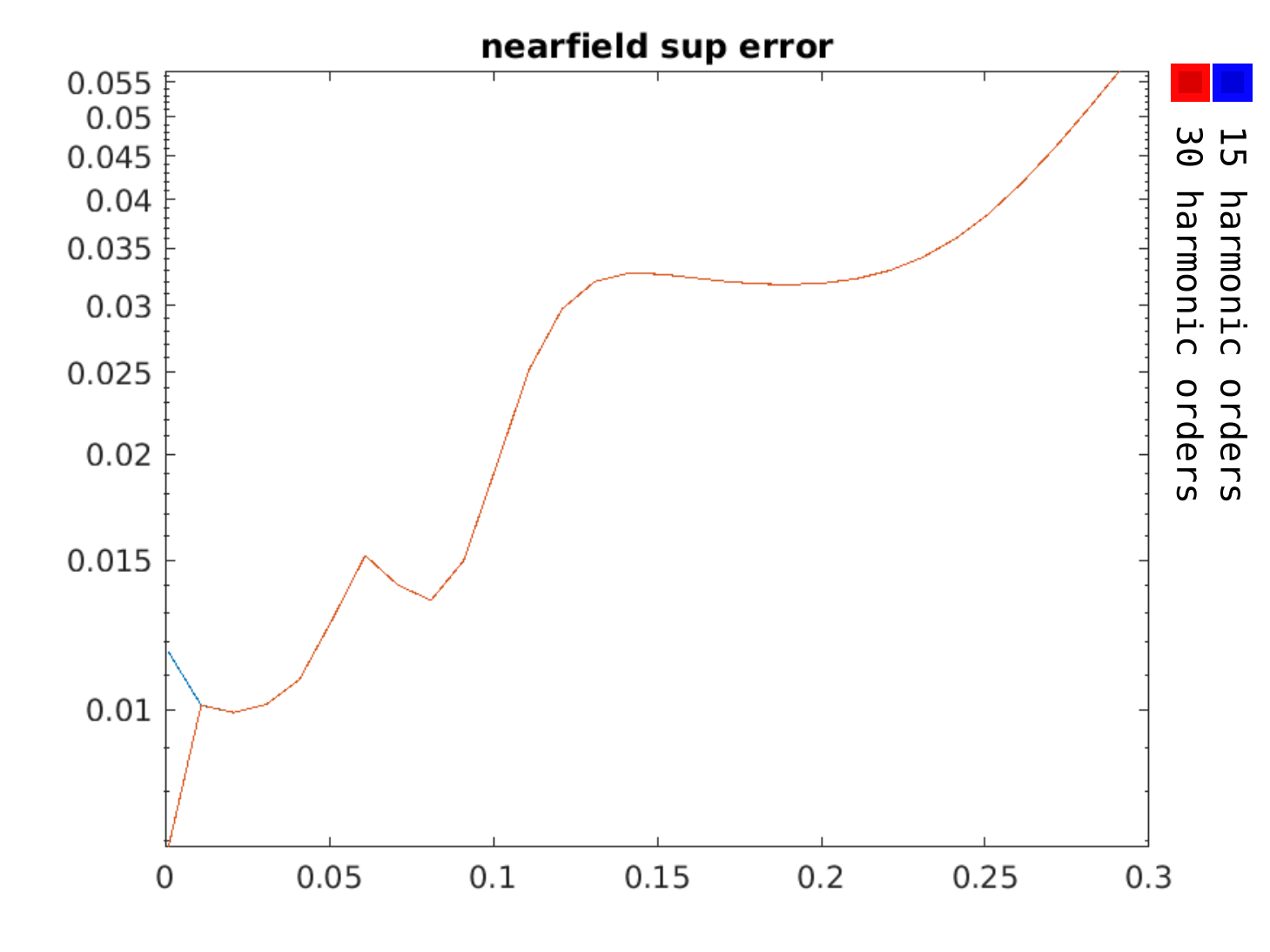}
		\label{chart:anr10hShift01_nearsup}
		\vspace{-0.5cm}
		\caption{Control further from the source}
	\end{subfigure}
	\caption{Relative supremum error in $D_1$ as a function of the distance between $D_1$ and $D_{a'}$}
	\label{2.44}
\end{figure}

\begin{figure}[!htb] \centering
	\begin{subfigure}{\figsizeB}
		\includegraphics[width=\textwidth]{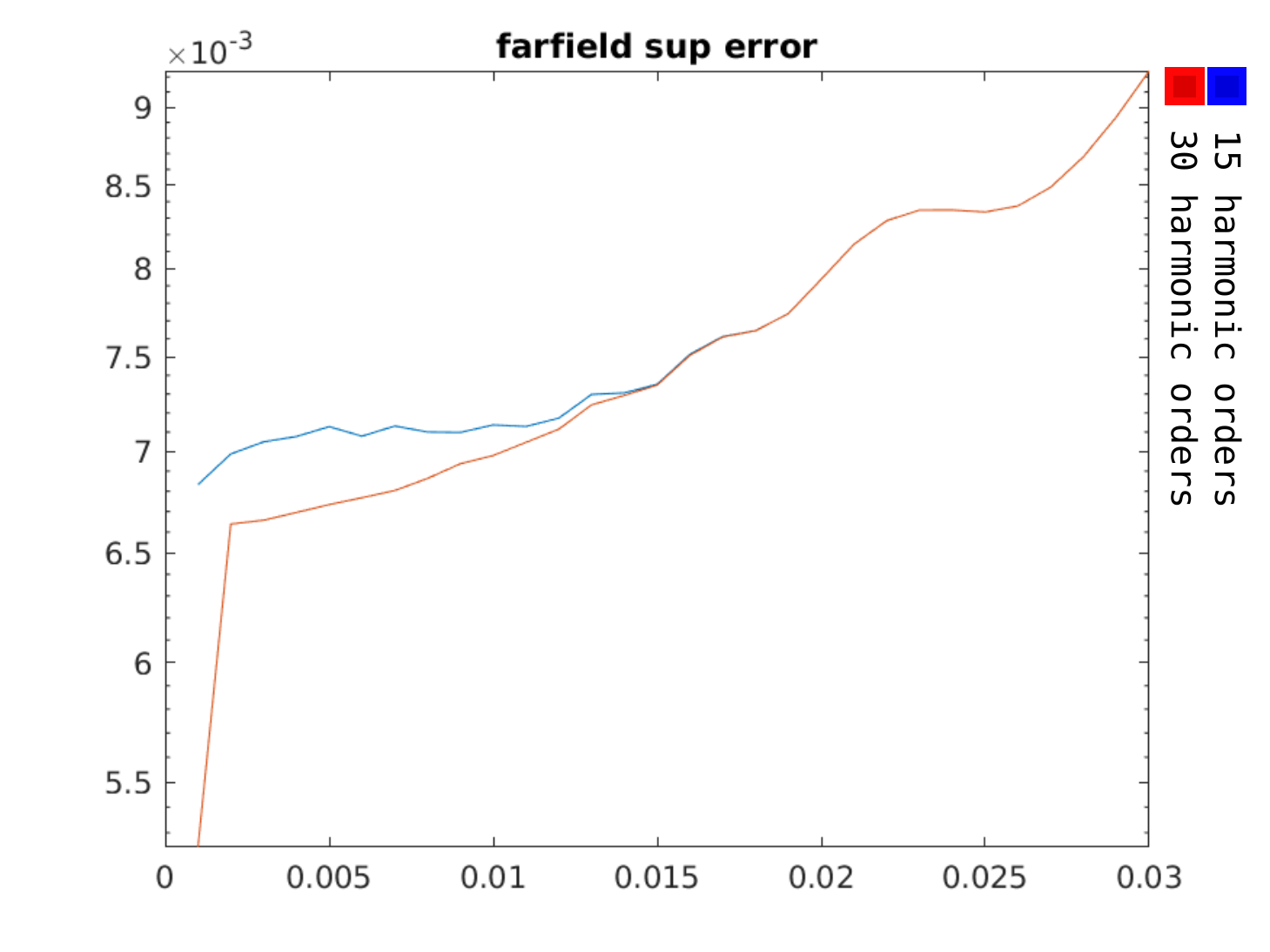}
		\label{chart:anr10hShift001_farsup}
		\vspace{-0.5cm}
		\caption{$D_1$ near the source}
	\end{subfigure}
	\begin{subfigure}{\figsizeB}
		\includegraphics[width=\textwidth]{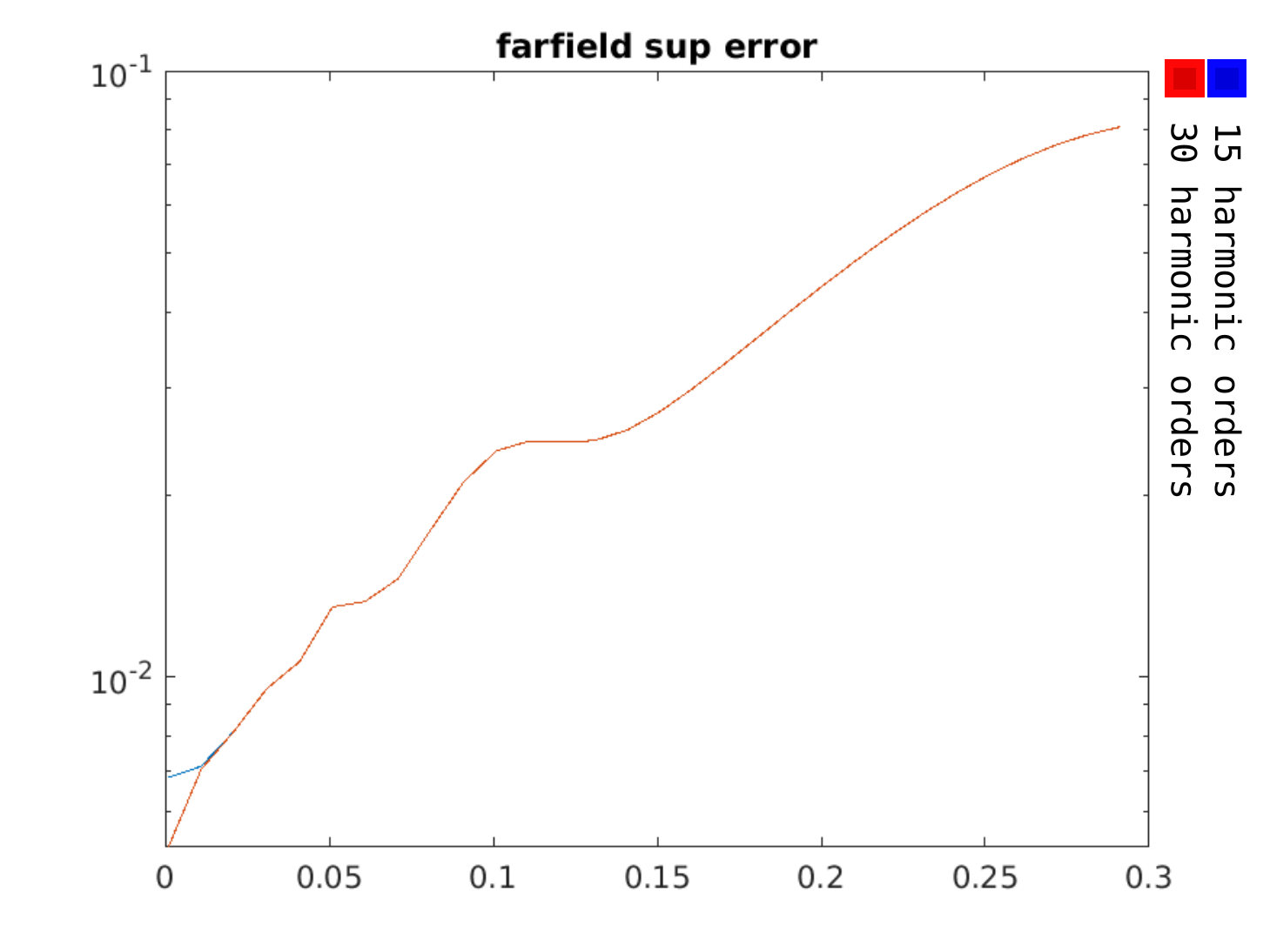}
		\label{chart:anr10hShift01_farsup}
		\vspace{-0.5cm}
		\caption{$D_1$ further from the source}
	\end{subfigure}
	\caption{Absolute supremum error on $\partial D_2$ as a function of the distance between $D_1$ and $D_{a'}$}
	\label{2.45}
\end{figure}

\begin{figure}[!htb] \centering
	\begin{subfigure}{\figsizeB}
		\includegraphics[width=\textwidth]{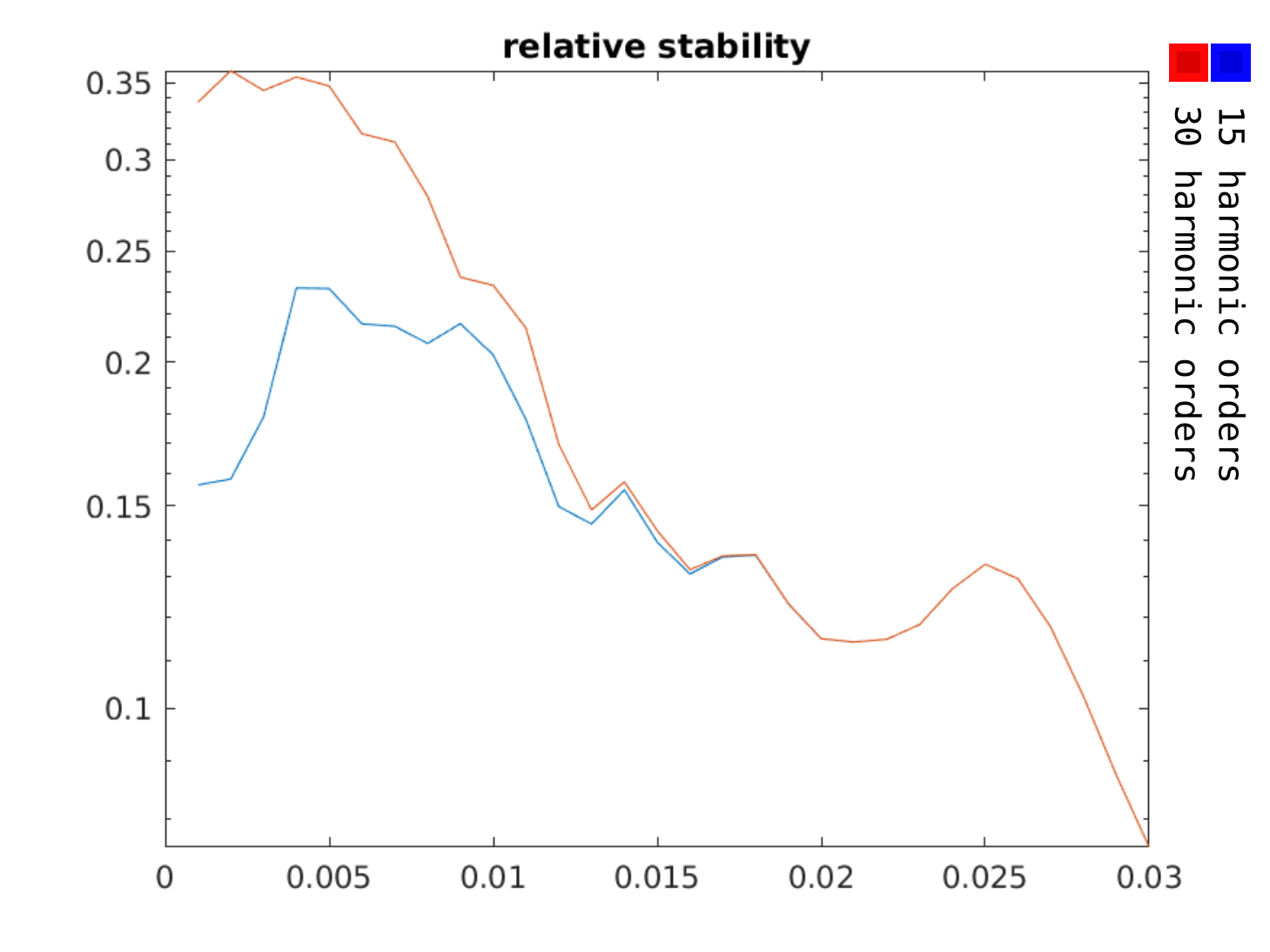}
		\label{chart:anr10hShift001_relstab}
		\vspace{-0.5cm}
		\caption{Control near the source}
	\end{subfigure}
	\begin{subfigure}{\figsizeB}
		\includegraphics[width=\textwidth]{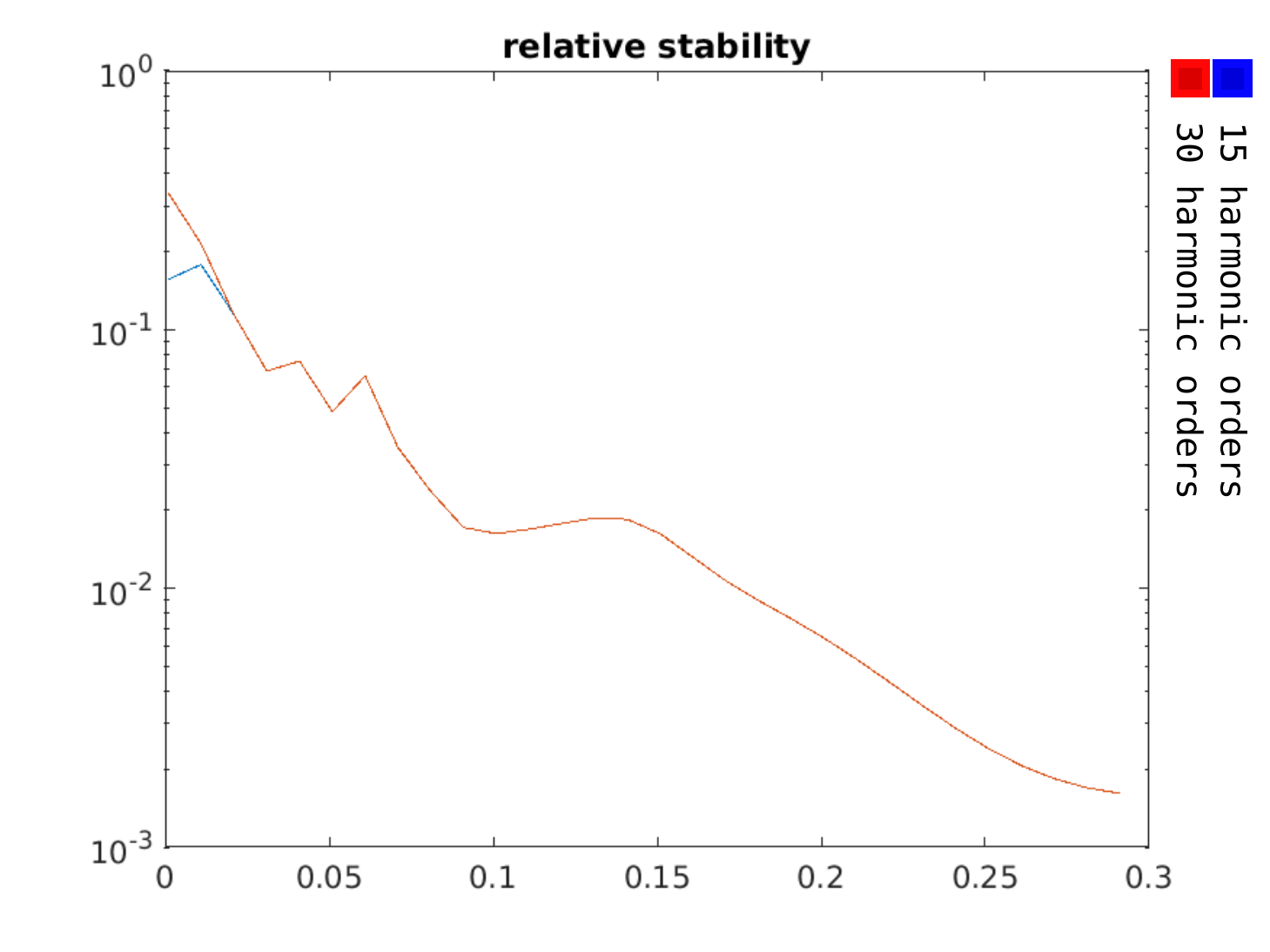}
		\label{chart:anr10hShift01_relstab}
		\vspace{-0.5cm}
		\caption{Control further from the source}
	\end{subfigure}
	\caption{Stability of the scheme as a function of of the distance between $D_1$ and $D_{a'}$.}
	\label{2.46}
\end{figure}

Figure \ref{2.44} shows the relative supremum error in $D_1$ and Figure \ref{2.45} shows the absolute supremum error on $\partial D_2$. In Figure \ref{2.46} we present our results concerning the stability of the scheme with respect to the distance between $D_1$ and $D_{a'}$. 

We observe that the accuracy error in $D_1$ is  better in the vicinity of the source but remains of order $O(10^{-2})$ throughout. On the other hand the far field absolute supremum error reaches undesirable levels when $D_1$ and $D_{a'}$ are at distances greater then $15$cm. Overall, except the near field when more harmonic orders produce a better accuracy we observe that the number of harmonics used does not seem to make a difference in the scheme behavior. In the same time the stability is worse when $D_1$ is in the vicinity of the source with the observation that the scheme is more stable when less harmonic orders are used.  

As a conclusion, for this sensitivity test, there seem to be a competition between accuracy and stability depending on the number of harmonics being used and the distance between $D_1$ and the fictitious source $D_{a'}$. We mention that in applications where the field to be approximated in region $D_1$ is prescribed a priori the stability analysis may not be so relevant but what may be challenging for a practical realization is the synthesis of a radiator with the complexity suggested by our simulations. As we suggested above, based on Remark \ref{rem1}, in order to better address the source synthesis challenge one could employ optimization procedures for finding the best possible shape for the physical source $D_a$ or one could consider different penalty functionals in the optimization scheme.

\subsubsection{Varying the outer radius of the near control}
\label{var2}
\begin{figure}[!htb] \centering
	\begin{subfigure}{\figsizeB}
		\includegraphics[width=\textwidth]{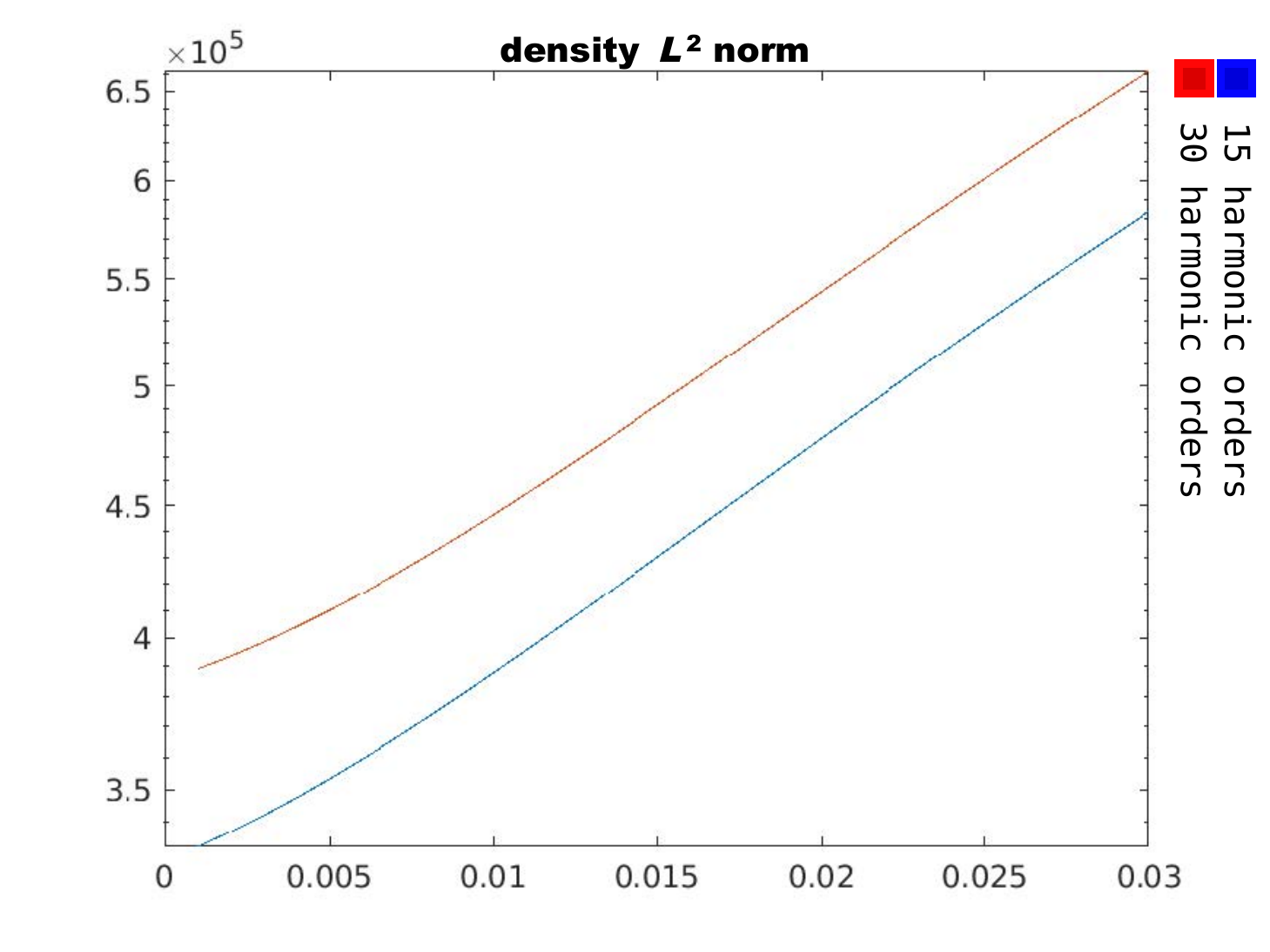}
		\label{2.47}
		\vspace{-0.5cm}
		\caption{  Outer radius in $[0.015, 0.03]$}
	\end{subfigure}
	\begin{subfigure}{\figsizeB}
		\includegraphics[width=\textwidth]{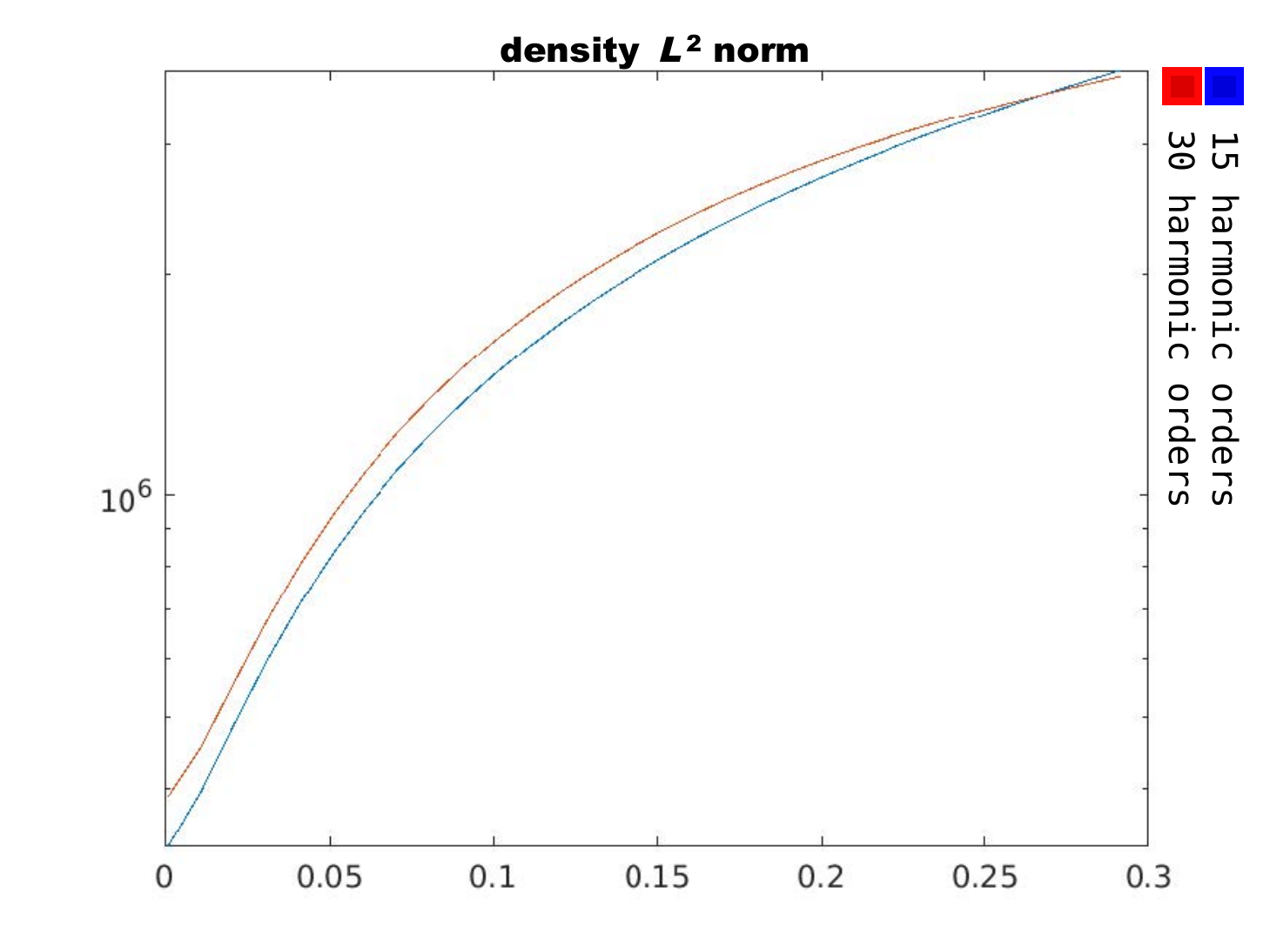}
		\label{2.48}
		\vspace{-0.5cm}
		\caption{ Outer radius in $[0.015, 0.3]$  }
	\end{subfigure}
	
	\caption{ $L^2$ norm of the source density $w_\alpha$ as a function of the increments in $D_1$ outer radius.}
	\label{2.47-2.48}
\end{figure}

In the next test of our sensitivity analysis, we consider the behavior of the above physically important quantities (see beginning of Section \ref{numresults}) with respect to increments of the outer radius of the sectorial near control region $D_1$ described at \eqref{D1}  while the distance between $D_1$ and the source is kept fixed. 
\begin{figure}[!htb] \centering
	\includegraphics[scale=0.4]{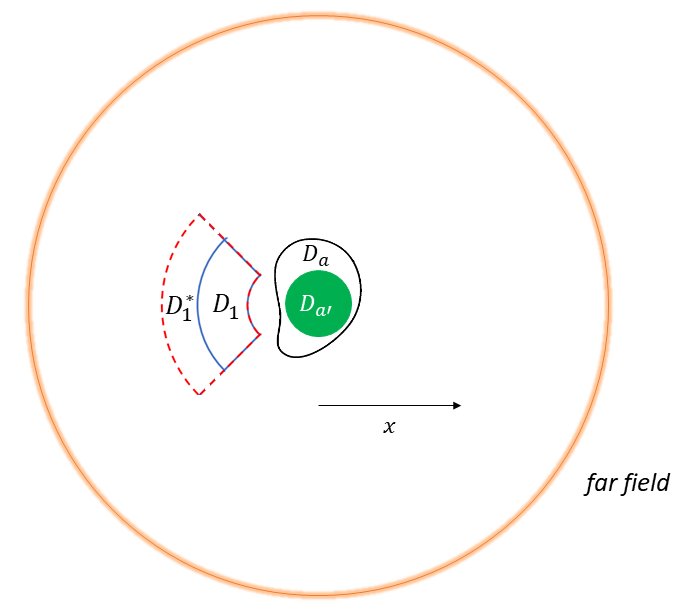}
	\caption{The initial near control $D_1$ and an iterate $D_1^*$ after increasing the outer radius. }
	\label{varyingouterrad}
\end{figure}
Figure \ref{varyingouterrad} shows an iterate $D_1^*$ of the near control after an increase in the outer radius. Figure \ref{2.47-2.48} shows the $L^2$ norm of the density $w_\alpha$ as a function of the increments in $D_1$ outer radius. We can observe that a larger $||w_\alpha||_{L^2}$ is required for the 30 harmonic orders as expected. Also, since $||w_\alpha||_{L^2}$ is an indicator of the source power, it can be seen that for larger control regions the power will tend to grows exponentially at approximately the same rate regardless of  the number of harmonics used.

\begin{figure}[!htb] \centering
	\begin{subfigure}{\figsizeB}
		\includegraphics[width=\textwidth]{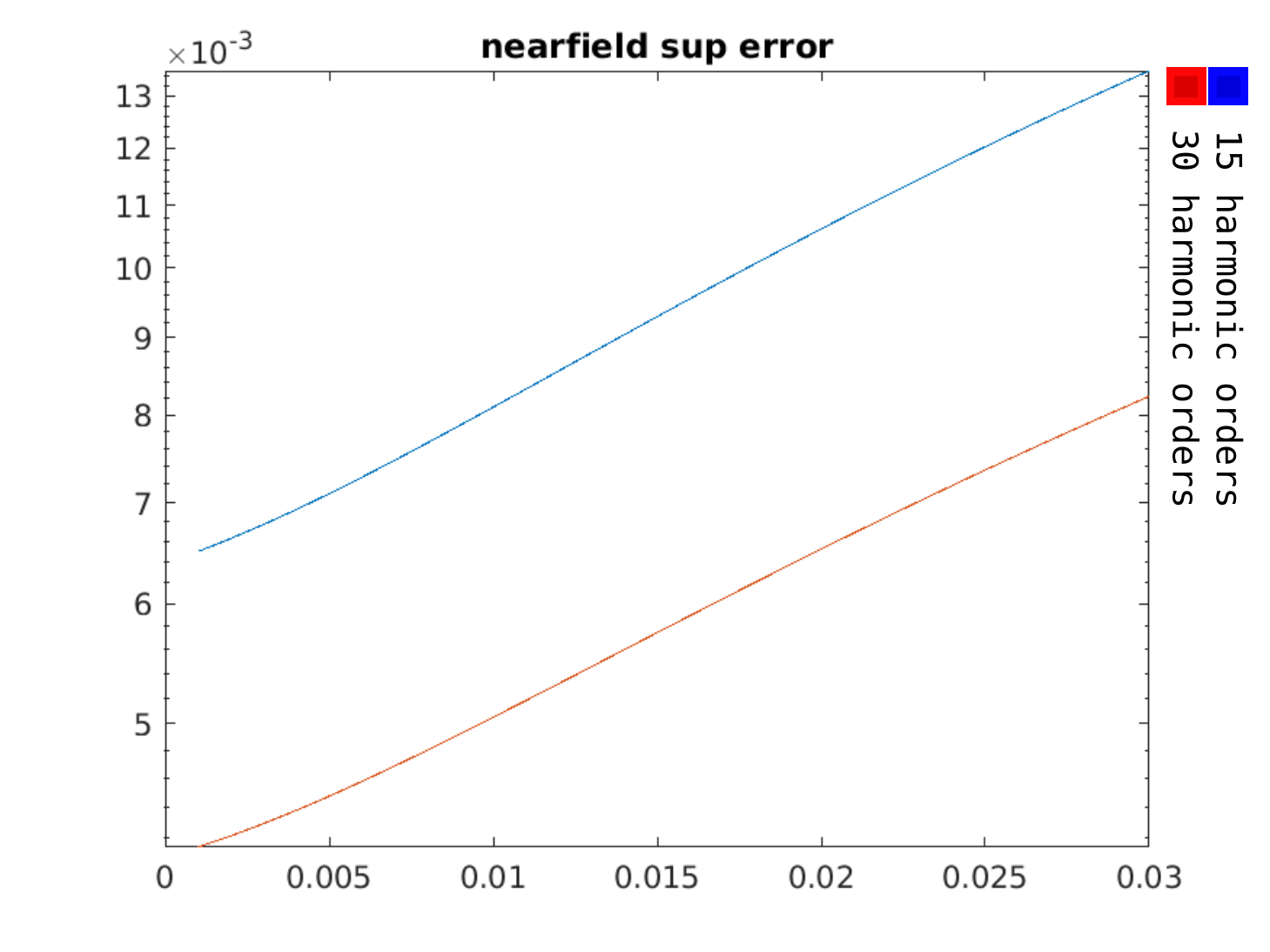}
		\label{chart:anr10hOuterR001_nearsup}
		\vspace{-0.5cm}
		\caption{  Outer radius in $[0.015, 0.03]$}
	\end{subfigure}
	\begin{subfigure}{\figsizeB}
		\includegraphics[width=\textwidth]{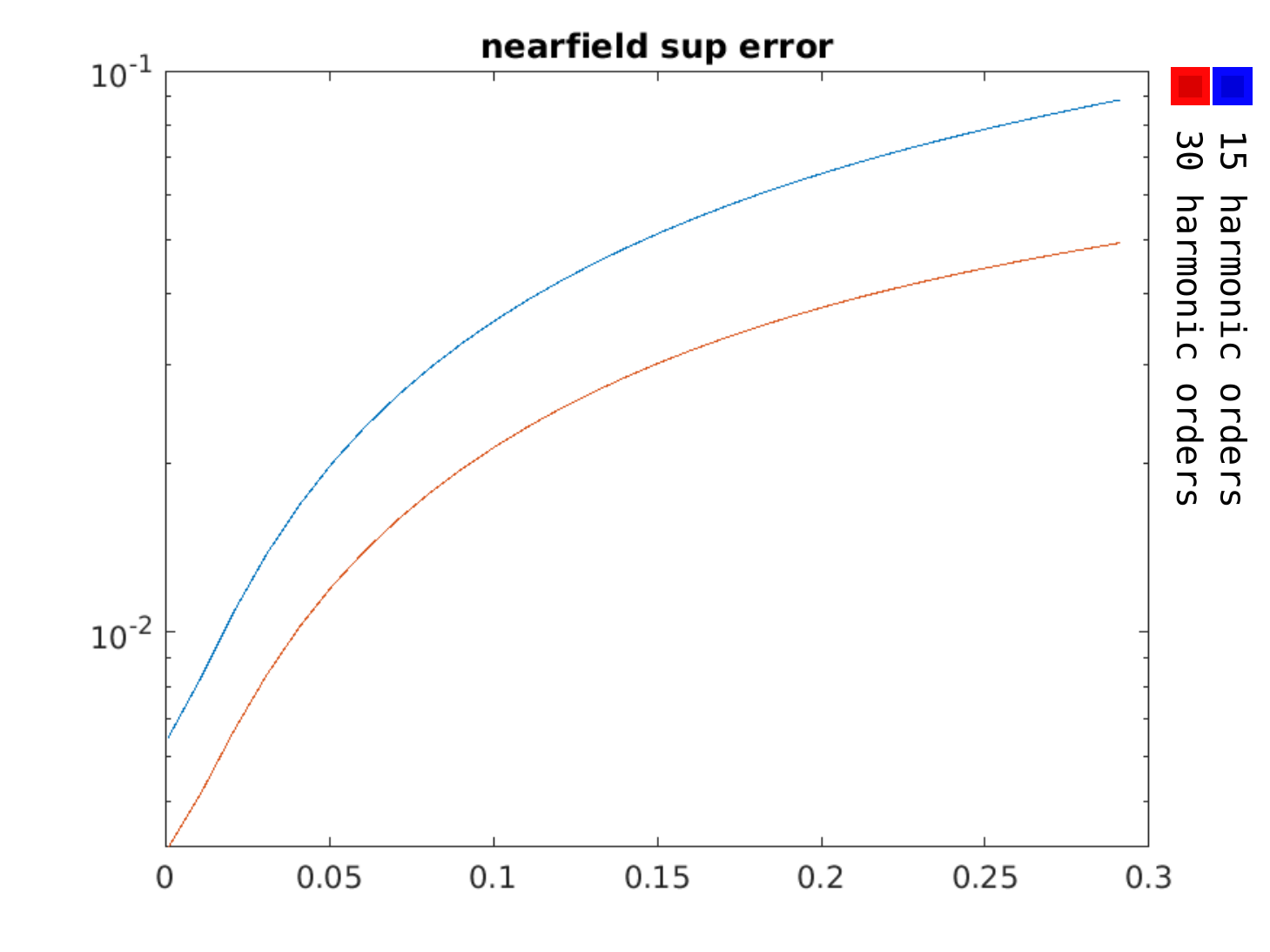}
		\label{chart:anr10hOuterR001_farsup}
		\vspace{-0.5cm}
		\caption{  Outer radius in $[0.015, 0.3]$}
	\end{subfigure}
	\caption{Relative supremum error in $D_1$ as a function of the increments in $D_1$ outer radius.}
	\label{2.49-2.50}
\end{figure}

\begin{figure}[!htb] \centering
	\begin{subfigure}{\figsizeB}
		\includegraphics[width=\textwidth]{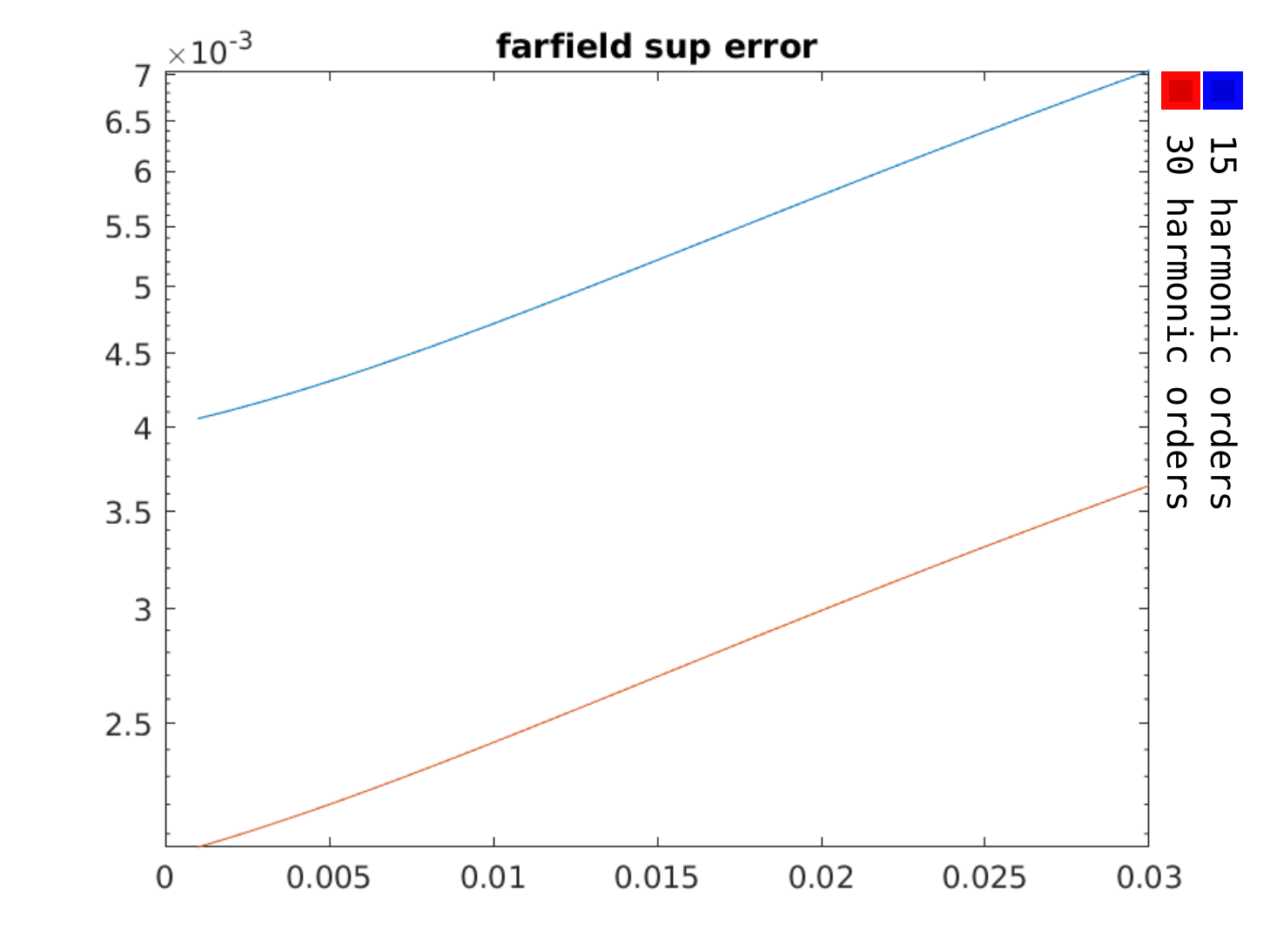}
		\label{chart:anr10hOuterR001_nearsup}
		\vspace{-0.5cm}
		\caption{  Outer radius in $[0.015, 0.03]$}
	\end{subfigure}
	\begin{subfigure}{\figsizeB}
		\includegraphics[width=\textwidth]{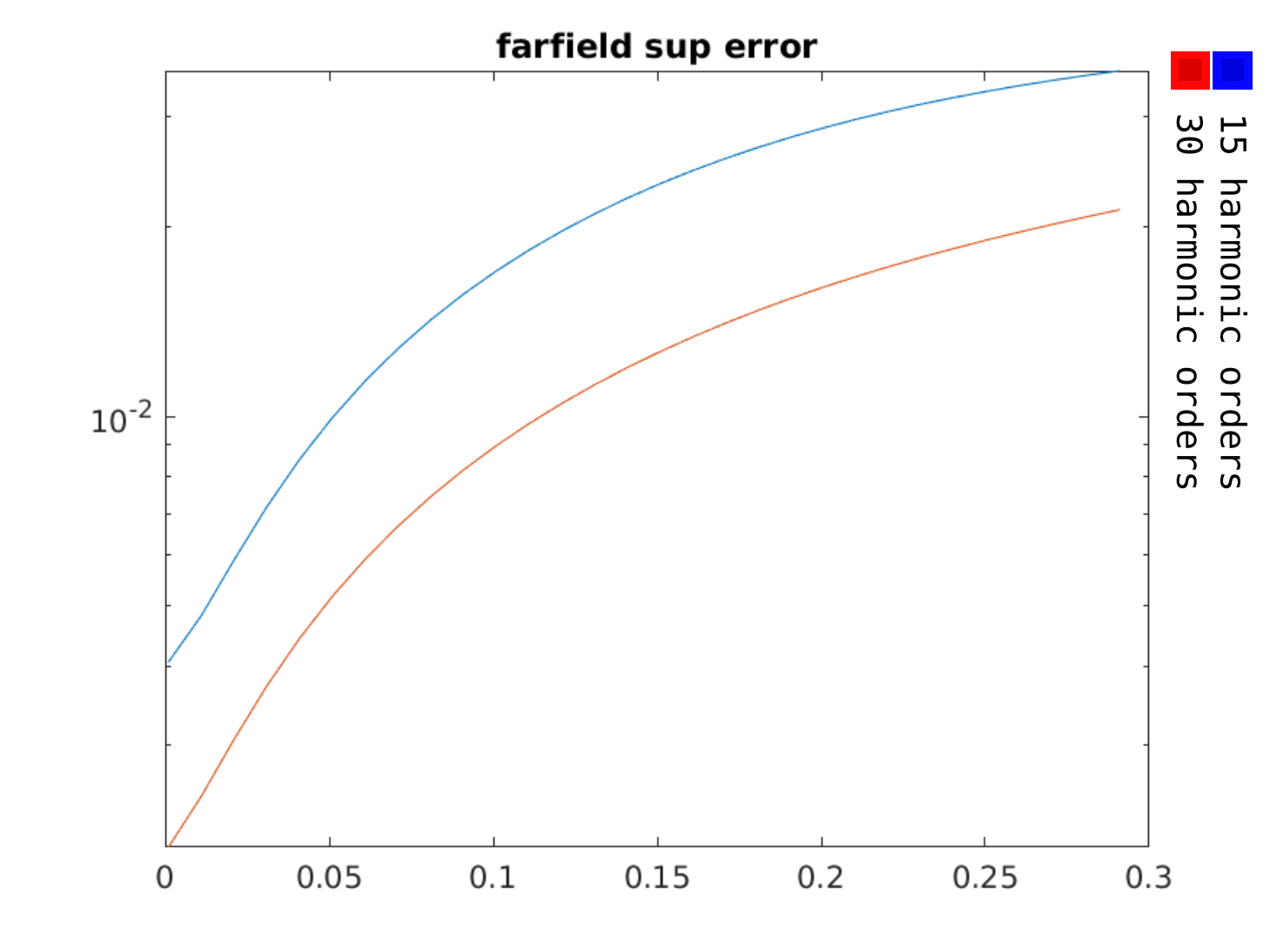}
		\label{chart:anr10hOuterR001_farsup}
		\vspace{-0.5cm}
		\caption{  Outer radius in $[0.015, 0.3]$}
	\end{subfigure}
	\caption{Absolute supremum error on $\partial D_2$ as a function of the increments in $D_1$ outer radius.}
	\label{2.49-2.50-far}
\end{figure}

In Figure \ref{2.49-2.50} and Figure \ref{2.49-2.50-far} one can observe that the relative supremum error in region $D_1$ and the absolute supremum error on $\partial D_2$ are of order $O(10^{-3})$ in the vicinity of the source and then reach order $O(10^{-2})$ (outer radius approximatively $2.5$cm) and grow to reach undesirable values when the region of control $D_1$ gets larger (outer radius approximatively $15$cm). This fact suggest the necessity of more harmonics to be used for greater accuracy in the context of larger control regions.

\begin{figure}[!htb] \centering
	\begin{subfigure}{\figsizeB}
		\includegraphics[width=\textwidth]{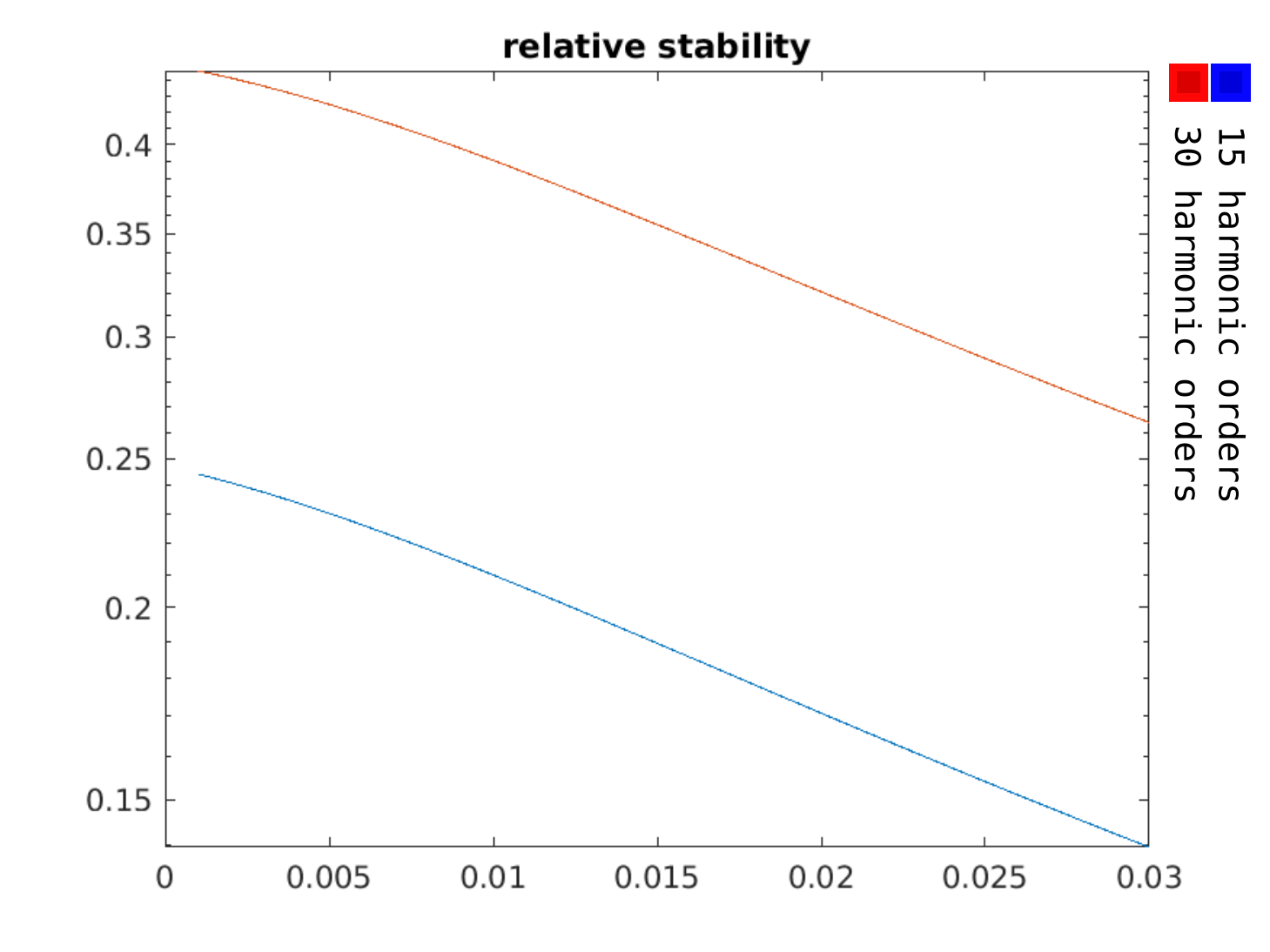}
		\label{chart:anr10hOuterR001_relstab}
			\vspace{-0.5cm}
			\caption{  Outer radius in $[0.015, 0.03]$}
	\end{subfigure}
	\begin{subfigure}{\figsizeB}
		\includegraphics[width=\textwidth]{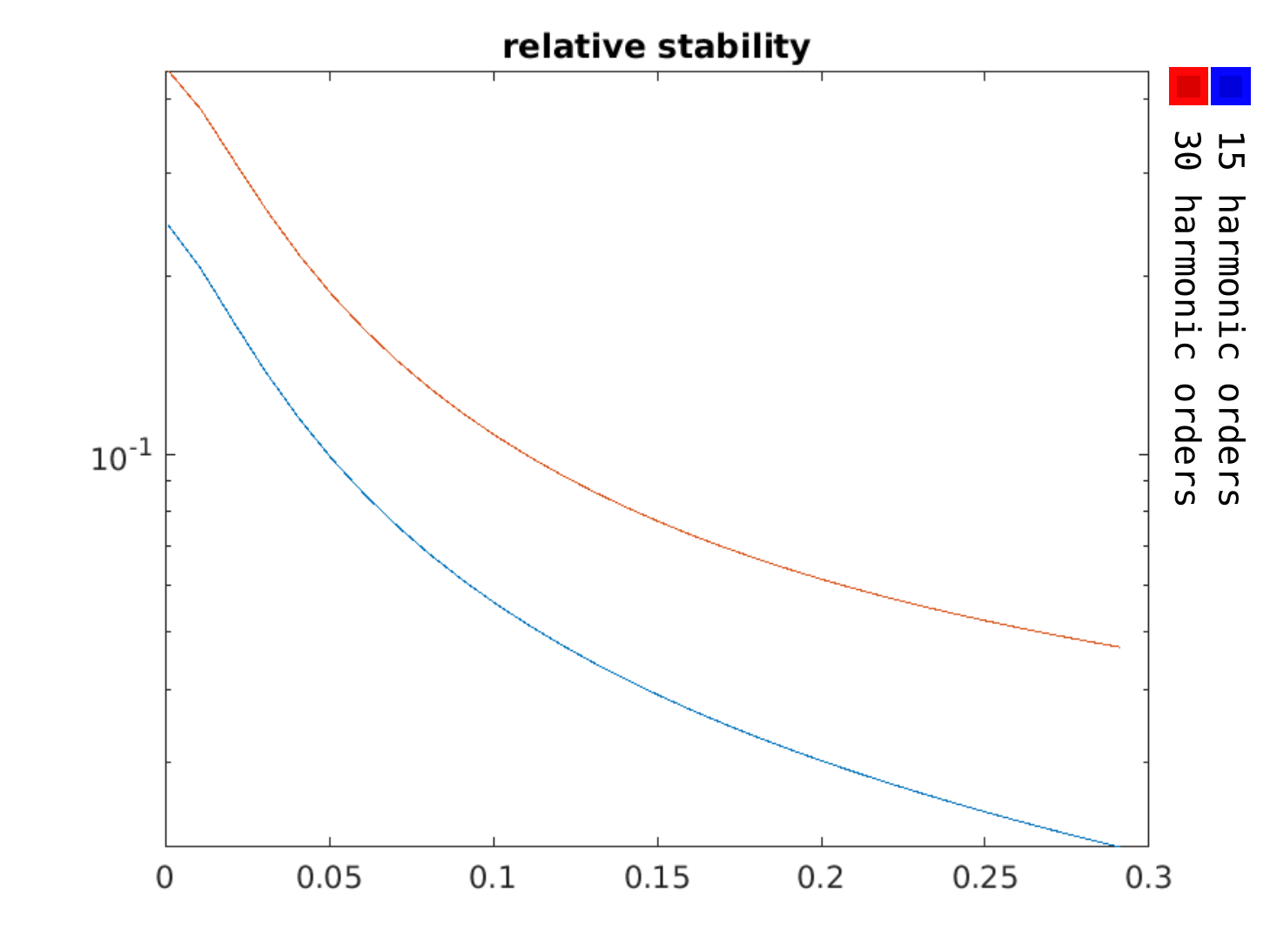}
		\label{chart:anr10hOuterR01_relstab}
		\vspace{-0.5cm}
		\caption{  Outer radius in $[0.015, 0.3]$}
	\end{subfigure}
	\caption{Stability measure as a function of the increments in $D_1$ outer radius.}
	\label{2.51}
\end{figure}

Figure \ref{2.51} shows the stability getting better when region $D_1$ gets larger. Notice also that the 15 harmonic order-configuration exhibited better stability even though it is less accurate. Using higher harmonic orders remains stable while maintaining accurate results. Again, as we mentioned above, for applications where the field to be approximated in $D_1$ is a priori known the stability is not an important issue but the complexity of the source to be synthesized remains a challenge which could be mitigated as we suggested at the end of Section \ref{var1}.

\subsubsection{Varying both the inner and outer radius of the near control}
\label{var3}
\begin{figure}[!htb] \centering
	\includegraphics[scale=0.4]{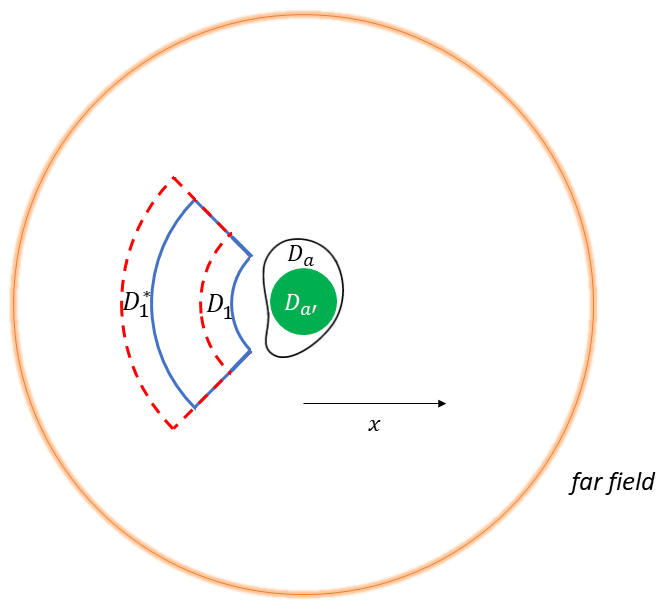}
	\caption{The initial near control $D_1$ and an iterate $D_1^*$ after increasing both the inner and outer radii by equal increments. }
	\label{varyingbothradii}
\end{figure}

The following results show the effect of increasing both the inner and outer radius of the near control, i.e, simultaneously moving the front and back sides of the near control away from the source. The increments in both radii are always kept equal. Figure \ref{varyingbothradii} illustrates two iterates $D_1$ and $D_1^*$ in this set of experiments.

\begin{figure}[!htb] \centering
	\begin{subfigure}{\figsizeB}
		\includegraphics[width=\textwidth]{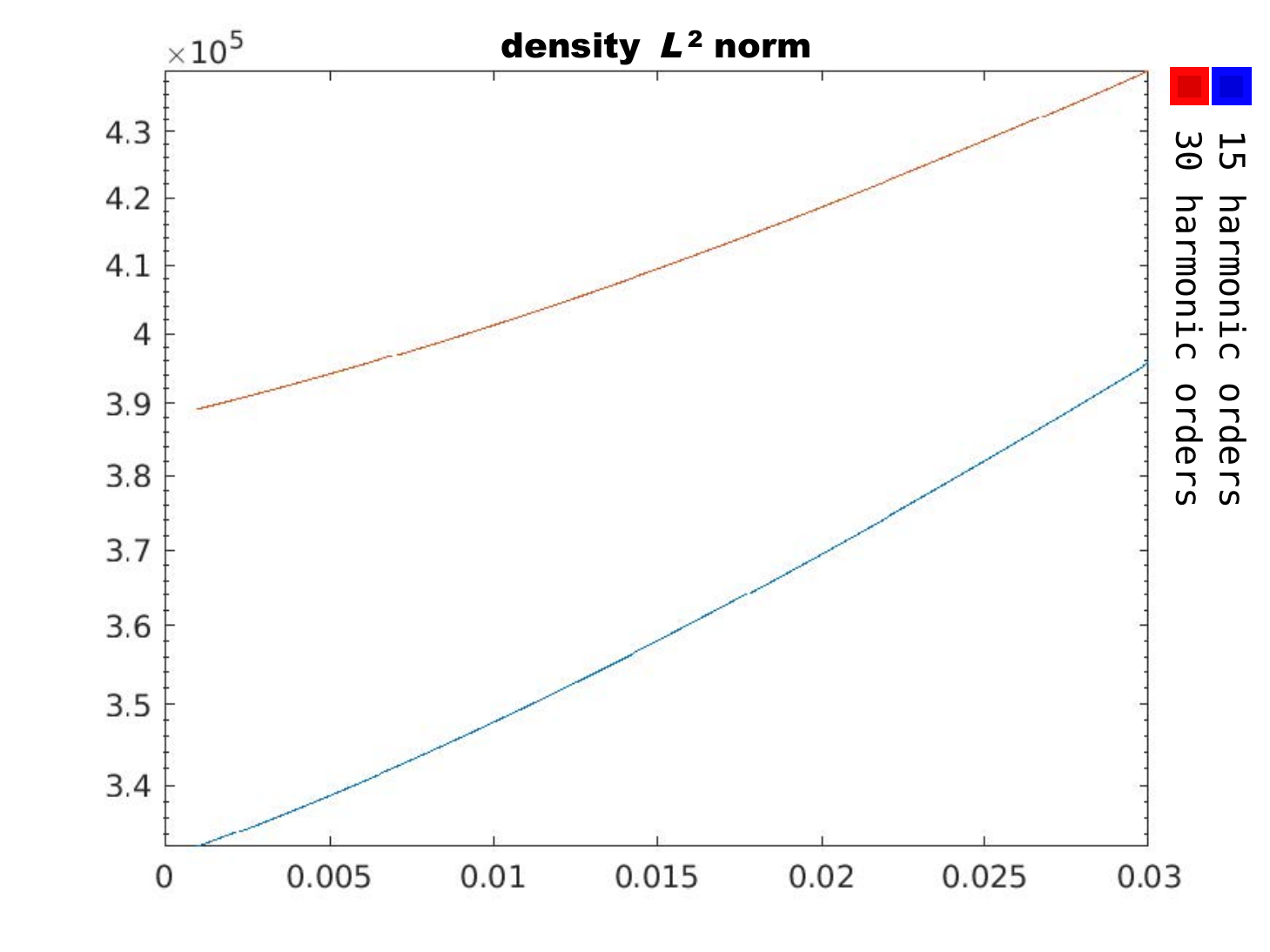}
		\label{chart:anr10hBothR001_antnorm}
			\vspace{-0.5cm}
		\caption{  Radii increments in $[0, 0.029]$}
	\end{subfigure}
	\begin{subfigure}{\figsizeB}
		\includegraphics[width=\textwidth]{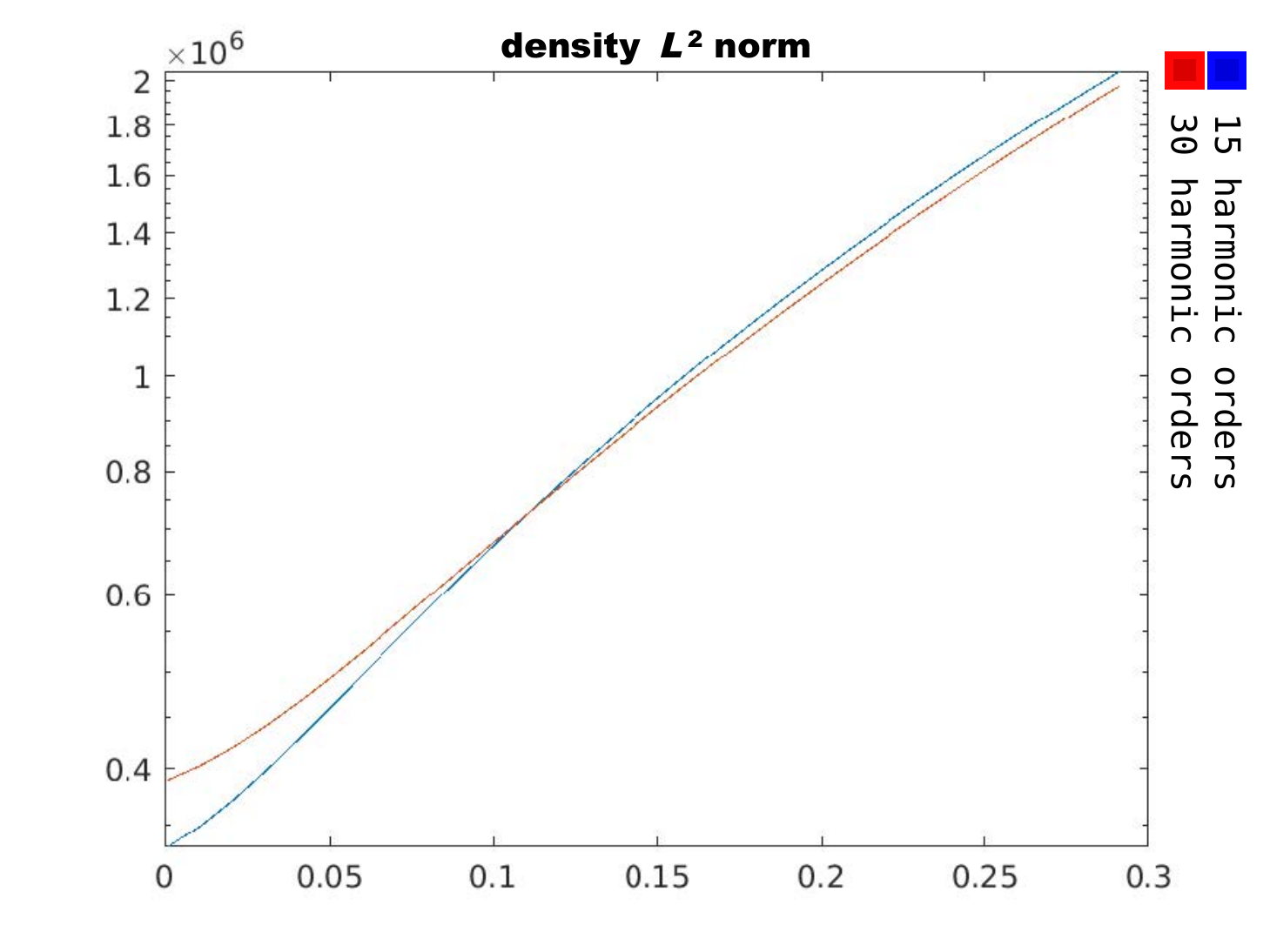}
		\label{chart:anr10hBothR01_antnorm}
		\vspace{-0.5cm}
		\caption{  Radii increments in $[0, 0.29]$}
	\end{subfigure}
	\caption{$L^2$ norm of the source density $w_\alpha$ as a function of the $D_1$ radii increment.}
	\label{2.52}
\end{figure}

Figure \ref{2.52} shows the $L^2$ norm of the density $w_\alpha$ as a function of the increments in both of $D_1$'s radii. As in the other cases discussed above, we observe that it grows exponentially with bigger values for more harmonics until the region is larger and further away where its value stabilizes for the two levels of harmonic orders used. 

\begin{figure}[!htb] \centering
	\begin{subfigure}{\figsizeB}
		\includegraphics[width=\textwidth]{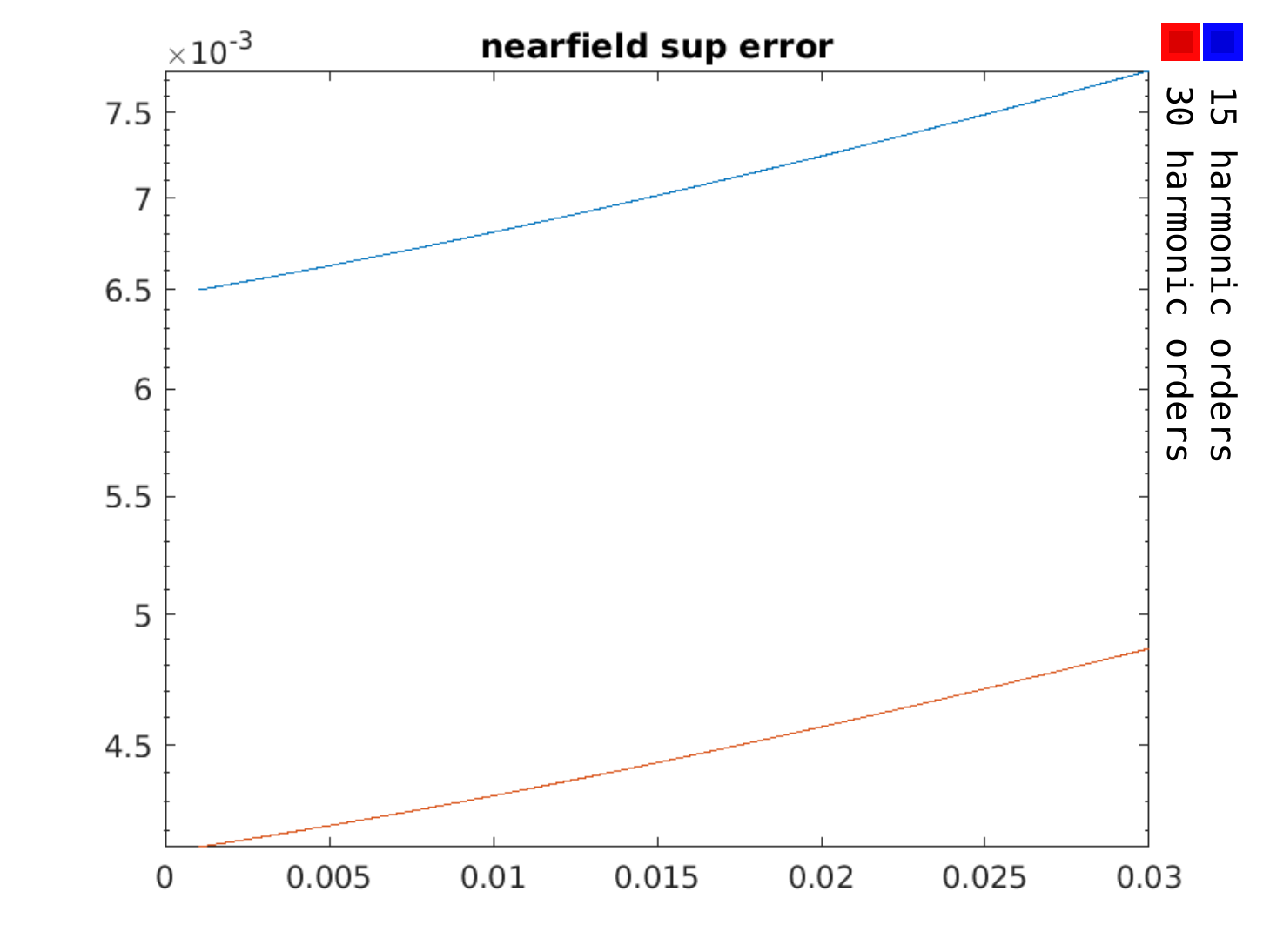}
		\label{chart:anr10hBothR001_nearsup}
			\vspace{-0.5cm}
		\caption{  Radii increments in $[0, 0.029]$}
	\end{subfigure}
	\begin{subfigure}{\figsizeB}
		\includegraphics[width=\textwidth]{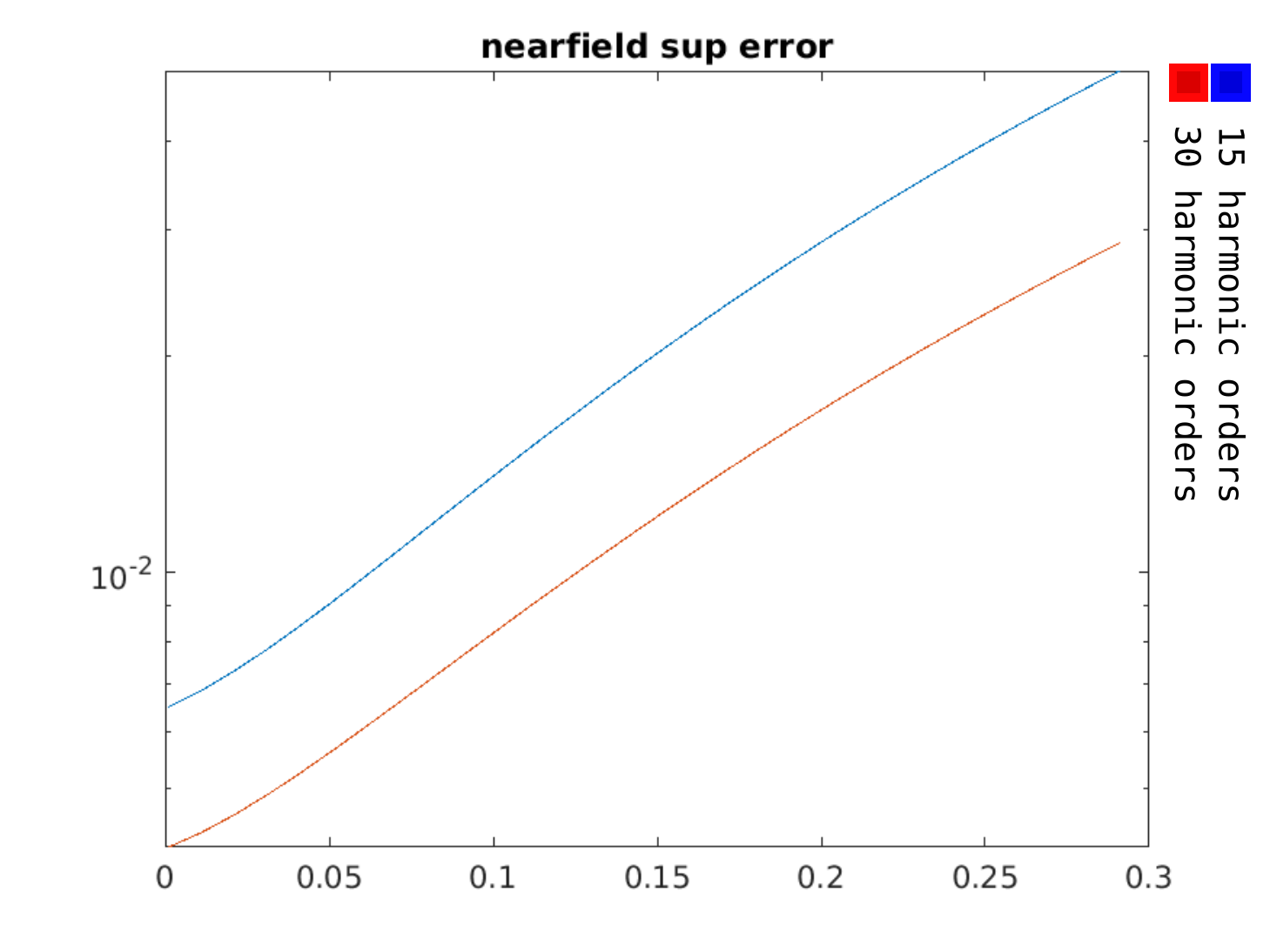}
		\label{chart:anr10hBothR01_nearsup}
		\vspace{-0.5cm}
	\caption{  Radii increments in $[0, 0.29]$}
	\end{subfigure}
	\caption{Relative supremum error as a function of $D_1$ radii increment. }
	\label{2.54}
\end{figure}

\begin{figure}[!htb] \centering
	\begin{subfigure}{\figsizeB}
		\includegraphics[width=\textwidth]{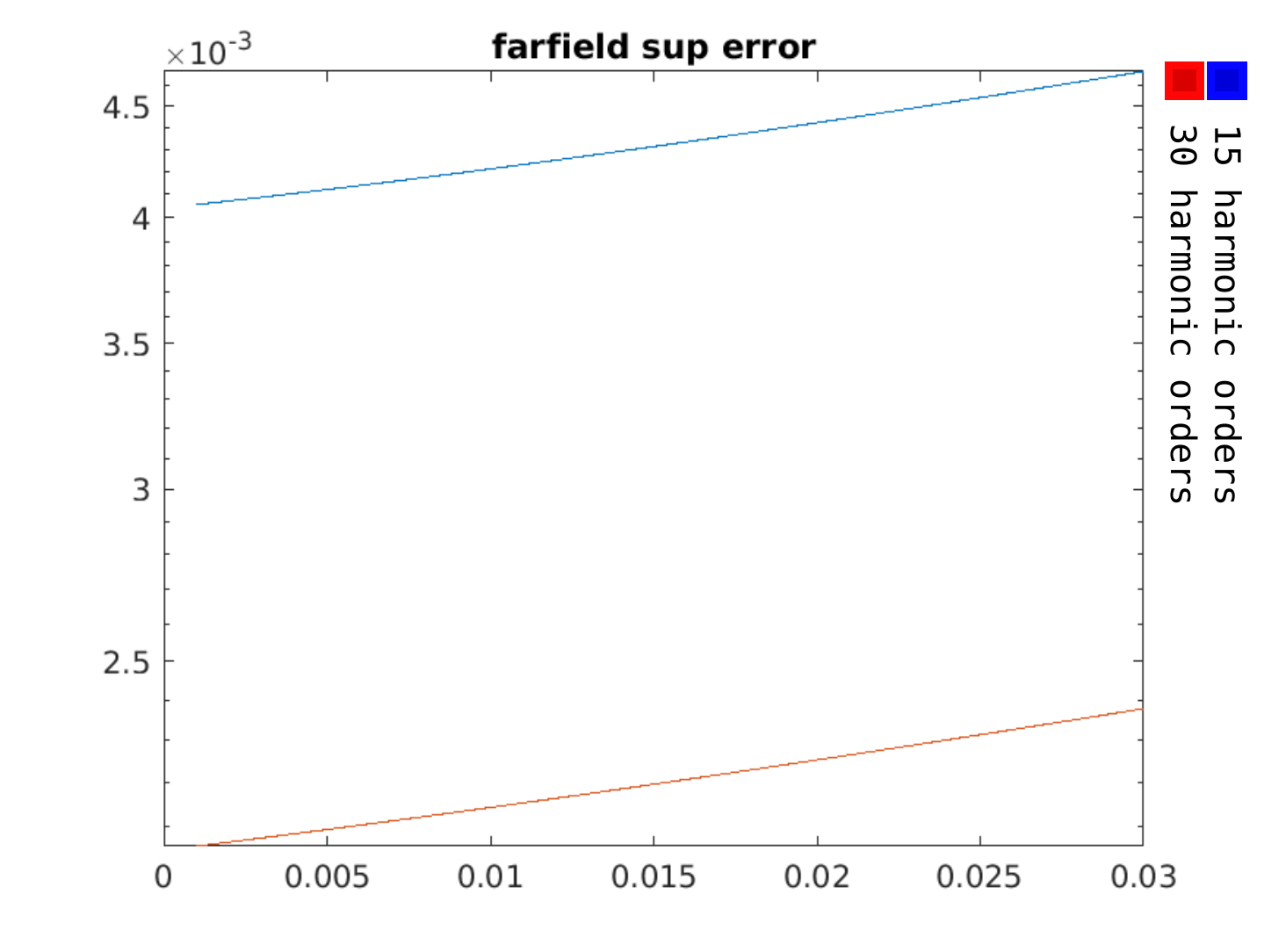}
		\label{chart:anr10hBothR001_farsup}
		\vspace{-0.5cm}
		\caption{  Radii increments in $[0, 0.029]$}
	\end{subfigure}
	\begin{subfigure}{\figsizeB}
		\includegraphics[width=\textwidth]{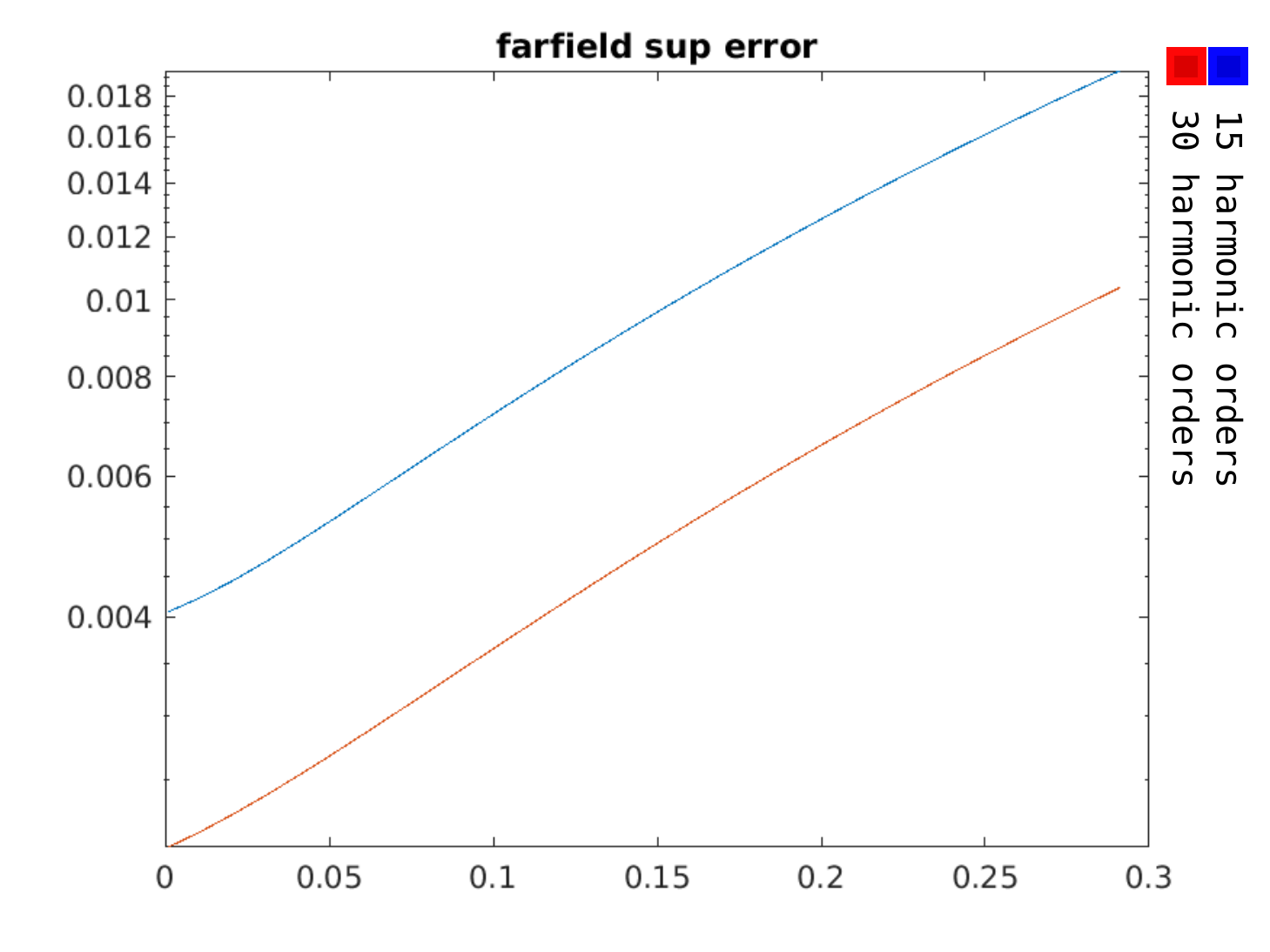}
		\label{chart:anr10hBothR01_farsup}
		\vspace{-0.5cm}
		\caption{  Radii increments in $[0, 0.29]$}
	\end{subfigure}
	\caption{Absolute supremum error on $\partial D_2$ as a function of $D_1$ radii increment.}
	\label{2.55}
\end{figure}

Again, the use of 15 harmonic orders produced less accurate results. In this regard, Figure \ref{2.54} shows the plots of the relative supremum errors in $D_1$. The relative supremum error in $D_1$ are of order $O(10^{-3})$ for the first increments and remains below $10^{-2}$ when the radii increments are less than $10$cm reaching undesirable levels for larger radii increments. Figure \ref{2.55} shows that over the entire range of increments used, the absolute supremum error in the far field stays of order $O(10^{-3})$ as long as the radii increments are less than $8$cm and reaches order $O(10^{-2})$ afterwards.
%

\begin{figure}[!htb] \centering
	\begin{subfigure}{\figsizeB}
		\includegraphics[width=\textwidth]{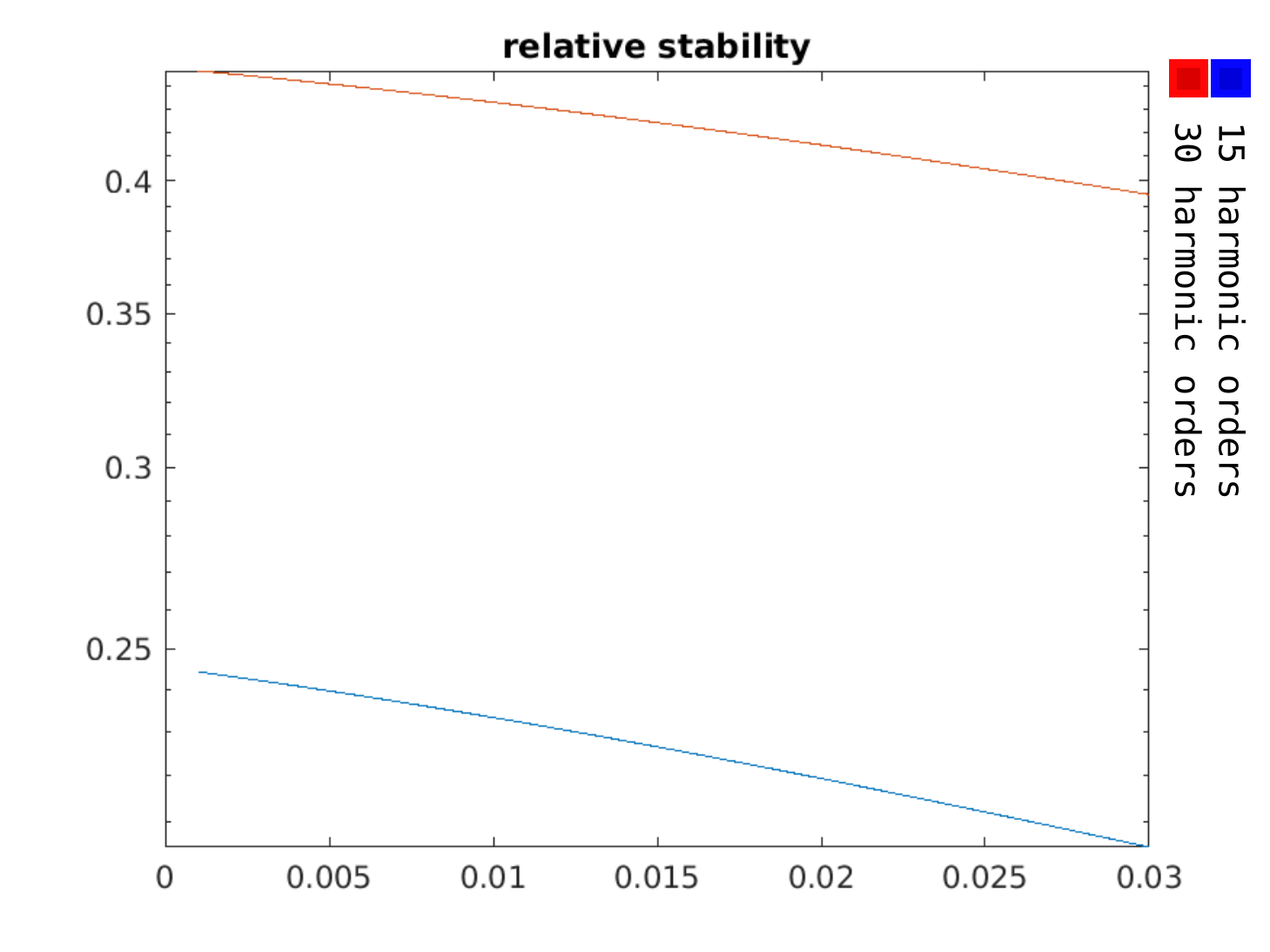}
		\label{chart:anr10hBothR001_relstab}
		\vspace{-0.5cm}
	\caption{  Radii increments in $[0, 0.029]$}
	\end{subfigure}
	\begin{subfigure}{\figsizeB}
		\includegraphics[width=\textwidth]{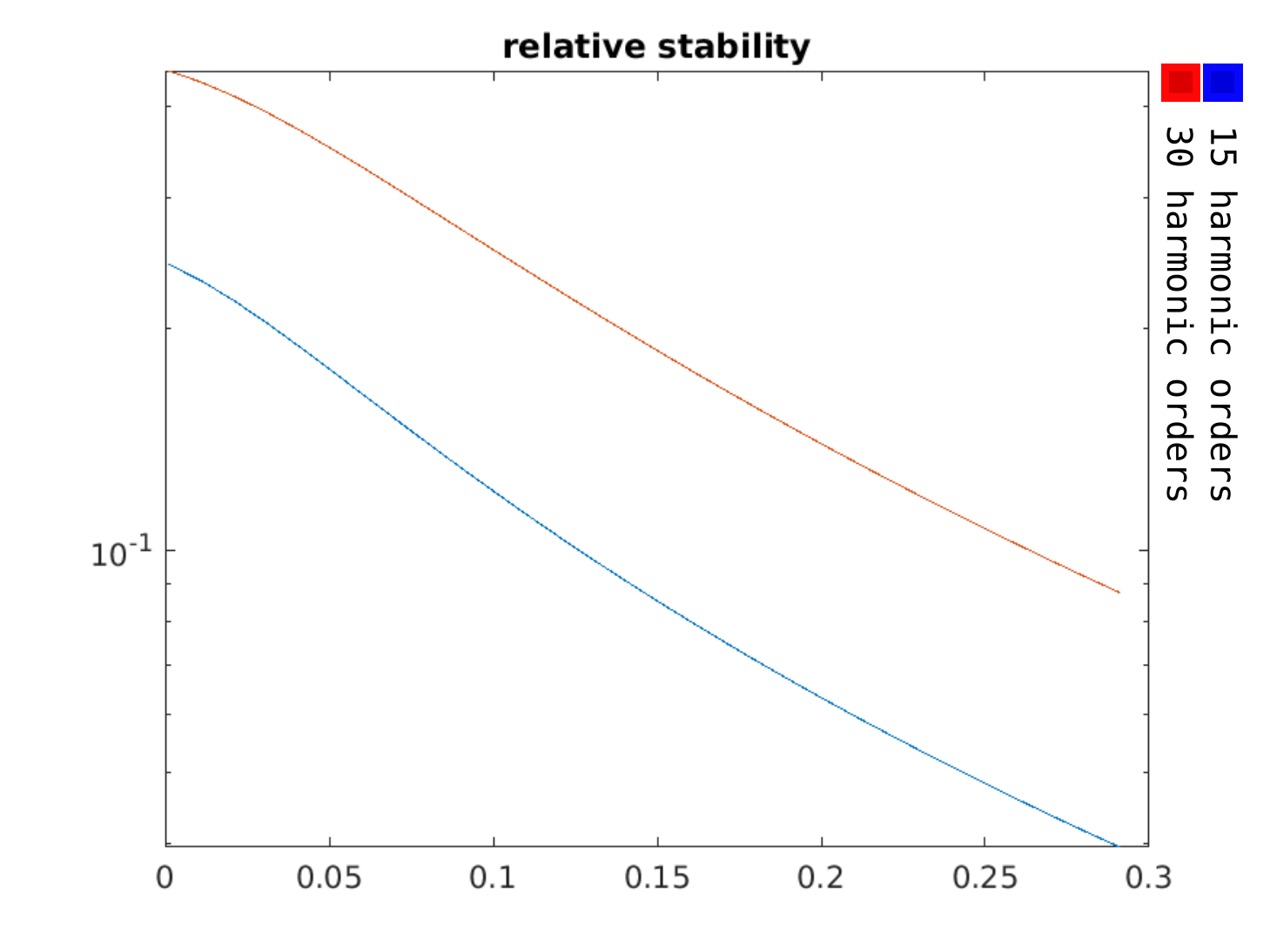}
		\label{chart:anr10hBothR01_relstab}
			\vspace{-0.5cm}
		\caption{  Radii increments in $[0, 0.29]$}
	\end{subfigure}
	\caption{Stability measure as a function of $D_1$ radii increment.}
	\label{2.56}
\end{figure}

The stability measure is shown in Figure \ref{2.56}. Notice that as the stability gets better with larger increments in the region size and distance from source.

\subsection{Multi region controls}
\label{SSnull}

\begin{figure}[!htb] \centering
	\includegraphics[scale=0.75]{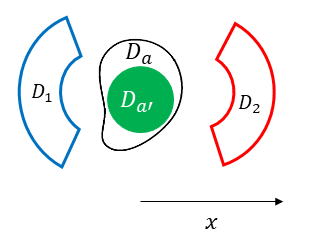}
	\caption{The initial geometry for the multi region controls.}
	\label{doubleregion}
\end{figure}

In the next set of tests, we consider problem \eqref{P1a}, \eqref{P1b} for the configuration described in (\ref{geom}-i) with a single fictitious source $D_{a'}=B_{0.01}({\mathbf{0}})$ and two compact control regions $D_1$ and $D_2$ (see Figure \ref{doubleregion}). For the initial geometry in our tests the primary control region $D_1$ is defined as before
\begin{equation}
\label{D1'}
D_1=\left \{(r,\theta, \phi) : r \in [0.011,0.015], \theta \in \left [-\frac{\pi}{4}, \frac{\pi}{4} \right ] , \phi \in \left [\frac{3\pi}{4}, \frac{5\pi}{4} \right] \right   \}
\end{equation} 
 while the secondary control region $D_2$ is given by \begin{equation}
\label{}
D_2\!=\!\left \{\!(r,\theta, \phi) : r \!\in \![0.011,0.015], \theta \in \left [-\frac{\pi}{4}, \frac{\pi}{4} \right ] , \phi \in \left [\frac{7\pi}{4}, 2\pi\right]\!\cup\!\left[0,\frac{\pi}{4} \right]  \!\right \} \!+ \!(0.09,0,0). 
\end{equation} 
As before, we mention that in our method only boundary controls are needed to imply smooth interior controls. The first question we address is the possibility to approximate the outgoing plane wave $u_1(\Bx)=e^{i\Bx \cdot (-10 \hat {\Be}_1)}$ in $D_1$  while imposing a small field in $D_2$. For this problem and the subsequent sensitivity tests the boundaries of the two regions are uniformly discretized with 2,400 points each.

\begin{figure}[!htb] \centering
	\!\!\!\!\vspace{-1cm}\begin{subfigure}{\figsizeC}
		\includegraphics[width=1.4\textwidth]{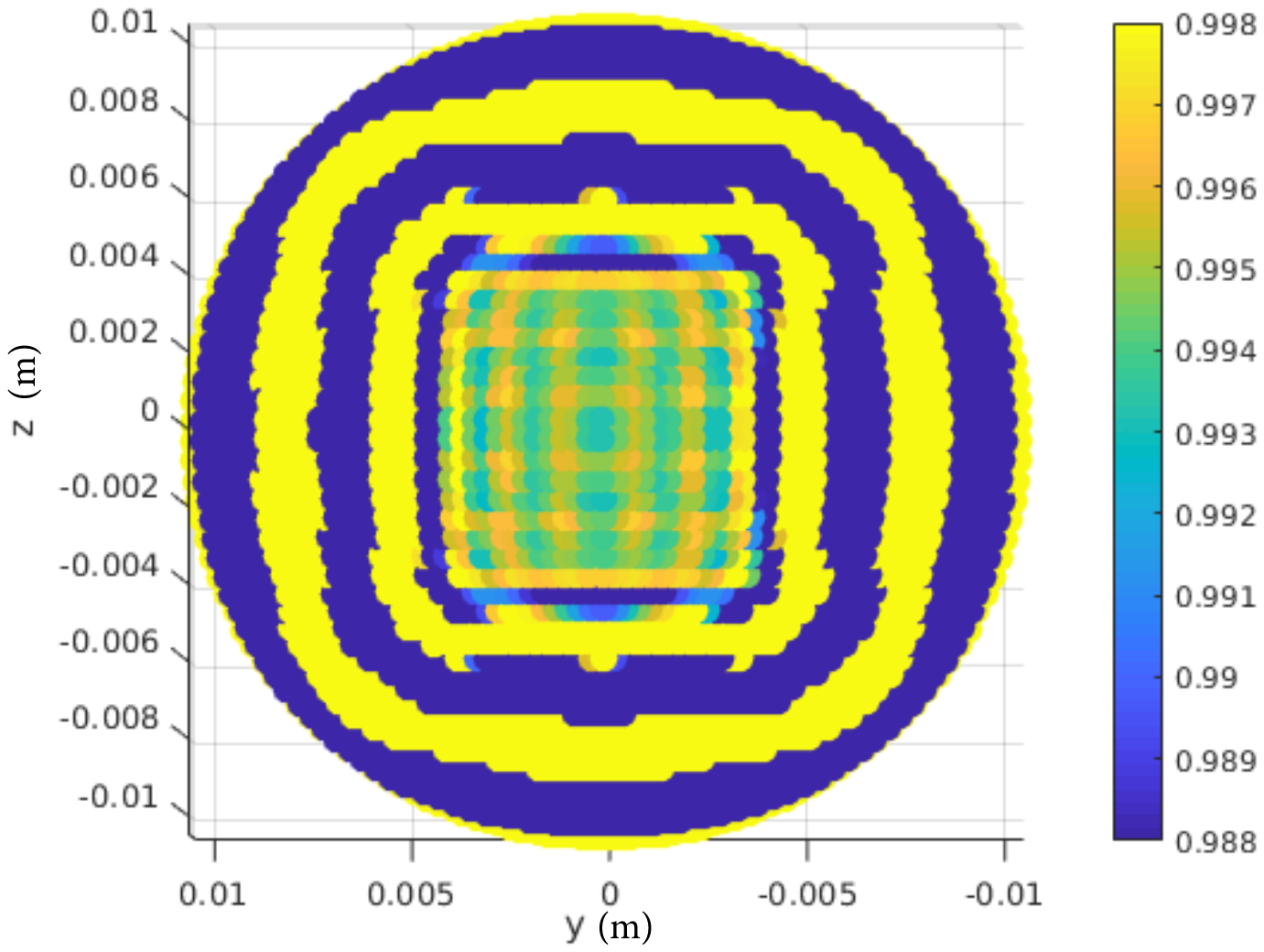}
		\label{kw_double_front}
		\vspace{-4cm}
		\caption{front}
	\end{subfigure}
	\;\;\;\;\;\;\;\;\;\;\;\;\;\;\begin{subfigure}{\figsizeC}
		\includegraphics[width=\textwidth]{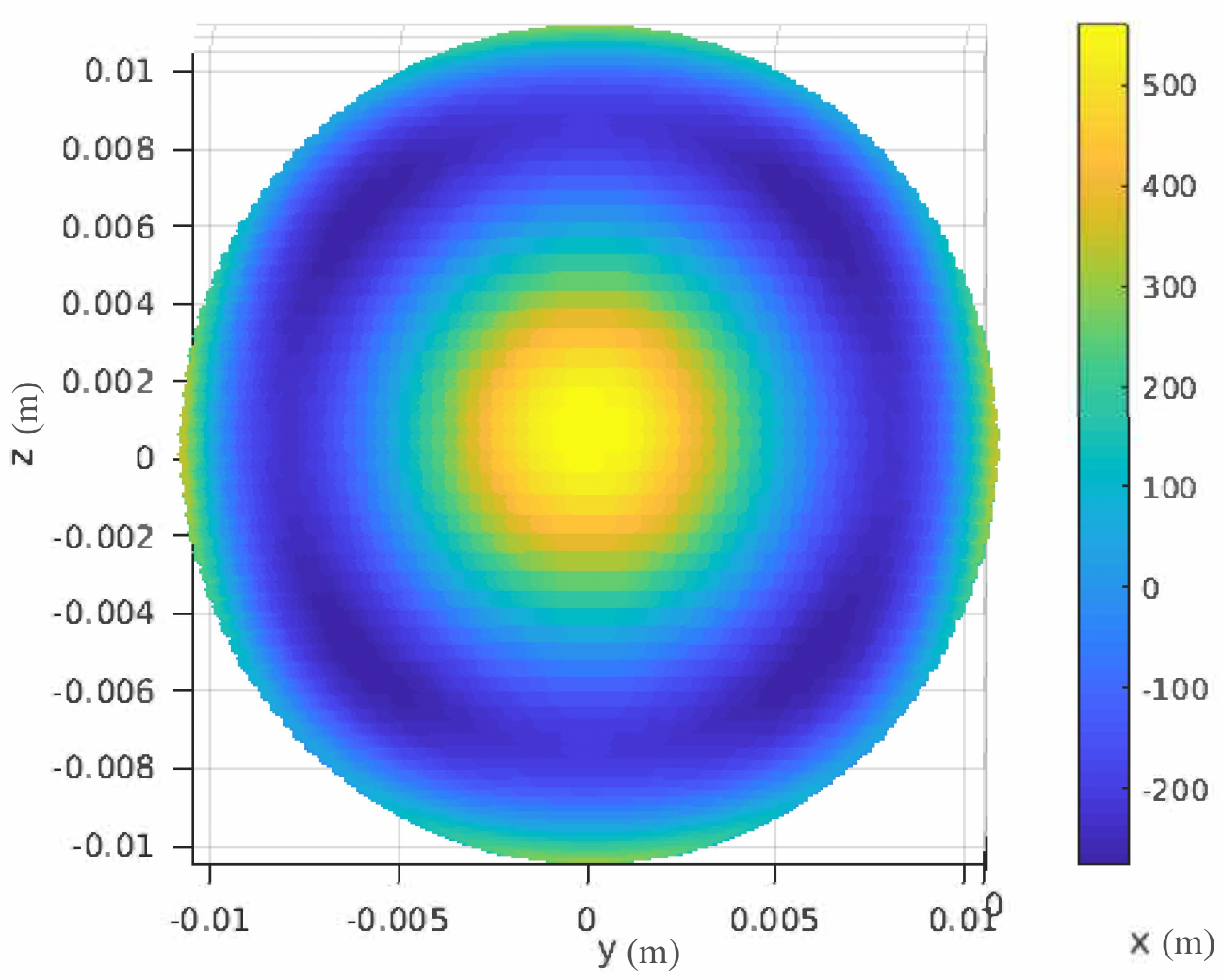}
		\label{kw_double_front}
		\caption{back}
	\end{subfigure}
	
	\begin{subfigure}{\figsizeB}
		\includegraphics[width=\textwidth]{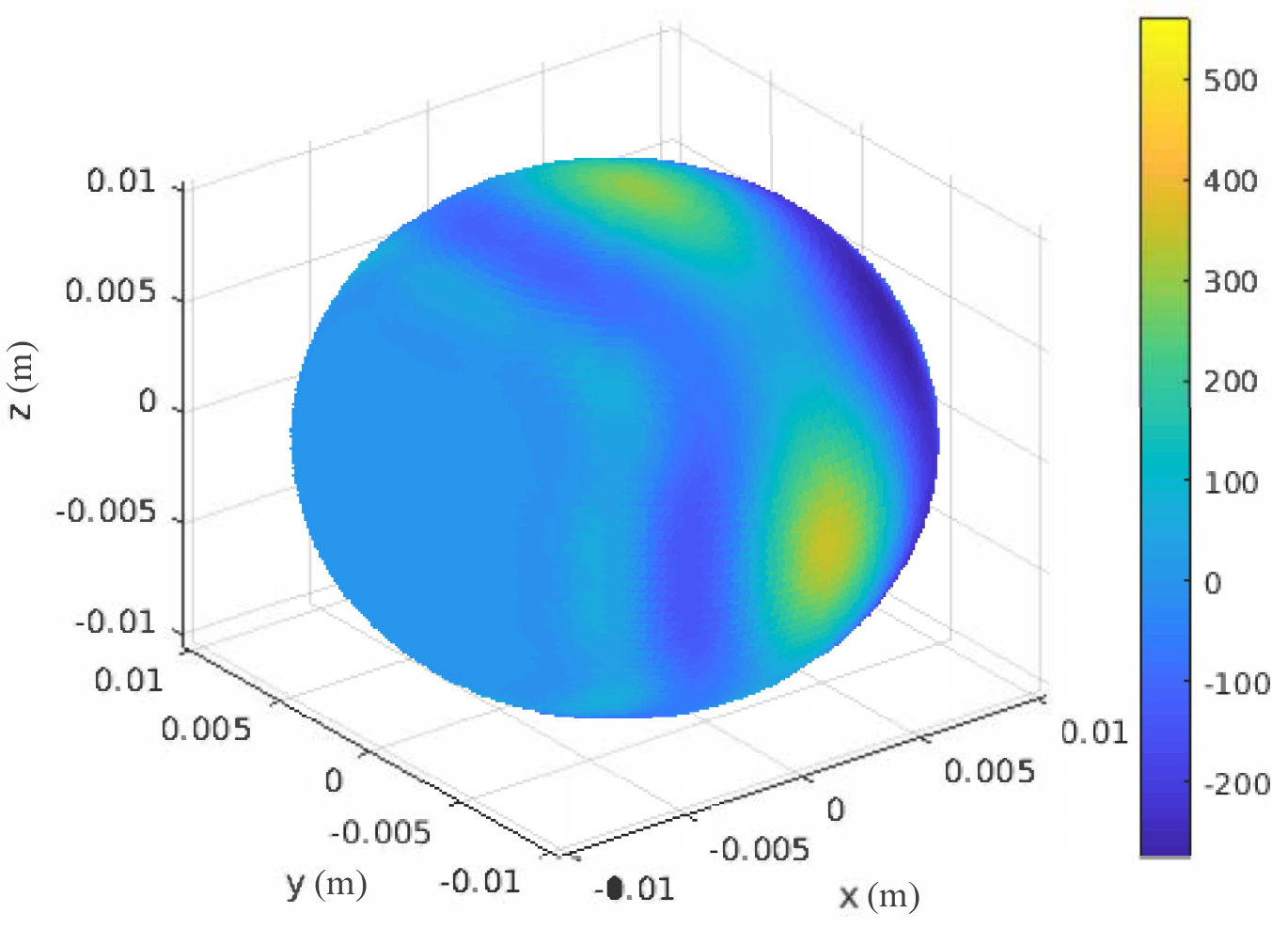}
		\label{kw_double_angle}
		\caption{side}
	\end{subfigure}
	
	\caption{Different views of the surface field pattern on $\partial B_{0.0105}$.}
	\label{kw_double}
\end{figure}

As in Remark \ref{rem1}, we recall that the actual physical source boundary is only constrained by the condition to embed the fictitious source and to have an empty intersection with $D_1\cup D_2$ and that, as it can be seen in \eqref{eqnpb}, the required boundary input on the physical source $D_a$ will be ${\cal D}(w_\alpha)$ restricted to $\partial D_a$. Thus, if for instance the real source $D_a$ is chosen to be the sphere $B_{0.0105}(\mathbf 0)$ then Figure \ref{kw_double} shows the necessary pattern required on $\partial D_a$ for a good control in $D_1$ and $D_2$. As observed before for the configuration (\ref{geom}-ii) discussed in Section \ref{SSANR}, we can see in Figure \ref{kw_double} that  the magnitude of the required boundary input on the real source $\partial B_{0.0105}(\mathbf 0)$ is a few orders smaller then the magnitude of the density $w_\alpha$ (which is $O(10^{5})$ as shown in Figure \ref{rotantnorm} below) and has less complexity in the pattern. This once more supports our claim that one could search and find an actual source boundary around the fictitious domain $D_{a'}$ so that while it does not intersect  $D_1\cup D_2$ it is such that the boundary input required on it needs less power for its instantiation and presents a lower level of complexity in its pattern. 

\begin{figure}[!htb] \centering
	\includegraphics[scale=0.7]{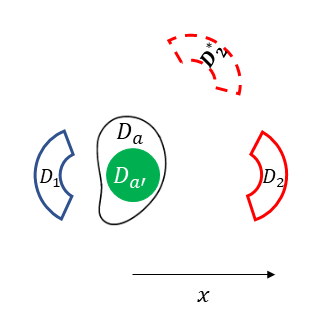}
	\caption{Illustration of two iterations showing the secondary region rotated around $D_a$.}
	\label{varyingposition}
\end{figure}

Next we will present the sensitivity of our scheme with respect to relative positions between the two control regions. In this regard, Figure \ref{varyingposition} shows the geometry of two iterations. In each iteration, $D_1$ is fixed while we vary the position of $D_2$ via counterclockwise rotation about $D_a$. 

\begin{figure}[!htb] \centering
	\begin{subfigure}{\figsizeB}
		\includegraphics[width=\textwidth]{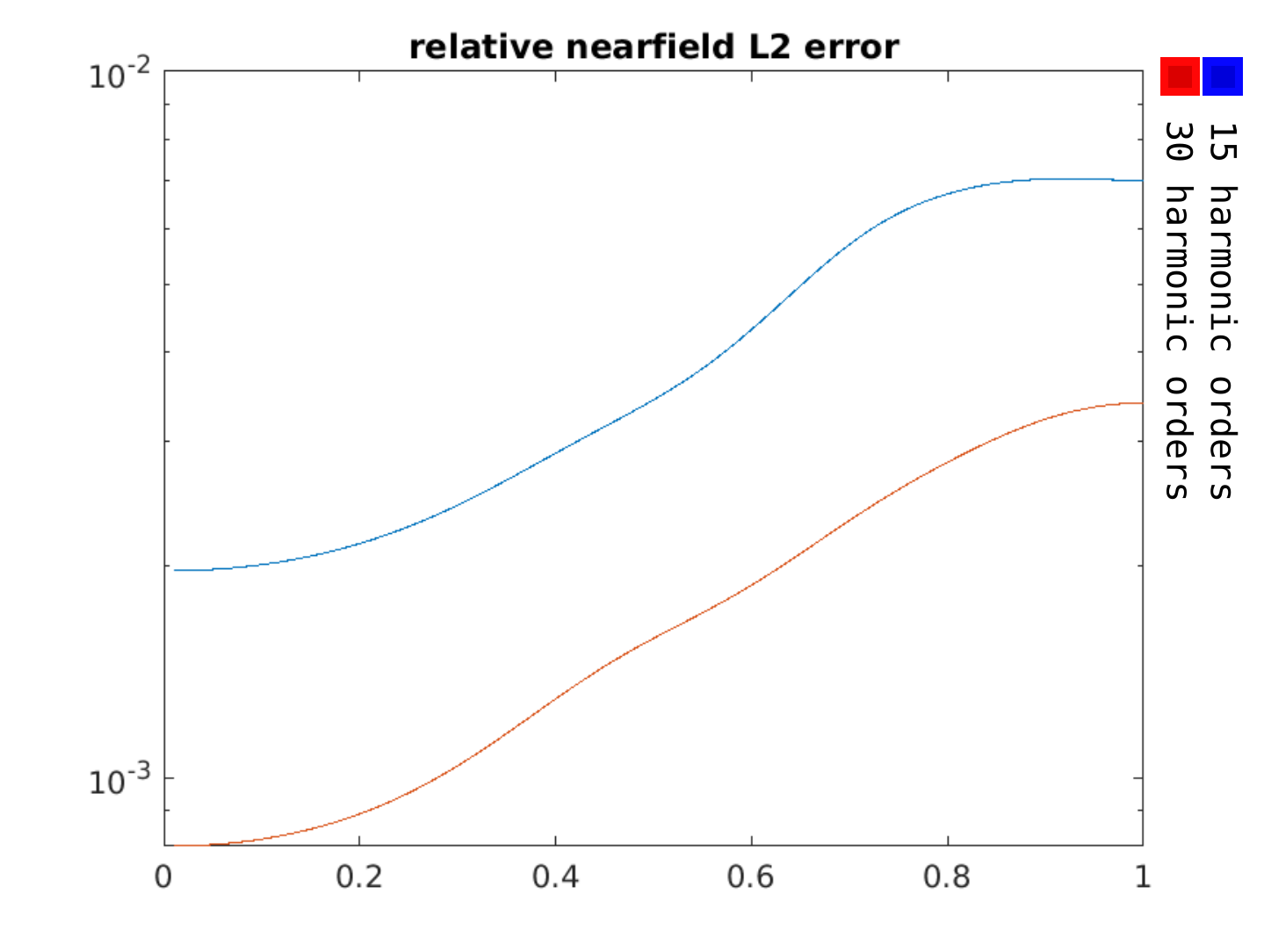}
		\label{chart:anr10hBothR001_rotnearrel}
		\vspace{-0.5cm}
		\caption{Relative $L^2$ error in $D_1$}
	\end{subfigure}
	\begin{subfigure}{\figsizeB}
		\includegraphics[width=\textwidth]{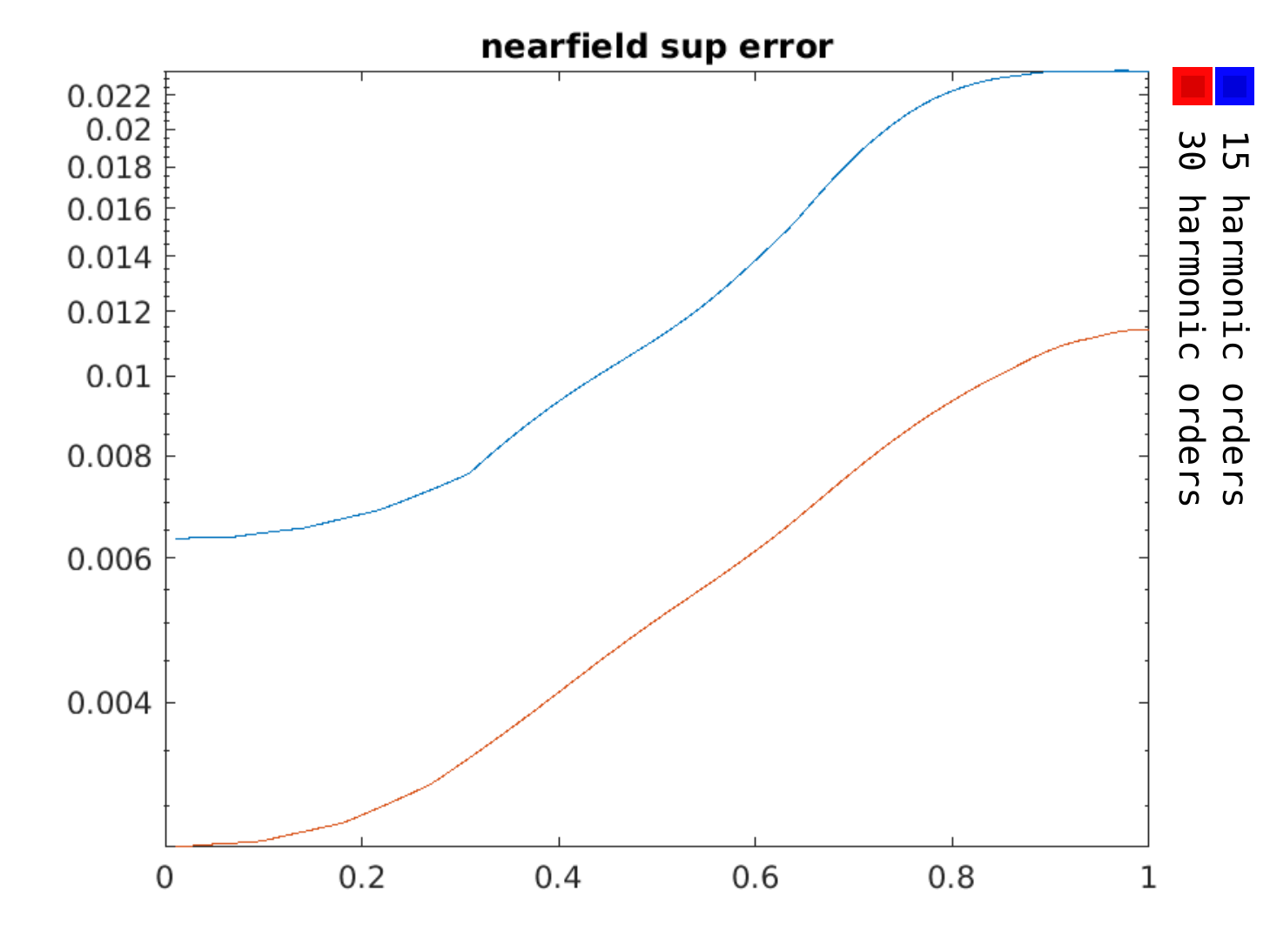}
		\label{chart:anr10hBothR01_rotnearsup}
		\vspace{-0.5cm}
		\caption{Relative sup error in $D_1$}
	\end{subfigure}
	
	\caption{Accuracy errors in $D_1$}
	\label{rotnearfielderr}
\end{figure}

\begin{figure}[!htb] \centering
	\includegraphics[scale=0.5]{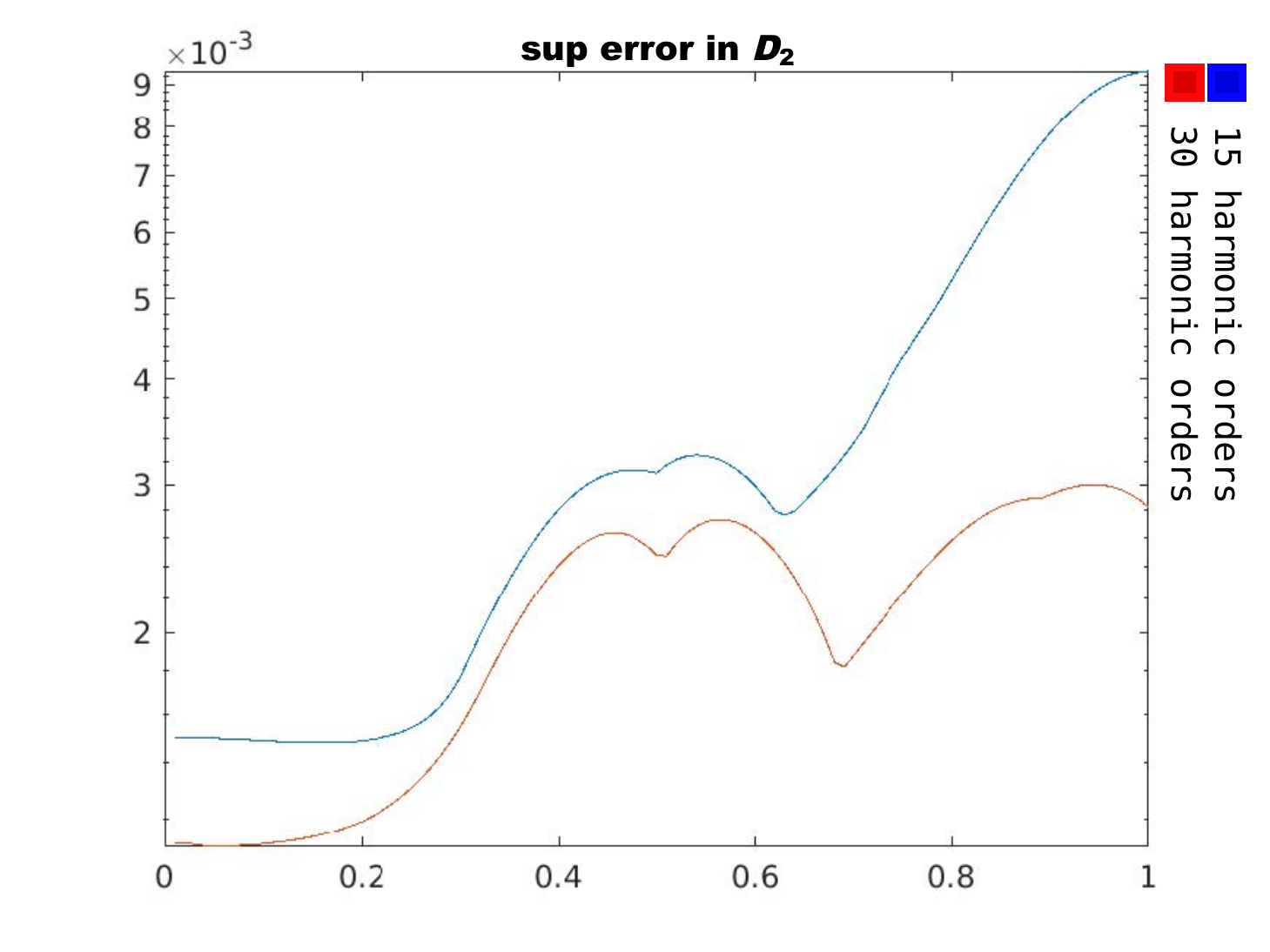}
	\label{chart:farsup}
	\caption{Supremum error in $D_2$ as a function of its angle of rotation.}
	\label{rotfarfielderr}
\end{figure}

 Figure \ref{rotantnorm} shows the $L^2$ - norm of the density $w_\alpha $, as a function of the angle of rotation. We can observe that its levels remain of order $O(10^6)$ for the entire rotation spectrum.
\begin{figure}[!htb] \centering
	\includegraphics[scale=0.52]{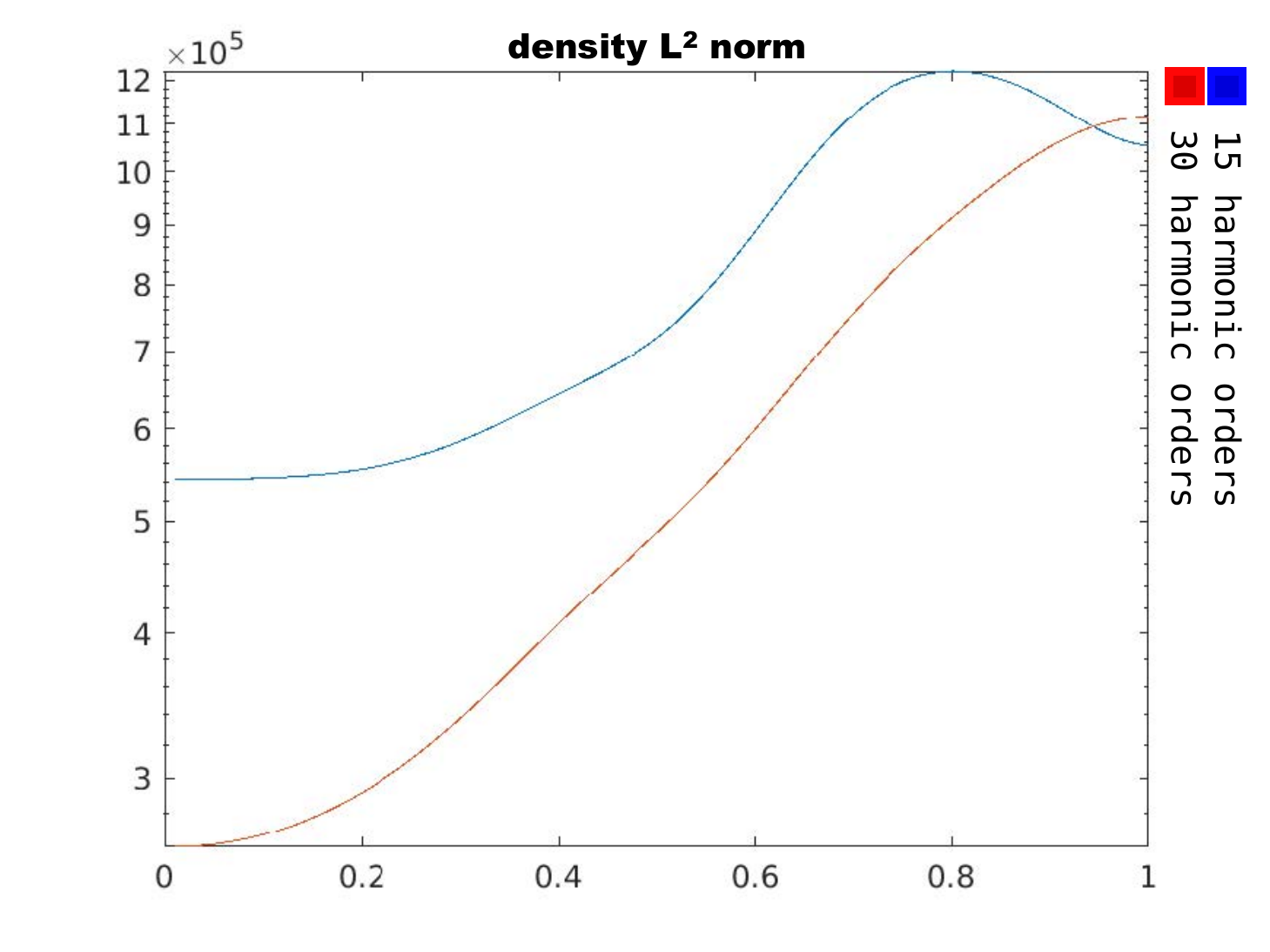}
	\caption{$L^2$ norm of $w_\alpha$ as a function of the secondary region's angle of rotation.}
	\label{rotantnorm}
\end{figure}
Figure \ref{rotnearfielderr} show the relative  $L^2$ and respectively relative supremum  error in region $D_1$. The performance is good with largest relative error for a 180 degree rotation but still of order $O(10^{-2})$. 
Figure \ref{rotfarfielderr} shows a small field generated in region $D_2$ as desired. One can see that the contrast between the field in region $D_1$ and the small field in region $D_2$ remains between approximately $40dB$ and $60dB$ with better performance for the case with more harmonic orders (where we use the $20\log_{10}$ convention for computing the decibel (dB) level).

\begin{figure}[!htb] \centering
	\includegraphics[scale=0.5]{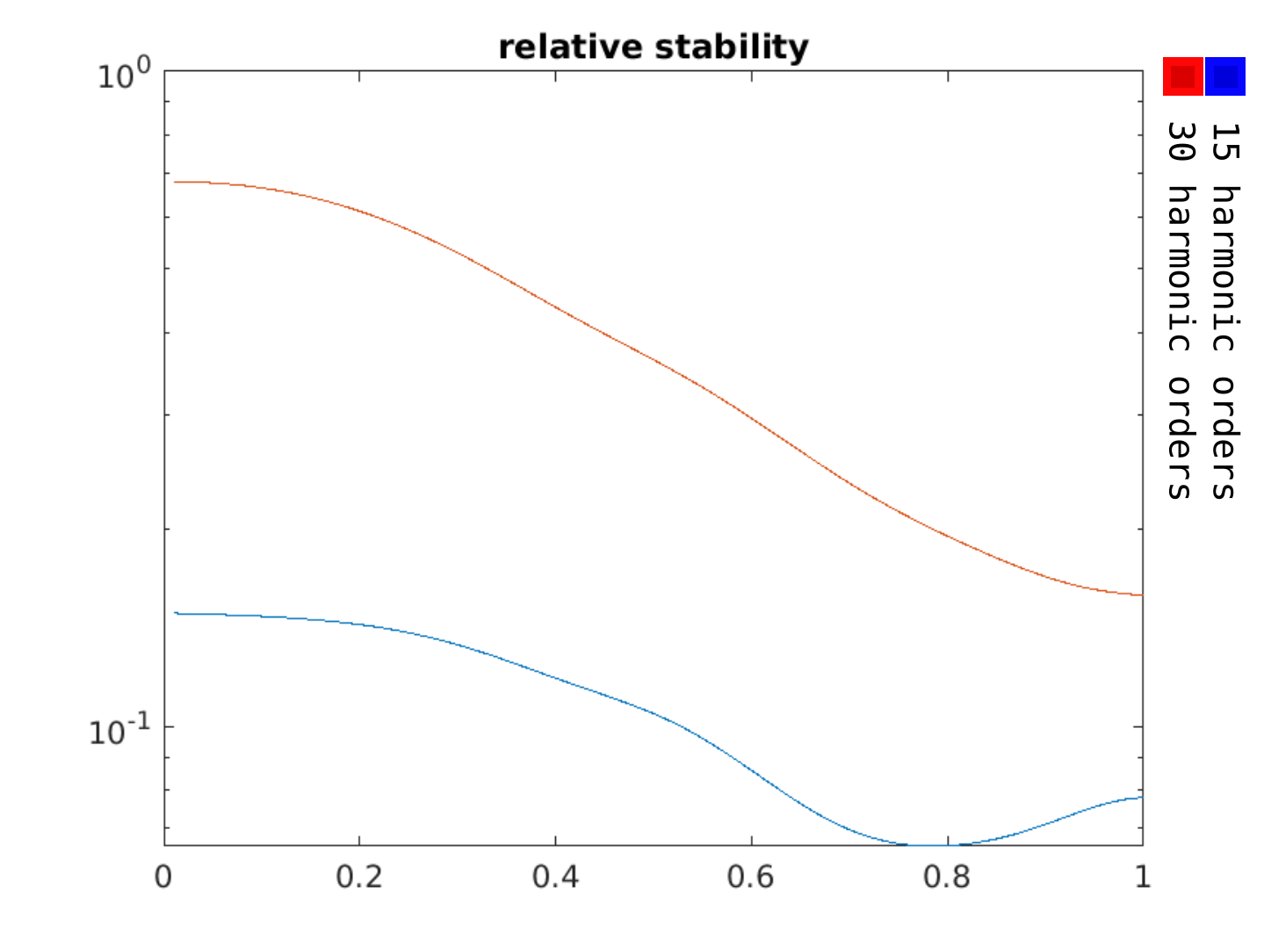}
	\caption{Stability measure as a function of the secondary region's angle of rotation.}
	\label{rotstab}
\end{figure}

The stability measure is plotted in Figure \ref{rotstab}. It tends to decrease as the angle of rotation approaches $\pi$ where it is  below $10^{-1}$. 

\begin{figure}[!htb] \centering
	\includegraphics[scale=0.8]{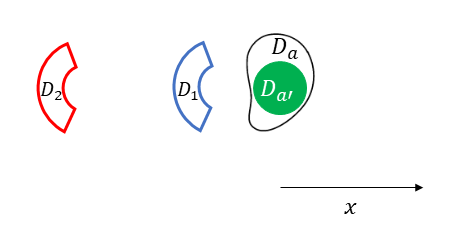}
	\vspace{-0.35cm}
	\caption{Sketch of the geometry where $D_1$ is acting as a near field obstacle to $D_2$.}
	\label{obstacleregion}
\end{figure}

We continue our sensitivity analysis with another important test. That is, we study the accuracy of our scheme in the context of problem \eqref{P1a}, \eqref{P1b} for the same configuration described at (\ref{geom}-i) with one fictitious source $D_{a'}=B_{0.01}({\mathbf{0}})$ and two compact control regions $D_1$ and $D_2$ described below and sketched in Figure \ref{obstacleregion}, 
\begin{equation}
\label{}
D_1=\left \{(r,\theta, \phi) : r \in [0.011, 0.015], \theta \in \left [-\frac{\pi}{8}, \frac{\pi}{8} \right ] , \phi \in \left [\frac{3\pi}{4}, \frac{5\pi}{4} \right] \right   \} +(-0.03,0,0)
\end{equation} 
and
\begin{equation}
\label{}
D_2=\left \{(r,\theta, \phi) : r \in [0.011,0.015], \theta \in \left [-\frac{\pi}{8}, \frac{\pi}{8} \right ] , \phi \in \left [\frac{3\pi}{4}, \frac{5\pi}{4} \right]\right \} + (-0.1,0,0)
\end{equation}
The novelty of the test consist in the fact that now we consider $u_1=0$ and $u_2(\Bx)=e^{i\Bx \cdot (-10 \hat {\Be}_1)}$ in problem \eqref{P1a}, \eqref{P1b}, i.e., we approximate an outgoing plane wave in region $D_2$ while having a very small field in the near field region $D_1$. This geometry is relevant for applications where the objective is to focus energy or communicate behind a near field obstacle. For this simulation we used 2400 control points on $\partial D_1$ and 38400 control points on $\partial D_2$. 

\begin{figure}[!htb] \vspace{-0.5cm}\centering
	\begin{subfigure}{\figsizeB}
		\includegraphics[width=\textwidth]{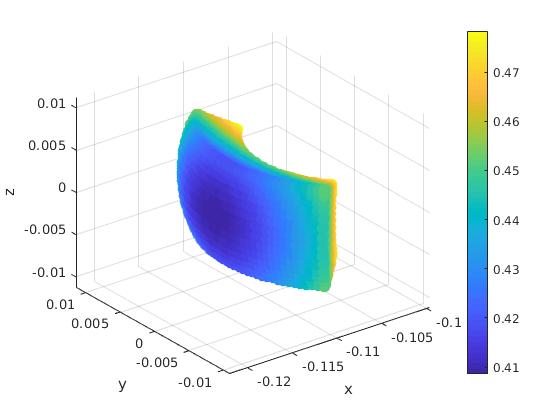}
		\vspace{-0.5cm}
		\caption{Field desired in $D_2$}
	\end{subfigure}
	\begin{subfigure}{\figsizeB}
		\includegraphics[width=\textwidth]{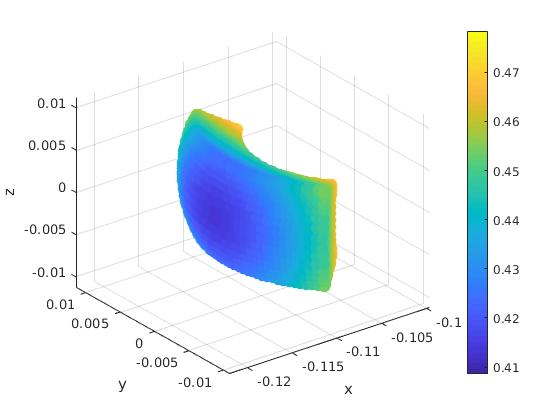}
		\vspace{-0.5cm}
		\caption{Field generated in $D_2$}
	\end{subfigure}
	
	\begin{subfigure}{\figsizeB}
		\vspace{0cm}
		\includegraphics[width=\textwidth]{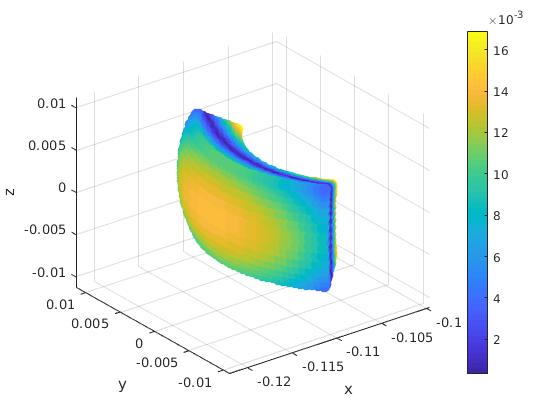}
		\vspace{-0.5cm}
		\caption{Accuracy error in $D_2$}
	\end{subfigure}
\begin{subfigure}{\figsizeB}
	\vspace{0cm}
	\includegraphics[width=\textwidth]{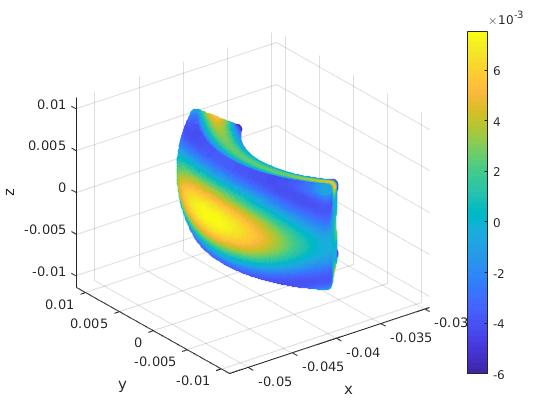}
	\vspace{-0.5cm}
	\caption{Field generated in $D_1$}
\end{subfigure}
\caption{Performance of the scheme in generating a null in $D_1$ and a plane wave in $D_2$.}
	\label{double_obs}
\end{figure}

The overall performance of our scheme is shown in Figure \ref{double_obs}. We can see the accuracy of the control in region $D_1$ with relative supremum error of order $O(10^{-2})$ and with small fields in region $D_2$ as desired. The contrast between the quiet region $D_2$ and the other control region $D_1$ is in this case over $40dB$. Our numerics suggest that more harmonics and more control points in the two control regions can allow us to extend these results to the situation when region $D_2$ is larger and when region $D_1$ is much closer to the source.

\begin{figure}[!htb] \centering
	\!\!\!\!\vspace{0cm}\begin{subfigure}{\figsizeC}
		\includegraphics[width=\textwidth]{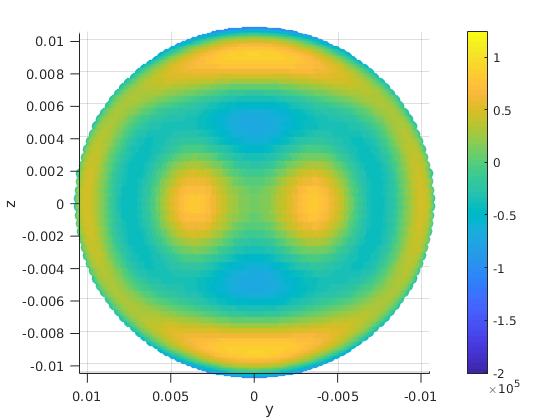}
		\label{kw_double_front}
		\vspace{0cm}
		\caption{front}
	\end{subfigure}
	\;\;\;\;\;\;\;\;\;\;\;\;\;\;\begin{subfigure}{\figsizeC}
		\includegraphics[width=\textwidth]{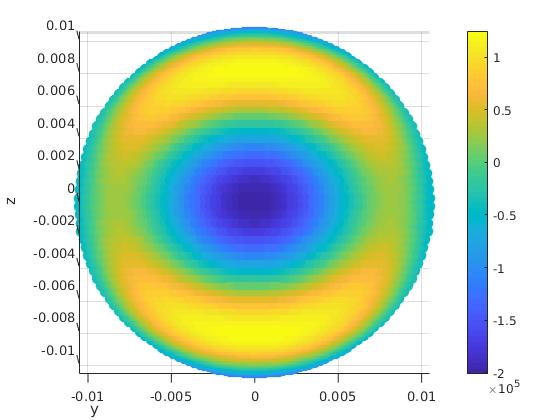}
		\label{kw_double_front}
		\caption{back}
	\end{subfigure}
	
	\begin{subfigure}{\figsizeB}
		\includegraphics[width=\textwidth]{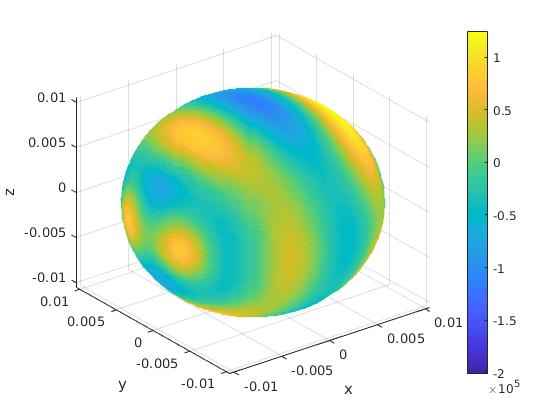}
		\label{kw_double_angle}
		\caption{side}
	\end{subfigure}
	
	\caption{Different views of the surface field pattern on $\partial B_{0.0105}({\mathbf{0}})$}
	\label{obs-pattern}
\end{figure}

Figure \ref{obs-pattern} presents the required pattern on the actual source $D_a$ in the particular example when $D_a=B_{0.0105}({\mathbf{0}})$. The front picture refers to the part of the surface facing the control regions, the back picture shows the opposite part of the surface and the side plot presents a panoramic view of the surface pattern. We can see that the required surface pattern is much simpler and without a large variation in the field oscillation as it was the case for the geometry (\ref{geom}-ii) considered in Section \ref{SSANR}. 

Based on the previous analysis, by superposition, we could in principle predict a good accuracy for our scheme for the case of approximating multiple prescribed Helmholtz potentials in given mutually disjoint exterior regions.

\subsection{Fourier synthesis} 
\label{fourier}
In this section we consider again problem \eqref{P1a}, \eqref{P1b} for the configuration described in (\ref{geom}-i) with a single fictitious source $D_{a'}=B_{0.01}({\mathbf{0}})$ and two compact control regions $D_1$ and $D_2$ (see Figure \ref{doubleregion}) and perform a Fourier synthesis test to better understand the sensitivity of our scheme with respect to the frequency change. Thus, we present an example of a time domain synthesis result obtained by superposition of inverse problem solutions for different wave numbers. Each individual inverse problem solution was obtained using the algorithm described in Section \ref{numframework}. These solutions were collected and used to determine the required Fourier coefficients to match a prescribed time domain pattern.

In the numerical simulations below, 21 equally spaced wave numbers ranging from 5 to 15 were used along with 30 harmonic orders. The primary control region is 
\begin{equation}
\label{D1_fourier}
D_1=\left \{(r,\theta, \phi) : r \in [0.011,0.015], \theta \in \left [-\frac{\pi}{4}, \frac{\pi}{4} \right ] , \phi \in \left [\frac{3\pi}{4}, \frac{5\pi}{4} \right] \right   \},
\end{equation} 
while the secondary region is
\begin{equation}
\label{D2_secondary_fourier}
D_2\!=\!\left \{\!(r,\theta, \phi) : r \!\in \![0.011,0.015], \theta \in \left [-\frac{\pi}{4}, \frac{\pi}{4} \right ] , \phi \in \left [\frac{7\pi}{4}, 2\pi\right]\!\cup\!\left[0,\frac{\pi}{4} \right]  \!\right \} \!+ \!(0.09,0,0). 
\end{equation}
Here we approximate on $D_1$ a superposition of plane waves given by  \begin{equation} \label{plane_wave_fourier}  u_1(\Bx,t)=\sum_{\ell=10}^{30} \dfrac{2}{\ell} \exp (i\Bx  \cdot (k_\ell \hat {\Be}_1) ) \exp (-ik_\ell c t),\end{equation} where $k_\ell$ are the wave numbers $\frac{\ell}{2}$,  while maintaining a small field in region $D_2$. Again, the boundaries of the two regions are each uniformly discretized with 2,400 points.

By using the superposition principle we produced an animation found at \href{https://drive.google.com/open?id=1m5nnofuIc56_aaU76HHjqkH4aRLRIsRi}{animation} \cite{animationfull} which shows a one-period propagation of the waves on a cross section of the control regions. The animation were generated using a 2000-point discretization of the angular time. It offers two columns, the left column describing the time propagation of the field generated by our scheme in $D_1$ at the top versus the time propagation of the desired field $u_1$ described at \eqref{plane_wave_fourier} at the bottom. The right column of the animation shows, during the same time interval, the very small field values in $D_2$. 
	
\begin{figure}[!h]
	\centering
	\includegraphics[scale=0.4]{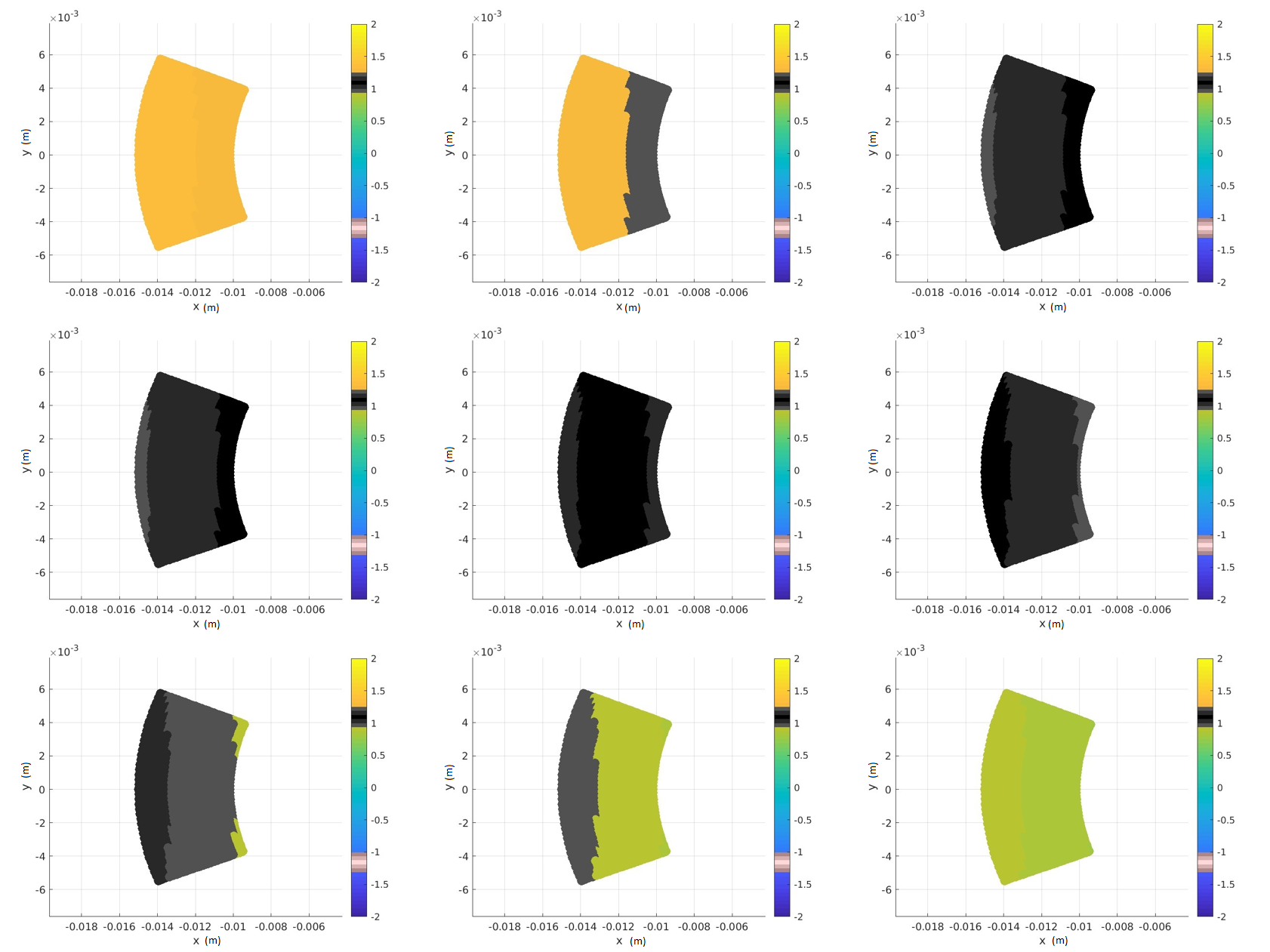}
	\caption{Time snapshots of a cross section of the near field at $kct$ values (left-right, top-down)$ \frac{37\pi}{2000},\frac{38\pi}{2000},\frac{39\pi}{2000}, \frac{40\pi}{2000},\frac{41\pi}{2000},...,\frac{45\pi}{2000}$.}
	\label{timesnapshots}
\end{figure}

Figure \ref{timesnapshots} shows snapshots of the cross section of the primary control region $D_1$ for various times. The color map on the panels for $D_1$ is slightly modified for the intervals $[-1.25, -1]$ and $[1, 1.25]$ to highlight the propagation of the waves. For better clarity we chose to highlight the contour line corresponding to a field value of around $1.15$ with a black stripe and the the contour line corresponding to a field value of around $-1.15$ with a white stripe. It can be observed that all throughout the period, there is a good match on the primary control. At the same time, as it can be seen from Figure \ref{timeerror} below, the field on the secondary control is maintained close to zero.

%

 Figure \ref{timeerror} shows the averaged relative error on the primary control (a) and the averaged absolute error on the secondary control (b) over the entire time period. Notice that the time-averaged relative error in region $D_1$ is of magnitude $10^{-3}$  while the absolute error on region $D_2$ stays of order $10^{-4}$ showing the accuracy of our scheme over a range of frequencies.
\begin{figure}[!htb] \centering
	\vspace{0cm}
	\begin{subfigure}{\figsizeB}
		\includegraphics[width=\textwidth]{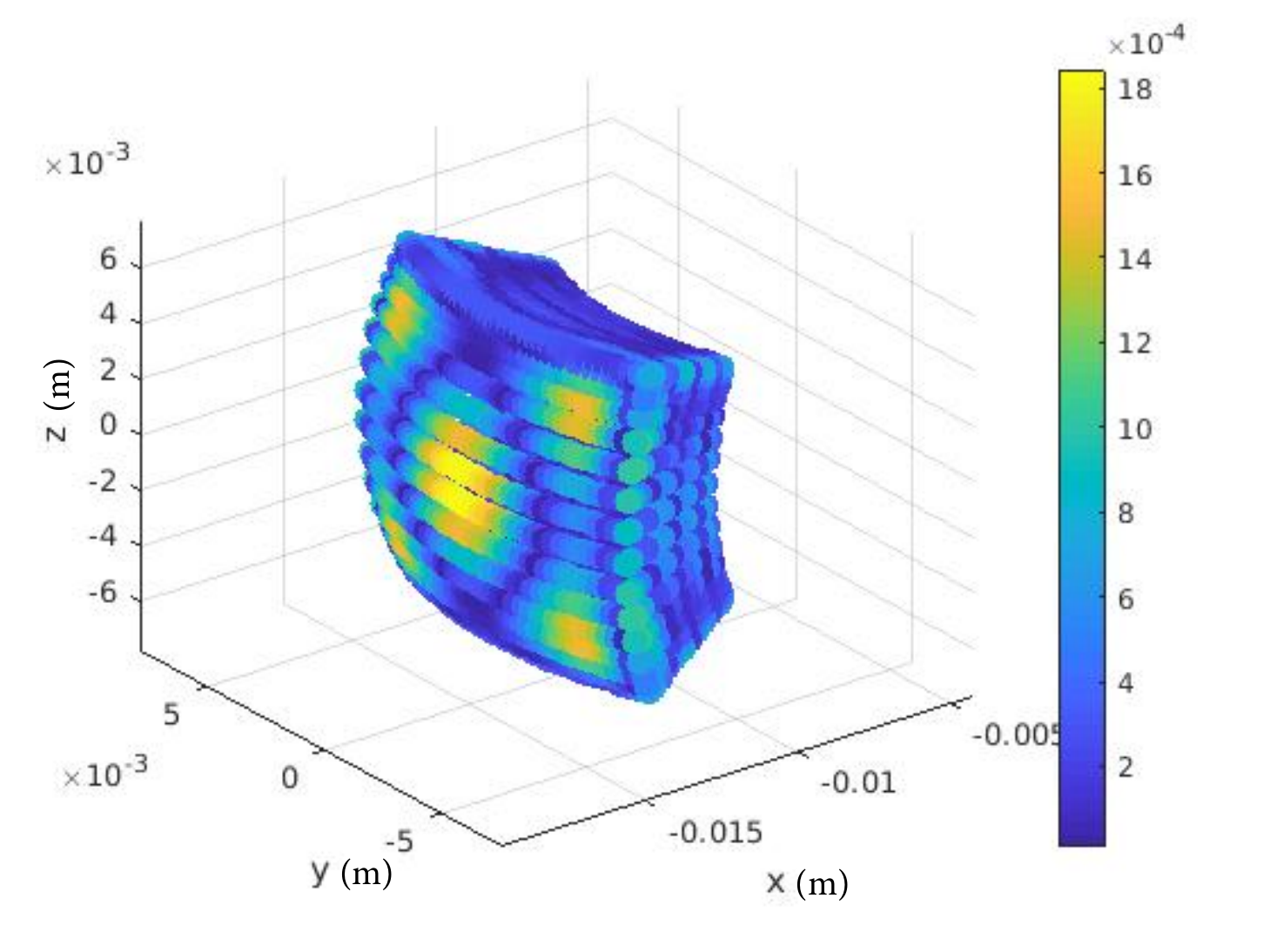}
		\caption{}
		\label{avg_e1.jpg}
	\end{subfigure}
	\begin{subfigure}{\figsizeB}
		\includegraphics[width=\textwidth]{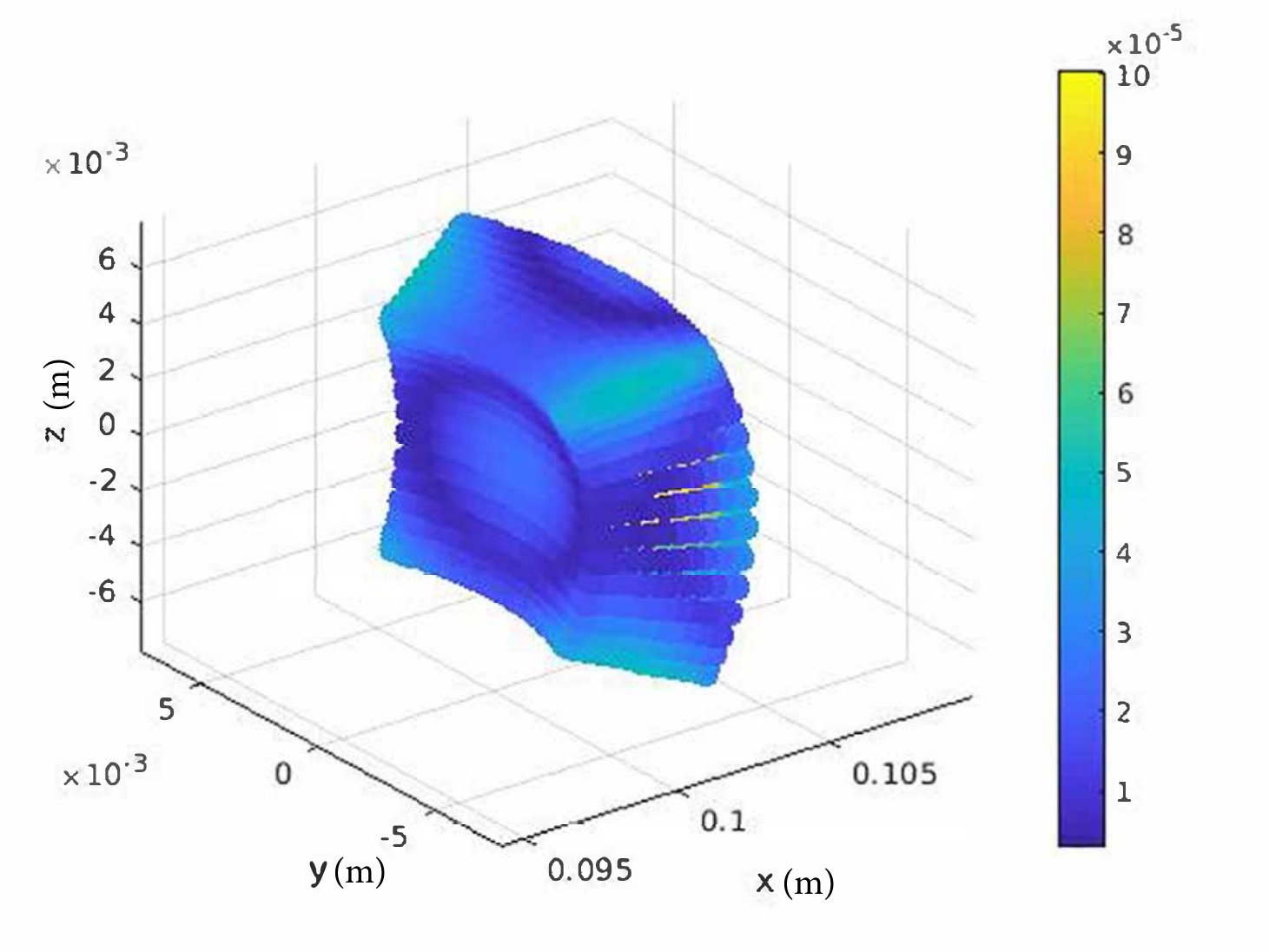}
		\caption{}
		\label{avg_e1.jpg}
	\end{subfigure}
	\caption{(a) Time averaged relative error in region $D_1$. (b) Time averaged absolute error in region $D_2$ over one period.}
	\label{timeerror}
\end{figure}

\section{Conclusions and Future Work} 

\label{future}

In this paper we presented an improved numerical strategy for the control scheme described in \cite{DO_EP} and employed it to perform  several sensitivity studies for the scheme. 

We have shown that for problem \eqref{P1a}, \eqref{P1b}, in the configuration (\ref{geom}-ii), it is possible to characterize sources with weak far field and controllable fields in a compact near field region $D_1$. More explicitly, we should that surface sources could be characterized so that their radiated field approximate an an {\textit {incoming}} plane wave in a region of their near fields while maintaining a very small field beyond a finite radius. The sensitivity studies performed showed that the control scheme produces accurate results with respect to variations in the outward shift, outer radius and in both the outer and inner radii of the near control. Also, our numerics seem to suggest that it may be in principle possible to characterize surface sources with weak far fields so that they approximate prescribed Helmholtz potentials in larger subregions of their near field. 

The control algorithm for problem \eqref{P1a}, \eqref{P1b} was also discussed in the configuration (\ref{geom}-ii). We considered first the situation when an outgoing plane wave was approximated in the near field region while maintaining a null in a prescribed region further away. The sensitivity analysis studied the scheme with respect to variations in the relative position of the control regions and showed very good accuracy results. For the same configuration we also studied the case when a null was to be created in a near-field region while approximating an outgoing plane wave in a prescribed region further away and showed a food performance of our scheme. These results together with the superposition principle suggest that, in principle, our scheme could characterize surface sources so that their radiated fields will approximate distinct prescribed Helmholtz potentials in given mutually disjoint exterior regions of space. 

 In an effort to test our results in the time domain we considered again problem \eqref{P1a}, \eqref{P1b} in the configuration (\ref{geom}-ii) and employed the superposition principle to obtain a Fourier synthesis matching an outgoing plane wave in the primary near field control region in front of the source while keeping a low signature on a region further away behind the source. Our positive result is an indication that it may be possible to characterize such surface controls for time domain signals.
 
 We believe that the study of the case when the source is an array, i.e.,  modeled as a union of compact mutually disjoint subdomains, is very important and our preliminary results in this regard are encouraging. We plan to present a detailed analysis of this case in a forthcoming paper. 


An important feature of our work is mentioned in Remark \ref{rem1} and Remark \ref{rem2}. Indeed, a first observation in this regard is that, as shown in Figure \ref{kw_anr}, Figure \ref{kw_double} and Figure \ref{obs-pattern}, one could employ suitable optimization strategies for the description of the most feasible physical source boundary and this will be part of our related next research efforts. Another important point is that, as mentioned in the introduction, although our research characterizes continuous surface sources for their instantiation one could use monopole and dipole approximating arrays in the spirit of \cite{Doicu}.

In the same spirit, a second important observation is that our analysis is not bounded to spherical fictitious regions $D_{a'}$ as is presented in this paper. In fact one could in principle chose to start with an arbitrary smooth region $D_{a'}$ by adapting the numerical analysis used to evaluate the propagator ${\mathcal D}$ and set up the optimization procedure. This will also be considered as part of our future investigations.

Last but not least, we mention that the sensitivity analysis performed in this paper together with the theory developed previously by our group in \cite{DO_EP} and \cite{DOactivemanipulation} are the foundation of our current efforts to extend the applicability of our scheme to the case of exterior surface control of vector fields which are solutions of the Maxwell equations. This part of our research is the current focus of our group and will be reported very soon.

	\section*{Acknowledgment}
	This work was initially supported by the Office of Naval Research 
	under the award N00014-15-1-2462 and then supported by the Army Research Office 
	under the award W911NF-17-1-0478. 

\bibliographystyle{elsarticle-num}
\bibliography{acousticpaper}

\ifx \bblindex \undefined \def \bblindex #1{} \fi
\begin{thebibliography}{10}
\expandafter\ifx\csname url\endcsname\relax
  \def\url#1{\texttt{#1}}\fi
\expandafter\ifx\csname urlprefix\endcsname\relax\def\urlprefix{URL }\fi
\expandafter\ifx\csname href\endcsname\relax
  \def\href#1#2{#2} \def\path#1{#1}\fi

\bibitem{DOactivemanipulation}
D.~Onofrei, Active manipulation of fields modeled by the helmholtz equation,
  Journal Of Integral Equations and Applications 26~(4) (2014) 553--579.

\bibitem{DO_EP}
D.~Onofrei, E.~Platt, On the synthesis of acoustic sources with controllable
  near fields, Wave Motion 77 (2018) 12--27.

\bibitem{fuller}
C.~R. Fuller, A.~H. von Flotow, Active control of sound and vibration, IEEE
  Control Systems Magazine 15~(6) (1995) 9--19.
\newblock \href {https://doi.org/10.1109/37.476383}
  {\path{doi:10.1109/37.476383}}.

\bibitem{peake}
N.~Peake, D.~G. Crighton, Active control of sound, Annu. Rev. Fluid Mech 32
  (2000) 137--164.

\bibitem{newactr5}
S.~Elliott, Signal Processing for Active Control, 1st Edition, Signal
  Processing and its Applications, Academic Press, 2000.

\bibitem{nelson}
P.~A. Nelson, S.~J. Elliott, Active Control of Sound, Academic London, 1992.

\bibitem{olson}
H.~Olson, E.~May, Electronic sound absorber, J. Acad. Soc. Am. 25 (1953)
  1130--1136.

\bibitem{ANC1}
J.~Cheer, S.~J. Elliott, Multichannel control systems for the attenuation of
  interior road noise in vehicles, Mech. Syst. Signal Process. 60-61 (2015)
  753–769.

\bibitem{kim}
Y.~Kim, J.~Choi, Sound Field Reproduction, Wiley-Blackwell, 2013, Ch.~6, pp.
  283--370.
\newblock \href {https://doi.org/10.1002/9781118368480.ch6}
  {\path{doi:10.1002/9781118368480.ch6}}.

\bibitem{rabenstein}
R.~Rabenstein, S.~Spors, Sound Field Reproduction, In: Benesty J., Sondhi M.M.,
  Huang Y.A. (eds) Springer Handbook of Speech Processing, Springer-Verlag,
  2008.

\bibitem{ahrens1}
J.~Ahrens, S.~Spors, Sound field reproduction using planar and linear arrays of
  loudspeakers, IEEE Transactions on Audio, Speech, and Language Processing
  18~(8) (2010) 2038--2050.

\bibitem{WFS1}
J.~Ahrens, Analytic Methods of Sound Field Synthesis, T-Labs Series in
  Telecommunications Services, Springer, 2012.

\bibitem{miller}
D.~A.~B. Miller, On perfect cloaking, Opt. Express 14 (2006) 12457--12466.

\bibitem{AC1}
Y.~Bobrovnitskii, A new impedance-based approach to analysis and control of
  sound scattering, J. Sound Vib. 297 (2006) 743--760.

\bibitem{AC2}
Y.~Bobrovnitskii, Impedance acoustic cloaking, New J. Phys 12 (2010) 043049.

\bibitem{AC3}
E.~Friot, R.~Guillermin, M.~Winninger, Active control of scattered acoustic
  radiation: A real-time implementation for a three-dimensional object, Acta
  Acust. Acust. 92 (2006) 278--288.

\bibitem{AC4}
A.~N. Noris, Acoustic cloaking, Acoustics Today 11~(1) (2015) 38--46.

\bibitem{broadband_vasquez}
F.~G. Vasquez, G.~W. Milton, D.~Onofrei, Broadband exterior cloaking, Optics
  Express 17~(17) (2009) 14800--14805.

\bibitem{active_vasquez}
F.~G. Vasquez, G.~W. Milton, D.~Onofrei, Active exterior cloaking, Phys. Rev.
  Lett. 103 (2009) 073901.

\bibitem{cheer}
J.~Cheer, Active control of scattered acoustic fields: cancellation,
  reproduction and cloaking, J. Acoust. Soc. Am. 140~(3) (2016) 1502--1512.

\bibitem{vasquez}
F.~G. Vasquez, G.~W. Milton, D.~Onofrei, Exterior cloaking with active sources
  in two dimensional acoustics, Wave Motion 48 (2011) 515--524.

\bibitem{loncaric2}
J.~Lon\v{c}ari\'c, V.~S. Ryaben'kii, S.~V. Tsynkov, Active shielding and
  control of noise, SIAM J. Appl. Math. 62~(2) (2001) 563--596.

\bibitem{loncaric}
J.~Lon\v{c}ari\'c, S.~V. Tsynkov, Quadratic optimization in the problems of
  active control of sound, Applied Numerical Mathematics 52 (2005) 381--400.

\bibitem{olivieri}
F.~Olivieri, F.~Fazi, S.~Fontana, D.~Menzies, P.~Nelson, Generation of private
  sound with a circular loudspeaker array and the weighted pressure matching
  method, IEEE/ACM Transactions on Audio, Speech, and Language Processing PP
  (2017) 1--21.

\bibitem{coleman}
P.~Coleman, P.~J.~B. Jackson, M.~Olik, M.~Moller, M.~Olsen, J.~A. Pedersen,
  Acoustic contrast, planarity and robustness of sound zone methods using a
  circular loudspeaker array, J. Acoust. Soc. Am. 135~(4) (2014) 1929--1940.

\bibitem{cai}
Y.~Cai, M.~Wu, J.~Yang, Sound reproduction in personal audio systems using the
  least-squares approach with acoustic contrast control constraint, J. Acoust.
  Soc. Am. 135~(2) (2014) 734--741.

\bibitem{zhang}
W.~Zhang, T.~D. Abhayapala, T.~Betlehem, F.~M. Fazi, Analysis and control of
  multi-zone sound field reproduction using modal-domain approach, J. Acoust.
  Soc. Am. 140~(3) (2016) 2134--2144.

\bibitem{Cheerps}
S.~J. Elliot, J.~Cheer, J.-W. Choi, Y.~Kim, Robustness and regularization of
  personal audio systems, IEEE Trans. of Audio, Speech, and Language process.
  20~(7) (2012) 2123--2133.

\bibitem{Poletti}
M.~A. Poletti, F.~M. Fazi, An approach to generating two zones of silence with
  application to personal sound systems, The Journal of the Acoust. Soc. of
  America 137 (2015) 598--605.

\bibitem{Onofrei-S}
D.~Onofrei, On the active manipulation of fields and applications. i - the
  quasistatic regime, Inverse problems 28~(10) (2012) 105009.

\bibitem{Doicu}
A.~Doicu, Y.~Eremin, T.~Wriedt, Acoustic and Electromagnetic Scattering
  Analysis Using Discrete Sources, Academic Press, 2000.

\bibitem{hubenthal_DO}
M.~Hubenthal, D.~Onofrei, Sensitivity analysis for active control of the
  helmholtz equation, Applied Numerical Mathematics 106 (2016) 1--23.

\bibitem{Colton-Kress}
D.~Colton, R.~Kress, Inverse Acoustic and Electromagnetic Scattering Theory,
  3rd Edition, Springer-Verlag, 2013.

\bibitem{stegun}
M.~Abramowitz, I.~A. Stegun, Handbook of Mathematical Functions with Formulas,
  Graphs and Mathematical Tables, 10th Edition, Applied Mathematics Series 55,
  National Bureau of Standards, 1972.

\bibitem{watson}
G.~N. Watson, A Treatise on the Theory of Bessel Functions, Cambridge, 1906.

\bibitem{animationfull}
{\url{https://drive.google.com/open?id=1m5nnofuIc56_aaU76HHjqkH4aRLRIsRi}}.

\end{thebibliography}

\end{document}